\documentclass[english,12pt,leqno]{article}
\usepackage{tikz}
\usetikzlibrary{decorations.markings}
\usetikzlibrary{calc}
\usepackage[latin1]{inputenc}
\usepackage{amsxtra}
\usepackage{amsmath}
\usepackage{amssymb}
\usepackage{amsfonts}
\usepackage[all]{xy}
\usepackage{mathrsfs}
\usepackage{amsthm}

\usepackage{babel}

\newtheorem{thm}[equation]{Theorem}

\newtheorem{cor}[equation]{Corollary}
\newtheorem{lem}[equation]{Lemma}
\newtheorem{prop}[equation]{Proposition}

\newtheoremstyle{example}{\topsep}{\topsep}%
     {}
     {}
     {\bfseries}
     {.}
     {2pt}
     {\thmname{#1}\thmnumber{ #2}\thmnote{ #3}}

   \theoremstyle{example}
   
   \newtheorem{Defi}[equation]{Definition}
   
   \newtheorem{rem}[equation]{Remark}

   \newtheorem{ex}[equation]{Example}
   
    \newtheorem{refo}[equation]{Reformulation}

\newtheoremstyle{example}{\topsep}{\topsep}%
     {}
     {}
     {\bfseries}
     {.}
     {2pt}
     {\thmname{#1}\thmnumber{ #2}\thmnote{ #3}}

   \numberwithin{equation}{section}

\def\CC{\mathbb{C}}
\def\DD{\mathbb{D}}

\def\LL{\mathbb{L}}

\def\RR{\mathbb{R}}
\def\ZZ{\mathbb{Z}}

\def\HH{\mathbb{H}}

\def\Fen{\mathfrak{F}}

\def\hen{\mathfrak{h}}

\def\Sen{\mathfrak{S}}

 \def\Gen{\mathfrak{G}}

\def\Ac{\mathcal{A}}
\def\Bc{\mathcal{B}}
\def\Cc{\mathcal{C}}
\def\Kc{\mathcal{K}}
\def\Dc{\mathcal{D}}
\def\Ec{\mathcal{E}}
\def\Fc{\mathcal{F}}
\def\Gc{\mathcal{G}}
\def\Ic{\mathcal{I}}
\def\Jc{\mathcal{J}}
\def\Lc{\mathcal{L}}
\def\Mc{\mathcal{M}}
\def\Nc{\mathcal{N}}
\def\Hc{\mathcal{H}}
\def\Oc{\mathcal{O}}
\def\Pc{\mathcal{P}}
\def\Qc{\mathcal{Q}}
\def\Rc{\mathcal{R}}
\def\Sc{\mathcal{S}}

\def\Xc{\mathcal{X}}

\def\be{\begin{equation}}

\def\cod{\on{codim}}

\def\dim{{\rm{dim}}}

\def\DSh{D^b\Sh }

\def\ee{\end{equation}}

\def\Hom{\on{Hom}}

\def\k {\mathbf k}

\def\lra{\longrightarrow}

\def\on{\operatorname}

\def\orr{\on{or}}

\def\OR{\on{or}}

\def\Perv{\on{Perv}}

\def\Rep{\on{Rep}}

\def\Sh{\on{Sh}}

\def\sgn{\on{sgn}}

\def\ul{\underline}

\def\Vfd{\on{Vect}^{\on{fd}}_\k}

\title{ Perverse sheaves over real hyperplane arrangements}

\author{{Mikhail Kapranov} \and {Vadim Schechtman}}

\begin{document}

\maketitle

\begin{abstract}
Let $\Hc$  be an arrangement of real hyperplanes in $\RR^n$. The complexification of $\Hc$ defines a natural stratification of 
$\CC^n$. We denote by $\Perv(\CC^n, \Hc)$  the category of perverse sheaves on $\CC^n$
 smooth with respect to this stratification. We give a description of $\Perv(\CC^n, \Hc)$  as the category of representations of an explicit quiver with relations, whose vertices correspond to real faces of $\Hc$ (of all dimensions). The relations are of monomial nature: they identify some pairs of paths in the quiver. They can be formulated in terms of the oriented matroid associated to $\Hc$.
\end{abstract}

\tableofcontents

 \addtocounter{section}{-1}

\section{Introduction.}

 \noindent {\bf (0.1)}  Let $X$ be a smooth complex algebraic variety and $\Xc = (X_\alpha)$ be a complex algebraic equidimensional 
  Whitney stratification
  of $X$. Fix a base field $\k$. One has then an abelian category $\Perv(X,\Xc)$ of $\Xc$-smooth perverse sheaves
  of $\k$-vector spaces on $X$, see \cite{BBD, KS}. Understanding this category is one of the central problems of
  topology of algebraic varieties. In many applications it is important to have an explicit description of $\Perv(X,\Xc)$
  as the category of $\k$-representations of some quiver with relations. General microlocal methods 
  \cite{beil-gluing, MV, gelfand-MV} provide a way to approach such a description in principle: an object
  $\Fc\in\Perv(X,\Xc)$ gives rise, for each $\alpha$, to a ``local system of Morse data"
  $\Lc_\alpha$ on some open part of the conormal bundle $T^*_{X_\alpha} X$, and the strategy is  then
  to try to glue the categories formed by the $\Lc_\alpha$ into a single abelian category by an inductive
  procedure, adding one stratum at a time. However, there are only a very few higher-dimensional examples
  where a complete description has been obtained: normal crossings \cite{galligo-GM}, Grassmannians with
  the Schubert stratification \cite{braden-grass}, rectangular matrices with the rank stratification
  \cite{braden-rank}.
  
  \vskip .3cm
  
 \noindent {\bf (0.2)}  In this paper we consider the case when $X=\CC^n$ and $\Xc$ is given by an arrangement $\Hc_\CC$ of
  linear hyperplanes in $X$ with real equations (so $\Hc_\CC$ is the complexification of an arrangement $\Hc$
  of hyperplanes in $\RR^n$). We denote the corresponding category $\Perv(\CC^n,  \Hc)$
  and give a complete, combinatorial description of it in terms of  the following data:
  \begin{itemize}
  \item[(1)] The  poset (partially ordered set) $(\Cc, \leq)$ of {\em faces}, i.e.,
  convex locally closed subsets of all dimensions into which $\Hc$ stratifies $\RR^n$.
  
  \item[(2)] The concept of {\em collinearity}. We call an ordered triple of three faces $(A,B,C)$ collinear if  there are points
  $a\in A, b\in B, c\in C$ such that $b$ lies in the straight line segment $[a,c]$. 
  
  Here the order is important. For example, a triple $(A,A,B)$ is 
  always collinear, whereas $(A,B,A)$ is not if $B\neq A$.   
  
  \end{itemize}
  
\noindent Note that these data can be recovered from the {\em oriented matroid} associated to $\Hc$,
  see \cite{BLSWZ} and Proposition 7.4 below. See also (0.5) below. 
  
  Let us denote by $\Rep^{(2)}(\Cc)$ the category formed by {\em double representations}
  of $\Cc$, i.e., of diagrams consisting of finite-dimensional $\k$-vector spaces $E_C, C\in\Cc$ and linear
  operators 
  \[
  \gamma_{C'C}: E_{C'}\lra E_C, \,\,\,\delta_{CC'}: E_C\lra E_{C'}, \,\,\, C'\leq C,
  \]
  such that the $\gamma_{C'C}$ form a representation of $(\Cc, \leq)$, and the $\delta_{CC'}$ form
  a representation of the opposite poset $(\Cc, \geq)$,  in $\k$-vector spaces.  Our main result,
  Theorem \ref{thm:main}, is that $\Perv(\CC^n, \Hc)$ is equivalent to the full subcategory in $\Rep^{(2)}(\Cc)$ 
  consisting of representations satisfying the following three conditions:
  
  \vskip .2cm
  
  \noindent {\bf Monotonicity:} {\em For any $C'\leq C$ we have $\gamma_{C'C}\delta_{CC'} =\on{Id}_{E_C}$.}
  
  \vskip .2cm
   
  This allows us to define, for any $A,B\in\Cc$, the  transition map
  \[
  \phi_{AB}=\gamma_{CB}\delta_{AC}: E_A\lra E_B,
  \]
  where $C$ is any cell $\leq A,B$. For example, $\phi_{AA} = \text{Id}_{E_A}$. 
  
  \vskip .2cm
  
  \noindent {\bf Transitivity:} {\em If $(A,B,C)$ is a collinear triple of faces, then $\phi_{AC}=\phi_{BC}\phi_{AB}$.}
  
  \vskip .2cm
  
  Note that this relation does {\it not} imply that 
  $\phi_{AB}\phi_{BA} = \on{Id_{E_B}}$.
  
  \vskip .2cm
  
  \noindent {\bf Invertibility:} {\em If $C_1, C_2$ are faces of the same dimension $d$, lying in the same
  $d$-dimensional subspace, on the opposite sides of a $(d-1)$-dimensional face $D$, then $\phi_{C_1C_2}$
  is an isomorphism.}
  
  \vskip .2cm
  
  Let $L\subset X$ be {\em a flat}, i.e. an intersection of some 
  hyperplanes from $\Hc_\CC$, of dimension $d$.  
  The last two conditions mean that the spaces $E_C$ and the maps 
  $\phi_{CC'}$ where $C, C'$ run through all $d$-dimensional 
  faces inside $L$, form a local system over $L^\circ = 
  L\setminus \cup_{H\in \Hc_\CC} H$.   This follows from a 
  description of the fundamental groupoid of $L^\circ$ given in
   Proposition \ref{prop:GK-groupoid}.

  \vskip .3cm
  
  \noindent {\bf (0.3)} The type of description of $\on{Perv}(\CC^n, \Hc)$ that we obtain, is quite different from 
  those appearing in most of the earlier approaches.  More precisely:
    \vskip .2cm
  
  \noindent {\bf (a)} It is of {\em monomial nature}: the conditions on the maps $\gamma_{C'C}, \delta_{CC'}$
  do not appeal to the operations of addition or multiplication by scalars, but only to composition of maps. 
  For comparison, in the most classical case of $X=\CC$ stratified by $\{0\}$ and $\CC\setminus \{0\}$, 
  the standard description \cite{beil-gluing, galligo-GM} is in terms of diagrams
  $
  \xymatrix{
 \{ \Phi \ar@<.7ex>[r]^v& \Psi \ar@<.7ex>[l]^u\}
  }
  $
 such that $\on{Id}_\Psi+vu$ is invertible, so it is not monomial.  (See \S 9 for the comparison of the
 two descriptions in this case). Note that with a monomial description one has the means to define
 what should be a {\em ``perverse sheaf of sets"}, or a {\em ``perverse stack of  categories"}
 on $\CC^n$ smooth with respect to $\Hc_\CC$. Cf. \cite{KS-schobers}. 
 
 \vskip .2cm
 
 \noindent {\bf (b)}  Since the strata $X_\alpha$ are, in our situation, generic parts  of the complex flats $L_\CC$ of $\Hc$,
 see \S \ref{sec:back-arr}D, the local systems $\Lc_\alpha$ of Morse data in the standard approach
 are defined on some open parts of $T^*_{L_\CC}\CC^n$. 
 Our linear algebra data  provide maps ($\phi_{C_1, C_2}$ in the invertibility condition)  which can be related,
 in the sense outlined in (d) below, to {\em half-monodromies} of  appropriate $\Lc_\alpha$ corresponding
  to paths joining neighboring cells and going around the wall in the complex domain. 
  This is different from a more straightforward approach when a local system is described by its
  monodromies corresponding to closed loops.

 It is a known phenomenon in the theory of quantum groups that a typical monodromy matrix $M$ of
 the Knizhnik-Zamolodchikov equation  has quite complicated matrix elements, whereas the  two half-monodromies $M_+$ and $M_-$ of
 which it is composed via $M=M_+ M_-$, are of much
 simpler (monomial) nature, cf. \cite{sv-quant} 
 
 \vskip .2cm
 
 \noindent {\bf (c)} From the purely topological point of view, our approach emphasizes not the
 fundamental groups but the {\em fundamental groupoids} of the complex strata, with as many base points
  as there are different faces in the same stratum. The simpler nature of relations is 
 achieved therefore by introducing a certain redundancy in our description: 
 to each complex stratum we associate not a single vector space (as in the standard approaches)
 but several isomorphic spaces $E_C$ corresponding to 
 different cells $C$  in the  stratum. Our approach
 can  thus be seen as a natural development of the work of Salvetti \cite{salvetti}, 
 Gelfand-Rybnikov \cite{gelfand-rybnikov} and Bj\"orner-Ziegler \cite{BZ} who
 studied the topology of complexified arrangements by real methods and formulated
 the results in terms of oriented matroids. 
 
 \vskip .2cm
 
 \noindent {\bf (d)} The spaces $E_C$ appearing in our description, are not obtained
 as the stalks
 of the local systems $\Lc_\alpha$  at appropriate points (as in most of the standard approaches).
 The indexing sets for the $E_C$ and the $\Lc_\alpha$ are already  quite different. 
 Instead, $E_C$ can be identified, non-canonically,
 with direct sums of several such stalks (involving several different $\Lc_\alpha)$.
  This can be surmised already from the monotonicity
 condition which implies $\dim(E_{C'})\geq\dim(E_C)$ for $C'\leq C$. In particular, the ``biggest"
 space $E_0$ has the dimension equal to the sum of the ranks of all the local systems, corresponding to
 all the strata. In our approach, $E_0$ can be identified with the 
 {\em space of hyperfunction solutions} of the holonomic $\Dc$-module $\Mc$ corresponding
 to the perverse sheaf. More generally, for any face $C$,  the space $E_C$ is the stalk, at any point $c\in C$, of the {\em sheaf
 of hyperfunction solutions} of $\Mc$, i.e., the space of such solutions
 defined in a small neighborhood of $c$ in $\RR^n$.

 The study of the spaces of hyperfunction solutions in 1 dimension (i.e., for $\CC$ stratified by $\{0\}$ and
 $\CC\setminus\{0\}$) goes back to the very origins of
 the theory of $\Dc$-modules as presented in  Kashiwara's 1971 Master Thesis \cite[Th. 4.2.7] {kashiwara},
 and  to the paper of Komatsu  \cite{komatsu} from the same year. These works identify the
 dimension of the
 space $E_0$ in this case with what in the ``standard" (much later) description   
 would be denoted by  
   $\dim(\Phi\oplus\Psi)$.  Higher-dimensional generalizations
 were found in the papers of
 Takeuchi \cite{takeuchi} and Sch\"urmann \cite{schurmann}  of which the first
 considers precisely the situation of the complexification of a real arrangement.
  Spaces which turn out to be identical with our $E_C$, have also appeared in the
  work of Bezrukavnikov, Finkelberg and one of the authors, \cite{BFS}, under the name ``generalized
  vanishing cycles".  In fact, it was conjectured in \cite{BFS}, p. 50, that  a perverse sheaf
  $\Fc$ can be uniquely reconstructed from the linear algebra data equivalent to our
  $(E_C, \gamma_{C'C}, \delta_{CC'})$. From this point of view, we not only prove the
  conjecture of \cite{BFS} but find an explicit characterization of the linear algebra
  data that can appear.   
  \vskip .2cm

 \vskip .3cm

 \noindent {\bf (0.4)} Our method is closest to that of Galligo, Granger and Maisonobe 
 \cite{galligo-GM}
 (who attribute the original idea to Malgrange). That is, we construct a version of Cousin resolution
 of a perverse sheaf $\Fc\in\Perv(\CC^n, \Hc)$ using a stratification of $\CC^n$ into ``tube cells"
 $C+i\RR^n$, $C\in\Cc$,
 so the terms of the resolution are direct sums of the sheaves of cohomology with supports
 in such tubes. It turns out that for each $C$ only  the sheaf
 $\Ec_C= \underline{\HH}^{\on{codim}(C)}_{C+i\RR^n}(\Fc)$ is non-zero,
 and $E_C$ can be identified with the space  of global sections of this sheaf. 
 The main technical point of our study is that the entire Cousin complex $\Ec^\bullet$ formed by the $\Ec_C$
 can be recovered from linear algebra data represented by the $E_C, \gamma_{C'C}$ and $\delta_{CC'}$,
 and that $\Hc_\CC$-smoothness and perversity of $\Ec^\bullet\simeq\Fc$ are precisely equivalent to
 the three conditions above.
 
 \vskip .3cm
 
 \noindent {\bf (0.5)} 
 The concept of collinearity of a triple of faces contains, as a particular case, the familiar
 condition
 \begin{equation*}
 l(w' w'') \,\,=\,\, l(w')+ l(w''), \quad w', w''\in W.  \leqno (*)
 \end{equation*}
 Here $W$ is the Weyl group of a root system $(\hen, \Delta)$. In this case we have the arrangement 
 $\Hc=\{\alpha^\perp\}_{\alpha\in\Delta}$ of root hyperplanes in $\hen$. Chambers (open faces) of $\Hc$
 form a $W$-torsor, so for any two chambers $A,B$ there is a unique $w_{AB}\in W$ such that $w_{AB}\cdot A=B$. 
 A triple of chambers $(A,B,C)$ is collinear if and only if $w'=w_{BC}$ and $w''=w_{AB}$ satisfy (*). 
 The transitivity property $\phi_{AC}=\phi_{BC}\phi_{AB}$ of our double representations is thus reminiscent of
 the classical Gindikin - Karpelevich 
factorization formula, cf. \cite[Th.1]{gk}, and of the cocycle property of the principal series intertwiners
 \cite{knapp, schiffmann, mw} (which is another manifestation of that formula).

\vskip .3cm  
 
 \noindent {\bf (0.6)} We are grateful to Misha Finkelberg 
 for useful discussions and to Persi Diaconis for 
 showing us some interesting references.  We would like to thank Pierre Schapira for
 useful correspondence at the early stage of this work. 
 We are particularly grateful to the referees for many remarks which helped
 us improve the paper. These remarks included suggestions for simplifying several arguments
 as well as pointing out some erroneous ones,  which we have corrected. 
 M.K. would like to thank Universit\'e Paul
 Sabatier for hospitality and financial support during a visit when a substantial part of
 this work was carried out.  His  work was  also 
 supported by World Premier International Research Center Initiative (WPI Initiative), MEXT, Japan.

 \section{Generalities. }
 
\noindent {\bf A. Postnikov systems.}
 Let $\Dc$ be a triangulated category. 
 A {\em left Postnikov system} in $\Dc$ is a diagram 
 of exact triangles of the form 
 \[
 \begin{tikzpicture}[scale=0.4]

 \node (A) at (0,0) {$A=A_{01\cdots n}$}; 
 \node (A1etcn) at (5,0) {$A_{1\cdots n}$}; 
 \node (A2etcn) at (10,0) {$A_{2\cdots n}$}; 
 \node (Adots) at (15,0) {$\cdots$}; 
 \node (An-1n) at (20,0) {$A_{n-1,n}$};
 \node (An) at (25,0) {$A_n$};
 \node (0) at (29,0) {$0$};
 
 \draw[->] (A1etcn) -- (A); 
 \draw[->](A2etcn) -- (A1etcn); 
 \draw[->] (Adots) -- (A2etcn); 
 \draw [->] (An-1n) -- (Adots); 
 \draw [->] (An) -- (An-1n); 
 \draw [->] (0) -- (An); 
 
 \node (A0) at (2.5,-5) {$A_0$}; 
 \node (A1) at (7.5,-5) {$A_1$}; 
 \node (ddots) at (13.5, -5) {$\cdots$}; 
 \node (ddots2) at (17.5,-5) {$\cdots$}; 
 \node (An-1) at (22.5,-5) {$A_{n-1}$}; 
 \node (An') at (27.5,-5) {$A_n.$}; 
 
 \draw[->] (A) -- (A0); 
 \draw[->] (A0) -- (A1etcn); 
 \draw[->] (A1etcn) -- (A1); 
 \draw[->] (A1) -- (A2etcn); 
  \draw[->] (A2etcn) -- (ddots);  
  \draw[->] (ddots2) -- (An-1n); 
    \draw[->]  (An-1n) -- (An-1); 
\draw[->] (An-1n) -- (An);
\draw[->] (An) -- (An');     
\draw[->] (An') -- (0); 
\draw [->]  (An-1) -- (An); 

\node at (0.6,-3) {$\alpha_0$}; 
\node at (4.4, -2.5) {$\beta_0$}; 
\node at (2.9,-2.5) {$+1$}; 
\node at (5.8, -3) {$\alpha_1$}; 
\node at (8,-2.5) {$+1$}; 
\node at (9.4, -2.5) {$\beta_1$}; 
\node at (11.1, -3) {$\alpha_2$}; 
\node at (17.4,-2.5) {$+1$}; 

\node at (20.4,-3) {$\alpha_{n-1}$}; 
\node at (22.8,-2.5) {$+1$};
\node at (24.8,-2.5) {$\beta_{n-1}$};
\node at (25.9,-3.5) {$\alpha_n$};

 \end{tikzpicture}
 \]
 It gives rise to a sequence of objects and morphisms in $\Dc$
 \be\label{eq:postcomleft}
 A_0 \buildrel \delta_0\over\lra A_1[1] \buildrel \delta_1\over\lra A_2[2] \buildrel \delta_2\over
 \lra \cdots \buildrel \delta_{n-1}\over\lra A_n[n], \quad \delta_i = \alpha_{i+1}\beta_i.
  \ee

 Similarly, a {\em right Postnikov system} in $\Dc$ is a diagram of exact triangles of the form
 \[
 \begin{tikzpicture}[scale=0.4]

 \node (A) at (-2,0) {$A=A_{01\cdots n}$}; 
 \node (A0etcn-1) at (4,0) {$A_{0\cdots n-1}$}; 
 \node (A0etcn-2) at (11,0) {$A_{0\cdots n-2}$}; 
 \node (Adots) at (16,0) {$\cdots$}; 
 \node (A01) at (21,0) {$A_{01}$};
 \node (A0) at (26,0) {$A_0$};
 \node (0) at (30,0) {$0$};
 
 \draw[->] (A) -- (A0etcn-1); 
 \draw[->](A0etcn-1) -- (A0etcn-2); 
 \draw[->] (A0etcn-2) -- (Adots); 
 \draw [->] (Adots) -- (A01); 
 \draw [->] (A01) -- (A0); 
 \draw [->] (A0) -- (0); 
 
 \node (An) at (2.5,-5) {$A_n$}; 
 \node (An-1) at (7.5,-5) {$A_{n-1}$}; 
 \node (ddots) at (14.5, -5) {$\cdots$}; 
 \node (ddots2) at (17.5,-5) {$\cdots$}; 
 \node (A1) at (22.5,-5) {$A_{1}$}; 
 \node (A0') at (27.5,-5) {$A_0.$}; 
 
 \draw[->] (An) -- (A); 
 \draw[->] (A0etcn-1) -- (An); 
 \draw[->] (An-1) -- (A0etcn-1); 
 \draw[->] (A0etcn-2) -- (An-1); 
  \draw[->] (ddots) -- (A0etcn-2);  
  \draw[->] (A01) -- (ddots2); 
    \draw[->]  (A1) -- (A01); 
\draw[->] (A0) -- (A1);
\draw[->] (A0') -- (A0);     
\draw[->] (0) -- (A0'); 
 
\node at (0.4,-3) {$\alpha'_n$}; 
\node at (4, -2.5) {$\beta'_n$}; 
\node at (2.6,-2.5) {$+1$}; 
\node at (5.8, -3) {$\alpha'_{n-1}$}; 
\node at (8,-2.5) {$+1$}; 
\node at (10.4, -2.5) {$\beta'_{n-1}$}; 
\node at (12.3, -3.5) {$\alpha'_{n-2}$}; 
\node at (18,-2.5) {$+1$}; 

\node at (21,-3) {$\alpha'_{1}$}; 
\node at (22.8,-2.5) {$+1$};
\node at (24.8,-2.5) {$\beta'_{1}$};
\node at (26.3,-3.5) {$\alpha'_0$};

 \end{tikzpicture}
 \]
 
 It gives rise to a sequence of objects and morphisms in $\Dc$
 \be\label{eq:postcomright}
 A_0 \buildrel \delta'_0\over\lra A_1[1] \buildrel \delta'_1\over\lra A_2[2] \buildrel \delta'_2\over
 \lra \cdots \buildrel \delta'_{n-1}\over\lra A_n[n], \quad \delta'_i = \beta'_{i+1}\alpha'_i.
  \ee

 \begin{prop}\label{prop:post-sys-compl} 
 Suppose we have a left resp. right, Postnikov system in $\Dc$. Then:
 
 (a)  The
 sequence \eqref{eq:postcomleft} resp. \eqref{eq:postcomright} is a complex in $\Dc$,
 i.e., $\delta_{i+1}\delta_i=0$, resp. $\delta'_{i+1}\delta'_i=0$. 
 
 (b) Suppose, in addition, that $\Dc=D^b(\Ac)$ is the bounded derived category of an abelian category
 $\Ac$ and assume that each $A_i[i]$ is quasi-isomorphic to an object $B_i\in \Ac$, so $B_i = H^i(A_i)$
 is the only cohomology object of $A_i$. Then the complex in $\Ac$
 \[
 B_0\lra B_1 \lra\cdots \lra B_n
 \]
 induced from  \eqref{eq:postcomleft} resp. \eqref{eq:postcomright} by passing to $H^0$, is
 an object of $D^b(\Ac)=\Dc$,  isomorphic to $A$. 
  \end{prop}
  
  \noindent {\sl Proof:} (a) follows because the composition of two consecutive morphisms in an
  exact triangle is equal to 0. Part (b) is proved by induction on the length of the Postnikov system. 
  \qed
  
  \begin {rem} In general, the left or right Postnikov system as above
   exhibits $A$ as a total object of the complex
   \eqref{eq:postcomleft} or \eqref{eq:postcomright}, see
  \cite{gelfand-manin}, Ch. 4, \S 2. 
  \end{rem}
  
  \vskip .3cm
  
  \noindent {\bf B. Filtered topological spaces.} 
By a {\em space} we  mean a 
 topological space homeomorphic to an open subset of a finite
 CW-complex. We fix a base field $\k$. For a space $X$ we denote by $\Sh_X$
 the category of sheaves of $\k$-vector spaces on $X$ and by $D^b\Sh_X$ the corresponding
 bounded derived category. For any map $f: X\to Y$ of spaces we have the standard functors
 $f^*, f^!: \DSh_Y\to\DSh_X$ and $f_*, f_!: \DSh_X\to\DSh_Y$.  In particular, we reserve
 the notation $f_*$ for the derived direct image functor, and denote the usual direct
 image functor on sheaves by $R^0f_*$. 
 If $j: Y\hookrightarrow X$ is an embedding of a locally closed subspace, then we denote
 \[
 \underline{R\Gamma}_Y(\Fc) = j_* j^! \Fc
 \]
 the complex of ``cohomology with support in $Y$".

 \vskip .2cm
 
 If $V$ is a $\k$-vector space, then we denote by $\underline V_X$ the constant sheaf on $X$
 with stalk $V$. If $i: Z\hookrightarrow X$ is the embedding of
  a closed subspace, then, by a slight abuse of notation,
 we consider $\underline V_Z$ as a sheaf on $X$ via the direct image functor $i_*$. 
  
 \vskip .2cm

  Let $X$ be a space and
 \[
 \Xc = \bigl\{ X_0 \subset X_1 \subset \cdots \subset X_n=X\bigr\} 
 \]
 be a filtration of $X$ by closed subspaces. We then have the locally closed subspaces
 $Y_d = X_d\setminus X_{d - 1}$ and denote
 \[
 i_d: X_d \hookrightarrow X, \quad j_d: Y_d \hookrightarrow X
 \]
 the embeddings. 
 
 Any complex of sheaves $\Fc\in D^b\Sh_X$  
 includes into two canonical Postnikov systems which we call the {\em cohomological}  (right) and
 {\em homological} (left) Postnikov systems of $\Fc$ relative to the filtration $\Xc$. The
 cohomological  system has the form
 
 \[
 \begin{tikzpicture}[scale=0.4]

 \node (A) at (-3,0) {$\Fc= i_{n!} i_n^*\Fc$}; 
 \node (A0etcn-1) at (4,0) {$(i_{n-1})_! i_{n-1}^*\Fc$}; 
 \node (A0etcn-2) at (11,0) {$(i_{n-2})_! i_{n-2}^*\Fc$}; 
 \node (Adots) at (16,0) {$\cdots$}; 
 \node (A01) at (21,0) {$i_{1!} i_1^*\Fc$};
 \node (A0) at (26,0) {$i_{0!}i_0^*\Fc$};
 \node (0) at (30,0) {$0$};
 
 \draw[->] (A) -- (A0etcn-1); 
 \draw[->](A0etcn-1) -- (A0etcn-2); 
 \draw[->] (A0etcn-2) -- (Adots); 
 \draw [->] (Adots) -- (A01); 
 \draw [->] (A01) -- (A0); 
 \draw [->] (A0) -- (0); 
 
 \node (An) at (2.5,-5) {$j_{n!} j_n^*\Fc$}; 
 \node (An-1) at (7.5,-5) {$(j_{n-1})_! j_{n-1}^*\Fc$}; 
 \node (ddots) at (14.5, -5) {$\cdots$}; 
 \node (ddots2) at (17.5,-5) {$\cdots$}; 
 \node (A1) at (22.5,-5) {$j_{1!} j_1^*\Fc$}; 
 \node (A0') at (27.5,-5) {$j_{0!}j_0^*\Fc$}; 
 
 \draw[->] (An) -- (A); 
 \draw[->] (A0etcn-1) -- (An); 
 \draw[->] (An-1) -- (A0etcn-1); 
 \draw[->] (A0etcn-2) -- (An-1); 
  \draw[->] (ddots) -- (A0etcn-2);  
  \draw[->] (A01) -- (ddots2); 
    \draw[->]  (A1) -- (A01); 
\draw[->] (A0) -- (A1);
\draw[->] (A0') -- (A0);     
\draw[->] (0) -- (A0'); 
 
 \node at (2.6,-2.5) {$+1$}; 
 
\node at (8,-2.5) {$+1$}; 
 
\node at (18,-2.5) {$+1$}; 

 \node at (22.8,-2.5) {$+1$};

 \end{tikzpicture}
 \]
 and gives rise to the complex in $\DSh_X$
 \be
 j_{0!} j_0^*\Fc \lra j_{1!} j_1^*\Fc[1] \lra\cdots \lra j_{n!} j_n^*\Fc[n]
 \ee
 with total object $\Fc$. The  homological system has the form
  \[
 \begin{tikzpicture}[scale=0.4]

 \node (A) at (-2,0) {$\Fc= i_{n*} i_n^!\Fc$}; 
 \node (A1etcn) at (5,0) {$(i_{n-1})_*i_{n-1}^!\Fc$}; 
 \node (A2etcn) at (12,0) {$(i_{n-2})_*i_{n-2}^!\Fc$}; 
 \node (Adots) at (17,0) {$\cdots$}; 
 \node (An-1n) at (22,0) {$i_{1*}i_1^!\Fc$};
 \node (An) at (27,0) {$i_{0*}i_0^!\Fc$};
 \node (0) at (31,0) {$0$};
 
 \draw[->] (A1etcn) -- (A); 
 \draw[->](A2etcn) -- (A1etcn); 
 \draw[->] (Adots) -- (A2etcn); 
 \draw [->] (An-1n) -- (Adots); 
 \draw [->] (An) -- (An-1n); 
 \draw [->] (0) -- (An); 
 
 \node (A0) at (1.5,-5) {$j_{n*}j_n^!\Fc$}; 
 \node (A1) at (7,-5) {$(j_{n-1})_*j_{n-1}^!\Fc$}; 
 \node (ddots) at (13.5, -5) {$\cdots$}; 
 \node (ddots2) at (17.5,-5) {$\cdots$}; 
 \node (An-1) at (22.5,-5) {$j_{1*}j_1^!\Fc$}; 
 \node (An') at (27.5,-5) {$j_{0*}j_0^!\Fc$}; 
 
 \draw[->] (A) -- (A0); 
 \draw[->] (A0) -- (A1etcn); 
 \draw[->] (A1etcn) -- (A1); 
 \draw[->] (A1) -- (A2etcn); 
  \draw[->] (A2etcn) -- (ddots);  
  \draw[->] (ddots2) -- (An-1n); 
    \draw[->]  (An-1n) -- (An-1); 
\draw[->] (An-1n) -- (An);
\draw[->] (An) -- (An');     
\draw[->] (An') -- (0); 
\draw [->]  (An-1) -- (An);

\node at (2.5,-2.5) {$+1$}; 
 \node at (8.5,-2.5) {$+1$}; 
 
\node at (18.4,-2.5) {$+1$};

\node at (23.8,-2.5) {$+1$};

 \end{tikzpicture}
 \]
 and gives rise to a complex in $\DSh_X$
 \be
 j_{n*}j_n^!\Fc \lra (j_{n-1})_*j_{n-1}^!\Fc[1] \lra \cdots\lra
 j_{0*}j_0^!\Fc[n]
 \ee
 with total object $\Fc$. 
 
 \vskip .3cm
 
 \noindent {\bf C. Verdier duality.} For a space $X$ we denote by
 \[
 \DD: \DSh_X\lra\DSh_X
 \]
 the Verdier duality functor, see \cite{KS, verdier}. 
 We recall  that $\DD$ interchanges the functors  $f_*$ and $f_!$,
as well as the functors $f^*$ and $f^!$ for any continuous map $f$ of spaces. 
 
 Let $X$ be a real analytic manifold of dimension $d$. 
 We denote by $\orr_X$ the orientation local system on $X$. This is the rank 1 local system
 of $\k$-vector spaces whose stalk at $x\in X$ is $H^d_c(U, \k)$, where $U$ is any open
  neighborhood of
 $X$ homeomorphic to a $d$-ball. In this case
 \[
 \DD(\Fc) = \underline{R\Hom}(\Fc, \orr_X[d]),
 \]
 and we denote
 \be
 \Fc^\bigstar = \DD(\Fc)[-d] =  \underline{R\Hom}(\Fc, \orr_X)
 \ee
 the shifted Verdier duality which has the advantage of preserving local systems
 in degree 0. 

 We further  denote by $D^b_{\on{constr}}\Sh_X\subset \DSh_X$ 
the subcategory formed by $\RR$-constructible complexes,
i.e., by complexes whose cohomology sheaves are $\RR$-constructible
\cite{KS}.  The functor $\DD$ preserves the subcategory $D^b_{\on{constr}}\Sh_X$ and 
its restriction there is a perfect duality. The same is true for $\bigstar$. 

If $\Xc$ is a filtration of $X$ by real analytic subsets, then $\DD$ takes the
cohomological Postnikov system for $\Fc$ to the homological Postnikov system for $\DD(\Fc)$
for any $\RR$-constructible complex $\Fc$. 

\vskip .3cm

\noindent {\bf D. Cellular sheaves and complexes.} For background on cellular spaces and sheaves
we refer the reader to \cite{curry, polishchuk,  vybornov1, vybornov2}. Here we present a 
(self-contained) synopsis of features needed in the rest of the paper.

 By a $d$-{\em cell} we mean a topological
space homeomorphic to an open $d$-ball. For a $d$-cell $\sigma$ we denote
\[
\OR(\sigma) = H^d_c(\sigma, \k)
\]
its 1-dimensional orientation space.  For constant sheaves on $\sigma$ the Verdier
duality has the form
\[
\DD(\underline V_\sigma) = \underline{V^*\otimes_\k \OR(\sigma)}_\sigma[d]. 
\] 
 By a {\em cellular space} we will mean a space with a filtration $\Xc$ by closed subspaces
 such that each $X_d= X_d\setminus X_{d-1}$ is a disjoint union of finitely many $d$-cells.  
 We denote $\Cc=\Cc_\Xc$ the set of all cells of $X$ together with the partial order
 \[
 \sigma'\leq\sigma \quad \Leftrightarrow \quad \sigma'\subset\overline\sigma.
 \]
 By a slight abuse of language we will refer to the relation $\leq$ as {\em cell inclusion}, although,
 set-theoretically, it is the closures of the cells that are included.

  For a cell $\sigma$
 of a cellular space $X$ we denote by $j_\sigma:\sigma\to X$ the corresponding embedding. 
 A cellular space $X$ will be called {\em regular}, if the closure of each $d$-cell is
 homeomorphic to a closed $d$-ball. A cellular space $X$ will be called
  {\em quasi-regular}, if it can be represented as $X'-X''$, where $X'$ is a regular cellular space
  and $X''$ is a closed cellular subspace (union of some cells of $X'$) which is then  also regular. 
  In the sequel all cellular spaces will be assumed
quasi-regular.

 For any $\k$-vector space $V$ we have an identification in $\DSh_X$:
 \[
 \underline V_{\overline \sigma} = j_{\sigma *} \underline V_\sigma = R^0j_{\sigma *} \underline V_\sigma.
 \]
 A {\em cellular sheaf} on $X$  (with respect to $\Xc$) 
 is, by definition, a sheaf $\Fc$ on $X$ such that each $j_\sigma^*\Fc$
 is a constant sheaf on $\sigma$ with finite-dimensional stalks.
  We denote $\Sh_{X,\Xc}$ the category of cellular sheaves on $X$
 with respect to $\Xc$, and by $D^b_\Xc(\Sh_X)\subset \DSh_X$ the full subcategory of complexes
whose cohomology sheaves lie in $\Sh_{X,\Xc}$. 

For a cellular sheaf $\Fc$ on $X$ we have the {\em linear algebra data}
 which consists of stalks $F_\sigma = H^0(\sigma, j_\sigma^*\Fc)$ and
 {\em generalization maps} (terminology  taken from \cite{gelfand-macpherson})
 \[
 \gamma_{\sigma'\sigma}: F_{\sigma'}\lra F_\sigma, \quad \sigma'\leq\sigma
 \]
 defined as follows. Take a point $x'\in\sigma'$. Then $F_{\sigma'} = H^0(U',\Fc)$,
 where $U'$ is a sufficiently small contractible neighborhood of $x'$ in $X$. Let $x\in U'\cap\sigma$.
 Then $F_\sigma= H^0(U,\Fc)$, where $U$ is a sufficiently small contractible neighborhood
 of $x$ in $X$. Taking $U$ small enough, we can assume $U\subset U'$. Then
$\gamma_{\sigma'\sigma}$ is given by the restriction map
\[
F_\sigma' = H^0(U', \Fc) \buildrel\on{Res}\over\lra H^0(U, \Fc) = F_\sigma.
\]
The following is by now well known.

\begin{prop}\label{prop:cellullar-LAD}
(a) For $\Fc\in\Sh_{X, \Xc}$ the data $(F_\sigma, \gamma_{\sigma'\sigma})$ form a
representation of the poset $(\Cc, \leq)$, i.e., a covariant functor $\Cc\to\on{Vect}^{\on{fd}}_\k$.

(b) The construction in (a) defines an equivalence of  $\Sh_{X, \Xc}$ with
$\Rep(\Cc)= \on{Fun}(\Cc, \Vfd)$.

(c) The natural functor
\[
\phi: D^b\Rep(\Cc) = D^b \Sh_{X,\Xc} \lra D^b_\Xc(\Sh_X)
\]
is an equivalence of triangulated categories. 

\end{prop}

\noindent {\sl Proof:} (a) and (b) are obvious. Cf. the case of
``simplicial complexes" considered in detail by Kashiwara
\cite{kashiwara-RH}. 
To see (c), note that the sheaf $j_{\sigma *} \underline \k$
corresponds via (b) to the injective object of $\Rep(\Cc)$, i.e., to the 
covariant functor  $R_\sigma: \tau \mapsto (\k \Hom_\Cc (\tau,\sigma))^*$
dual to the representable contravariant one.  
 Therefore for $p\in\ZZ$ 
 \be\label{eq:formula}
 \begin{gathered}
 \Hom_{D^b\Sh_{X,\Xc}} (j_{\sigma *} \underline \k, j_{\tau *} \underline \k[p])\,\,=\,\,
\on{Ext}^p_{\Sh_{X,\Xc}}(j_{\sigma *} \underline \k, j_{\tau *} \underline \k) \,\,=
\\
 \,\,=\,\,
\on{Ext}^p_{\on{Rep}(\Cc)}(R_\sigma. R_\tau) \,\,=\,\,
\begin{cases}
\k, \text{ if } p=0, \,\, \sigma\geq\tau;\\
0, \text {otherwise}. 
\end{cases}
\end{gathered}
\ee
On the other hand, adjointness between $j_\tau^*$ and $j_{\tau *}$ gives that
\[
\Hom_{D^b\Sh_{X}}(j_{\sigma *} \underline \k, j_{\tau *} \underline \k[p] )\,\, =\,\, 
\Hom_{D^b\Sh_{\tau}} (j_\tau^* j_{\sigma *}\underline \k, \underline\k[p])
\]
 and so it is given by the RHS of \eqref{eq:formula} as well. This means that the morphism
\[
\phi_{\Fc, \Gc}: \Hom_{D^b\Sh_{X,\Xc}}(\Fc, \Gc) \lra \Hom_{D^b\Sh_{X}}(\Fc, \Gc)
\]
is an isomorphism whenever both $\Fc$ and $\Gc$ are shifts of sheaves of the form $j_{\sigma_*}\underline\k$,
 $\sigma\in\Cc$. Now, each object of $\Sh_{X,\Xc}$ has a finite resolution by direct sums of the injective objects
 $j_{\sigma *} \underline \k$. This means that
   $\phi_{\Fc, \Gc}$ is an isomorphism for any $\Fc, \Gc\in D^b \Sh_{X,\Xc}$, i.e., that $\phi$ is fully faithful. 
 To see that $\phi$ is essentially surjective, note that the cohomological Postnikov system
of any object  $\Fc\in D^b_\Xc(\Sh_X)$ shows that $\Fc$ lies in the smallest triangulated
category containing all the $j_{\sigma *}\underline\k$. 
\qed

\begin{rem}\label{rem:lin-alg-data}
An object of $D^b\Rep(\Cc)$ corresponding to $\Fc\in D^b_\Xc(\Sh_X)$ by (c)
can be seen as ``linear algebra data for $\Fc$ at the level of complexes", i.e.,
as a choice of actual complexes of vector spaces $F_\sigma^\bullet$ quasi-isomorphic to
$R\Gamma(\sigma, j_\sigma^* \Fc)$ and actual morphisms of complexes 
$\gamma_{\sigma'\sigma}: F^\bullet_{\sigma'}\to F^\bullet_\sigma$
 forming a representation of $(\Cc,\leq)$ at the level of complexes. 
 \end{rem}
 
 \vskip .2cm

 \begin{prop}\label{prop:VD-explicit} 
 Let $X$ be a quasi-regular cellular space of dimension $n$, and $\Fc$
 a cellular complex represented by $(F_\sigma^\bullet, \gamma_{\sigma'\sigma})$.
 Then the the complex $\DD(\Fc[n])$ is quasi-isomorphic to  the total complex of the double
 complex
 \[
 \bigoplus_{\cod(\sigma)=0} \underline{ F_\sigma^*\otimes\OR(\sigma)}_{\overline\sigma}
 \lra 
  \bigoplus_{\cod(\sigma)=1} \underline{ F_\sigma^*\otimes\OR(\sigma)}_{\overline\sigma}
\lra\cdots
\]
Here each summand is a complex of constant  sheaves on the closed cell. \qed
 
 \end{prop}
 
 This is a combinatorial version of the Cousin complex. 
 
 \vskip .3cm
 
 \noindent{\sl Proof:}
 Let first $\Fc$ be a cellular sheaf (not a complex) on $X$, so the  $F_\sigma^\bullet = F_\sigma$ are
 vector spaces. 
  The cohomological Postnikov system for $\Fc$ relative to $\Xc$ gives
 a complex in $D^b\Sh_X$
 \be\label{eq:coh-sys-F}
 \bigoplus_{\dim(\sigma)=0} j_{\sigma !}( \underline {F_\sigma}_\sigma)
 \lra \bigoplus_{\dim(\sigma)=1} j_{\sigma !} (\underline {F_\sigma}_\sigma) [1] \lra \cdots
 \ee
 whose total object is $\Fc$. Applying the exact functor $\DD[-n]$ to \eqref{eq:coh-sys-F}, we get
 a complex in $D^b\Sh_X$ described in the Proposition, with total object $\Fc$. By 
  By Proposition 
  \ref{prop:post-sys-compl}(b), this complex represents $\DD(\Fc[n])$. 
  
  The general case follows  by compatibility of $\DD[-n]$ with 
  forming total complexes of double
 complexes. \qed
 
 \begin{rem}
If $\Fc$ is a cellular sheaf as above and $X$ is regular, then 
 \[
 R\Gamma (j_{\sigma !} \underline{F_{\sigma}}_\sigma) = R\Gamma_c (\sigma, F_\sigma) = F_\sigma \otimes \orr(\sigma) [- \dim(\sigma)]
 \]
 and so from \eqref{eq:coh-sys-F},
  we get a complex  in $D^b\Vfd$
 \[
 \bigoplus_{\dim(\sigma)=0} F_\sigma\otimes\OR(\sigma) \lra 
 \bigoplus_{\dim(\sigma)=1} F_\sigma\otimes\OR(\sigma) \lra \cdots
 \]
 whose total object  is $R\Gamma(X,\Fc)$. By Proposition 
  \ref{prop:post-sys-compl}(b), this complex represents
  $R\Gamma(X,\Fc)$. This is the ``standard cellular
  cochain complex" of $\Fc$.  Compare also with the well known fact  that the geometric realization
of the simplicial nerve of $(\Cc, \leq)$ is, for a regular $X$,  homeomorphic
 to $X$.   
  \end{rem}

\section{Background on real arrangements and the three stratifications}\label{sec:back-arr}

\noindent {\bf A. Faces and sign vectors.}
Let $V$ be a  finite-dimensional real vector space and 
$\Hc$ be an arrangement of linear hyperplanes in $V$. 
Note that, in particular, $0\in H$ for any hyperplane $H\in\Hc$. 
We choose, once and for all, a linear equation $f_H: V\to\RR$ for each $H\in\Hc$
(the essential concepts that we define will not depend on the choice).

 By $\Lc=\Lc_\Hc$ we denote the poset of {\em flats}
of $\Hc$, i.e., of linear subspaces of the form $\bigcap_{H\in \Ic} H$ for various subsets $\Ic\subset\Hc$. 
Note that $\Lc$ contains $V$ (for $\Ic=\emptyset$). We will  assume that  
 $\{0\}\in\Lc$.  This can be always achieved by quotienting by the smallest flat of $\Hc$,
 without changing the combinatorial structure of the arrangement. 

In the sequel we also  assume $V=\RR^n$ for simplicity. 
We denote by
\[
\sgn: \RR \lra \{+, -, 0\}
\]
the standard sign function. By a {\em sign vector} we will mean a sequence $(s_H)_{H\in\Hc} \in\Hc^{\{+,-,0\}}$
assigning a ``sign" to each element of $\Hc$. Each $x\in\RR^n$ gives rise to a sign vector
$(\sgn f_H(x))_{H\in\Hc}$. Level sets of this vector function subdivide $\RR^n$ into 
locally closed subsets called {\em faces}. Thus $x$ and $x'$ lie in the same face, iff
$\sgn f_H(x)=\sgn f_H(x')$ for each $H\in\Hc$. We will sometime identify a face $C$ with the corresponding
sign vector
\[
C \leftrightarrow (C_H)_{H\in\Hc}, \,\,\, C_H = \sgn (f_H|_C) \in \{+,-, 0\}. 
\]
Faces are convex subsets of $\RR^n$, each given by a system of linear equations and strict
linear inequalities \cite{BFS, salvetti, schechtman-varchenko, varchenko}. Open faces will be called
{\em chambers}. We denote by $\Cc = \Cc_\Hc$ the poset of faces, with the partial order given by
$C'\leq C$, if $C'\subset \overline C$. 

Note that the cells form a quasi-regular cell decomposition of $\RR^n$
which we also denote $\Cc$ (taking the one-point  compactification of  $\RR^n$ to a sphere, we
 embed $(\RR^n, \Cc)$ into a regular cellular space).

We also write
\[
C' <_1 C, \text{ if }  C'\leq C \text{ and } \dim(C') = \dim(C)-1
\]
and call this relation {\em codimension 1  inclusion} of faces.

Alternatively, let us introduce a partial order $\leq$ on $\{+,-,0\}$ in which the only nontrivial inequalities are
\be\label{eq:sign-ord}
0\leq +, \,\,\, 0\leq -
\ee
(and $+$ and $-$ are incomparable). This is the order of inclusion of the subsets $\{x\geq 0\}$, $\{x\leq 0\}$ and $\{x=0\}$
in $\RR$. Further introduce on $\Hc^{\{+,-,0\}}$ the Cartesian product partial order. Then
\[
C\leq D \,\,\, \Leftrightarrow \,\,\, C_H\leq D_H \text { for each } H\in\Hc, 
\]
i.e., $\Cc \hookrightarrow \Hc^{\{+,-,0\}}$ is an embedding of posets. 
Faces form a regular cellular stratification of $\RR^n$, which we also denote $\Cc$. Note that $\Cc$
does not depend on the choice of equations $f_H$, both as a poset and as a stratification.
For $x\in\RR^n$ we denote by $\sigma(x)=\sigma_\Hc(x)\in\Cc$ the face containing $x$.
Similarly, if $X\subset\RR^n$ is a subset contained in one face, we denote by $\sigma(X)=\sigma_\Hc(X)$
this face.

\vskip .3cm

\noindent {\bf B. Functoriality of faces.}
Any $L\in\Lc$ gives rise to two induced hyperplane arrangements:
\begin{itemize}
\item
The arrangement $\Hc\cap L$ in $L$, formed by the $H\cap L$, $H\in \Hc$, $H\not\supset L$.

\item The arrangement $\Hc/L $ in $V/L$, formed by the $H/L$, $H\in \Hc$, $H\supset L$. 
\end{itemize}

\noindent We have therefore the face stratifications $\Cc_{\Hc\cap L}$ and $\Cc_{\Hc/L}$
of these arrangements.

 For any face $C\in \Cc$ we denote  by $\LL(C)$ the  $\RR$-linear subspace in $\RR^n$ spanned by
$C$ and denote by 
\[
\pi_C: \RR^n\lra\RR^n/\LL(C)
\]
the projection. Note that $\pi_C$ induces an isomorphism of posets
\be\label{eq:overline-pi-C}
\overline \pi_C: \Cc^{\geq C} \buildrel\simeq \over\lra \Cc_{\Hc/\LL(C)}, \quad  D\mapsto D/\LL(C).
\ee
Here $\Cc^{\geq C}$ is the {\em star} of $\Cc$, i.e.,
the poset of $D\in\Cc = \Cc_\Hc$ such that $D\geq C$. For such a $D$, we denote
 $D/\LL(C) = \pi_C(D)$ (which is a chamber of $\Hc/\LL(C)$). 
 
 \vskip .3cm
 
 \noindent {\bf C. Composition of faces.} Let $(\Sigma, \leq)$ be a poset.
 As usual, we write $a<b$, if $a\leq b$ and $a\neq b$. 
  Following
 \cite{BZ}, introduce a binary operation $\star$ on $\Sigma$, called {\em composition}, by
 \[
 a\star b = \begin{cases}
 b, & \text{if } a<b,\\
 a, & \text{otherwise.}
 \end{cases}
 \]
 This operation is associative but not commutative in general. If $\leq$ is a total order, then $a\circ b=\max(a,b)$
 is commutative. Thus $\star$ can be seen as a generalization of the maximum to partially ordered sets. 
 
 We are interested in the case when $\Sigma=\{+,-, 0\}$ with the partial order 
 \eqref{eq:sign-ord}.

Define an associative operation $\circ$ on the set $\Hc^{\{+,-,0\}}$ of sign vectors
 by putting
 \[
 (s\circ t)_H = s_H \star t_H, \,\,\, H\in\Hc.
 \]
  
 \begin{prop}\label{prop:composition}
 Let $C,D\in\Cc$ be two faces with corresponding sign vectors $(C_H)_{H\in\Hc}$ and $(D_H)_{H\in\Hc}$. Then:
 \begin{itemize}
 
 \item[(a)] The sign vector $(C_H\circ D_H)_{H\in\Hc}$ corresponds to a (necessarily unique) face $C\circ D$. We thus
 obtain an associative binary operation $\circ$ on $\Cc$. 
 
 \item[(b)]  Explicitly, choose any $c\in C, d\in D$. Then
 \[
 C\circ D \,\,=\,\,\sigma\bigl( (1-\epsilon) c + \epsilon d\bigr), \,\,\, 0<\epsilon \ll 1
 \]
 is the cell containing a small displacement of $c$ in the direction of $d$.
 
 \item[(c)]
Alternatively, $C\circ D$ is the minimal cell $K\in\Cc$ such that $K\geq C$ and $K+\LL(C) \supset D$. 
 \end{itemize}
 \end{prop}
 
 \begin{rem}
 The operation $\circ$ on $\Cc$ is the cornerstone of the ``covector" axiomatization of oriented matroids, see
 \cite{BLSWZ}, Axioms 4.1.1. An arrangement of hyperplanes in $\RR^n$ gives  a {\em representable
 oriented matroid}, see {\em loc. cit.} \S 2.1. In this way, isomorphism classes of representable oriented matroids are
 in bijection with combinatorial equivalence classes of real hyperplane arrangements.

  The characterization in (b) is also taken from
 {\em loc. cit.} \S 4.1.

The operation $\circ$ was   introduced (in a slightly different context) by Tits under the name of "projection", cf. \cite{tits}, 3.19. It plays a basic 
role in the studies of random walks on hyperplane arrangements, 
cf. \cite{brown-diacon} and references therein.    
  
 \end{rem}
 
 \noindent {\sl Proof of Proposition \ref{prop:composition}:} Parts (a) and (b) follow from the next lemma, whose
 verification is left to the reader.
 
 \begin{lem}
 Let $x, y\in\RR$. Then, for $0<\epsilon \ll 1$, we have
 \[
 \sgn\bigl( (1-\epsilon) x + \epsilon y\bigr) \,\,\,=\,\,\, \sgn(x) \star \sgn(y).\qed
 \]
 \end{lem}
 
 Part (c) follows from (b). Indeed, (b) implies that $(C\circ D)+\LL(C)\supset D$. 
 Conversely, suppose $K\geq C$ and $K+\LL(C)\supset D$. We claim that $K\geq C\circ D$, that is, $K_H\geq C_H\star D_H$ for
 each $H\in\Hc$. We already know that $K_H\geq C_H$, since $K\geq C$. We need to prove that
 whenever $D_H > C_H$, we also have $K_H\geq D_H$. But $D_H > C_H$ means that $C_H=0$ and $D_H\neq 0$,
 that is, $C\subset H$ but $D\not\subset H$. The statement that $K_H\geq D_H$ means therefore
 $K_H=D_H$, that is, $K$ and $D$ lie on the same side of $H$. But this is clear since $D\subset K+\LL(C)$
 and $\LL(C)\subset H$ since $C\subset H$.\qed
 
 \vskip .2cm
 
 For future reference let us note yet another characterization of $C\circ D$. Consider the image
 $\pi_C(D)\subset \RR^n/\LL(C)$. This image lies in one face but may  not itself  be a face of the quotient arrangement 
 $\Hc/\LL(C)$. We denote $\sigma_{\Hc/\LL(C)} (\pi_C(D))\in\Cc_{\Hc/\LL(C)}$ the cell of $\Hc/\LL(C)$ containing
 $\pi_C(D)$.
 
 \begin{prop}\label{prop:circ-pi}
  \[
 C\circ D \,\,=\,\,\overline\pi_C^{-1} \sigma_{\Hc/\LL(C)} (\pi_C(D)),
 \]
 where $\overline\pi_C$ is the isomorphism of posets from \eqref{eq:overline-pi-C}. 
 \end{prop}
 
 \noindent {\sl Proof:} This is a reformulation of Proposition \ref{prop:composition}(c).\qed
 
 \vskip .2cm
 
 We further note the following monotonicity properties of the composition.
 
 \begin{prop}\label{prop:circ-monotone}

 \begin{itemize}
\item[(a)] If $D'\leq D$, then for any $C$ we have $C\circ D'\leq C\circ D$.

\item[(b)]  For any $C$ and $D$ the flat $\LL(C\circ D)$ is the minimal flat $L\in \Lc_\Hc$ containing
both $\LL(C)$ and $\LL(D)$. In particular, if $C'\leq C$, then for any $D$ we have
$\LL(C'\circ D) \subset \LL(C\circ D)$. 
 
 \end{itemize}
 \end{prop}
 
 \noindent {\sl Proof:} (a) follows from the obvious fact that the operation $a\star b$ on a poset $\Sigma$
 is monotone in the second argument. Part (b) follows because
 \[
 \LL(C\circ D) \,\,\,= \,\,\,\bigcap_{C_H\star D_H=0} H,
 \]
 and $C_H\star D_H=0$ means that $C_H=D_H=0$, i.e., $C,D\subset H$. \qed
 
 \begin{rem}
 Note that $\circ$ is not monotone in the first argument, That is, the condition $C'\leq C$ does not imply that
 $C'\circ D \leq C\circ D$. To obtain a counterexample, it suffices to take $C'=0$. Then $C'\circ D=D$ for any $D$.
 So if $D\not \geq C$, then $0\circ D \not\leq C\circ D$. 
 
 \end{rem}

   \vskip .3cm
   
   \noindent{\bf D. Complexified arrangement and its stratifications.}
   Let $V_\CC=\CC^n$ be the complexification of $V=\RR^n$. 
   For $z\in\CC^n$ we denote by $\Re(z), \Im(z)\in\RR^n$ its real and imaginary parts. 
For $L\in \Lc$ we denote by $L_\CC\subset \CC^n$ its complexification. 
Accordingly, we write $\LL_\CC(X)$ for the $\CC$-linear span of $X\subset\CC^n$. 
 We have therefore the
{\em complexified arrangement} $\Hc_\CC$  of complex hyperplanes in $V_\CC=\CC^n$,
formed by the $H_\CC, H\in\Hc$. We will be interested in three
natural stratifications of $\CC^n$ induced by $\Hc$.

\vskip .2cm

\noindent The {\bf complex stratification }$\Sc^{(0)}$ consists of the open parts of the complexified
flats
 \[
L_\CC^\circ = L_\CC \setminus\bigcup_{H\not\supset L} H_\CC. 
\]

\noindent The {\bf ${\bf s}^{(2)}$-stratification $\Sc^{(2)} = \Cc+i\Cc$} 
consists of direct product cells $C+iD$, $C, D\in\Cc$. It makes $\CC^n$ into a quasi-regular cellular space. 

\vskip .2cm

\noindent The {\bf  ${\bf s}^{(1)}$-stratification $\Sc^{(1)}$} consists of cells $[C,D]$
defined for any {\em face interval}, by which we mean a pair $C,D\in\Cc$ such that $C\leq D$. 
 By definition, $[C,D]$
consists of $z\in\CC^n$ such that:
\begin{itemize}
\item[(a)]  $\Im(z) \in C$.
\item[(b)]  $\Re(z)$ lies in $\pi_C^{-1}(\pi_C(D))$. 
\end{itemize}
Note that $\dim [C,D] = \dim(C) + \dim(D)$. Notice also that for different pairs $(C,D)\neq (C', D')$ we have
$[C,D] \cap [C', D']=\emptyset$, i.e., the strata of $\Sc^{(1)}$ are precisely labelled by the pairs $(C\leq D)$. 
The stratification  $\Sc^{(1)}$ also makes $\CC^n$ into a quasi-regular cellular space.

\vskip .2cm

The stratifications $\Sc^{(\nu)}$, $\nu=1,2$,  were introduced and studied by Bj\"orner and Ziegler
\cite{BZ}. They use a slightly different definition, based on the following  two complex generalizations
of the sign  function: the obvious one
\[
{\bf s}^{(2)}: \CC \lra \{+,-,0\}^2, \quad {\bf s}^{(2)} (a+bi) = (\sgn (a), \sgn(b)),
\] 
and the less obvious one
\[
{\bf s}^{(1)}: \CC \lra\{ i,j,+,-,0\}, \quad {\bf s}^{(1)} (a+bi) 
=\begin{cases}
i & \text{if } b>0,\\
j & \text{if } b<0,\\
+ &  \text{if } b=0 \text{ and } a>0,\\
- &  \text{if } b=0 \text{ and } a<0,\\
0&  \text{if } b=a=0.
\end{cases}
\]

\begin{prop}\label{prop:sign-strat}
Two vectors $z,w\in\CC^n$ lie in the same stratum of $\Sc^{(\nu)}$, $\nu=1,2$, if and only if 
we have ${\bf s}^{(\nu)}(f_H(z))={\bf s}^{(\nu)}(f_H(w))$ for each $H\in\Hc$. 
\end{prop}

\noindent {\sl Proof:} Obvious for $\nu=2$. For $\nu=1$ this is the content of Theorem 5.1(ii)
of \cite{BZ}. \qed

\vskip .2cm

We will view each $\Sc^{(\nu)}$, $\nu=0,1,2$, as a poset, i.e., as the set of strata with the partial
order $\leq$ defined by $S'\leq S$ iff $S'\subset\overline S$.

\begin{prop}\label{prop:3-strat}
\begin{itemize}
 \item[(a)]  Denoting by $\prec$ the relation of refinement of stratifications, we have
\[
\Sc^{(2)} \prec \Sc^{(1)} \prec \Sc^{(0)}. 
\]
In particular, we have order preserving maps of posets $\Sc^{(2)} \to \Sc^{(1)} \to \Sc^{(0)}$ describing which stratum of
$\Sc^{(\nu-1)}$ contains a given stratum of $\Sc^{(\nu)}$, $\nu=1,2$. These maps are as follows:

\item[(b)]  The stratum of $\Sc^{(0)}$ containing a cell $[C,D]\in \Sc^{(1)}$, $C\leq D$,
is $\LL_\CC(D)^\circ$. 

\item[(c)] The stratum of $\Sc^{(1)}$ containing a cell $D+iC \in \Sc^{(2)}$, is $[C, C\circ D]$. 

\item[(d)]  For two cells $[C',D']$ and $[C,D]$ of $\Sc^{(1)}$ we have $[C',D']\leq [C,D]$, if and only if
$C'\leq C$ and $C\circ D'\leq D$. 
 \end{itemize}
\end{prop}
 
\vskip .2cm

\noindent {\sl Proof:} (a) Note that  the stratification $\Sc^{(0)}$ can also be described
in the style of Proposition \ref{prop:sign-strat},  by using
 the sign-type
function
\[
{\bf s}^{(0)}: \CC \lra \{ 0, * \}, \quad {\bf s}^{(0)}(z)=\begin{cases}
0, & \text{if } z=0,\\
*, & \text{if } z\neq 0.
\end{cases}
\]
The statement follows since the three stratifications of $\CC$ induced by the three sign functions
${\bf s}^{(\nu)}$, $\nu=0,1,2$, refine each other as claimed. 

 (b) The $\RR$-linear span of $\pi_C^{-1}(D/\LL(C))$,
the allowable range for $\Re(z)$, $z\in [C,D]$,  is $\LL(D)$. So $[C,D]$
cannot be contained in any $\CC$-linear subspace strictly smaller than $\LL_\CC(D)$.
Also, by construction it is indeed contained in $\LL_\CC(D)$, since the allowable
range for $\Im(z)$ is $C\subset \overline D$. 

(c) This follows by definition of $C\circ D$. 

(d)  This is  Proposition 5.2 of \cite{BZ}. \qed

 \section{Constructible complexes   on arrangements}
 
 We keep the notation of the previous section, in particular, concerning the three stratifications
 $\Sc^{(\nu)}$ of $\CC^n$, $\nu=0,1,2$,  induced by the arrangement $\Hc$. 
 
  \vskip .3cm
  
  \noindent {\bf A. Constructible sheaves.}
A  sheaf $\Fc\in \Sh_{\CC^n}$ will be called $\Sc^{(\nu)}$-{\em smooth},
if it is locally constant on each stratum of $\Sc^{(\nu)}$. 
We will say that $\Fc$ is $\Sc^{(\nu)}$-{\em constructible},
if it is $\Sc^{(\nu)}$-smooth with finite-dimensional stalks. 
A complex $\Fc\in \DSh_{\CC^n}$ will be called $\Sc^{(\nu)}$-smooth
(resp. $\Sc^{(\nu)}$-constructible), if 
if all the cohomology sheaves of $\Fc$ are $\Sc^{(\nu)}$-smooth
(resp. $\Sc^{(\nu)}$-constructible).
We denote
by 
 $D^b_{\Sc^{(\nu)}}\Sh_{\CC^n}\subset \DSh_{\CC^n}$
the full subcategory of   $\Sc^{(\nu)}$-constructible complexes.

In particular, for $\nu=2,1$, since $\Sc^{(\nu)}$  is cellular, an $\Sc^{(\nu)}$-smooth
sheaf   $\Fc$ is uniquely defined by its cellular stalks and 
generalization maps which we denote,   respectively, as follows:
\be\label{eq:S-1-LAD}
\begin{gathered}
\Fc\bigl|_{D+iC}, \quad \gamma^\Fc_{D'+iC', D+iC}: \Fc\bigl|_{D'+iC'}\lra
 \Fc\bigl|_{D+iC}, \,\,\, D'\leq D, \,\, C'\leq C;
 \\
\Fc_{[C,D]}= \Fc\bigl|_{[C,D]}, \quad \gamma^\Fc_{[C',D'], [C,D]}: \Fc_{[C',D']}\lra \Fc_{[C,D]}, \,\,\,
[C',D']\leq[C,D].
\end{gathered}
\ee
 We use similar notations for the case when $\Fc$ is an $\Sc^{(\nu)}$-smooth complex.
In this case the stalks are complexes of vector spaces and we can assume,
by using Proposition \ref{prop:cellullar-LAD}(c), that the generalization maps are
morphisms of complexes forming a representation of the poset of cells. 

\vskip .3cm

\noindent {\bf B. From $\Sc^{(1)}$-smoothness to $\Sc^{(0)}$-smoothness.}
We will need to compare the conditions of smoothness with respect to different stratifications.

\begin{Defi}
An inclusion of ${\bf s}^{(1)}$-cells $[C',D']\leq [C,D]$ will be called {\em elementary}, if one of the two following
cases hold:

\begin{itemize}
\item[(1)] $C'\leq C$ and $D'=D$, so we have a flag $C'\leq C\leq D$.

\item[(2)] $C=D$, $\LL(D)=\LL(D')$, $\dim(C')=\dim(D')-1$, and $C'\leq D$. In other words, $C'$ is a codimension 1
``wall" separating $D$ and $D'$.

\end{itemize}
\end{Defi}

\begin{prop}\label{prop:S1-to-S0}
 Let $\Fc$ be an $\Sc^{(1)}$-smooth sheaf. The following are equivalent:
\begin{itemize}
\item[(i)] $\Fc$ is $\Sc^{(0)}$-smooth.

\item[(ii)] The map $\gamma^\Fc_{[C',D'], [C,D]}$ is an isomorphism for each elementary inclusion
$[C',D']\leq [C,D]$. 
\end{itemize}
\end{prop}

\noindent {\sl Proof:} The tautological equivalent of (i) is
\begin{itemize}
\item[(iii)] The map $\gamma^\Fc_{[C',D'], [C,D]}$ is an isomorphism for each inclusion $[C',D']\leq [C,D]$
such that $\LL(D')=\LL(D)$.
\end{itemize}
Indeed, by Proposition \ref{prop:3-strat}(b), inclusions in (iii) are precisely all inclusions of
${\bf s}^{(1)}$-cells in the same ${\bf s}^{(0)}$-stratum. Clearly (iii) is stronger than (ii).
So we need to prove that (ii) implies (iii). 

To prove (iii), it is enough to fix $L\in\Lc$ and to concentrate on
  inclusions with $\LL(D)=\LL(D')=L$. For this, we do not need to consider
any faces outside $L$, so we can and will  assume that $L=\RR^n$ and $D$ and $D'$ are $n$-dimensional,
i.e., are chambers.

 Let
  $\delta(D,D')$ be the {\em chamber distance}
between $D$ and $D'$, i.e., the minimal length of a sequence
\[
D=D_0, D_1, \cdots, D_l=D'
\]
such that each $D_i$ is a chamber, and each $D_p, D_{p+1}$ have a codimension 1 face in common. 
Thus $\delta(D,D')=0$ means that $D'=D$. 

Let us prove (iii) by induction on $\delta(D,D')$.  If $\delta(D,D')=0$, i.e., $D'=D$, then
our inclusion is of type (1).  
Consider now an arbitrary inclusion
\[
[C',D']\leq [C,D],\,\,\, C'\leq C, \,\, C\circ D'=D,
\]
with $D',D$ being chambers (open faces).

\begin{figure}
 \centering
  \begin{tikzpicture}[scale=0.3]
 \node (C') at (0,0){$\bullet$}; 
 \node (C) at (14,7){}; 
 \node (c) at (7,3.5){$\bullet$};
 \node at (7.5,3) {$c$};
 \node at (14.5,7.5) {$C$}; 
 \draw (C') -- (C); 
 \node (C1) at (10,10){};
 \draw (C') -- (C1); 
 \node at (10.5,10.5) {$C_1$}; 
 \node (N1) at (6,13){};
 \draw (C') -- (N1); 
 
 \node (N2) at (-5,11){};
 \draw (C') -- (N2); 
 \node (N3) at (-8,7){};
 \draw (C') -- (N3); 
 \node at (-1,9.5){$\cdots$};
 \node (d') at (-5,6){$\bullet$};
 \draw (c) -- (d');
 \node at (-5,7){$d'$};
 \node at (10,7){$D$};
 \node at (6.5,10){$D_1$};
 \node at (-7,10){$D'$}; 
 \node at (0.5, -1) {$C'$};
 \end{tikzpicture}

  \caption{Reduction of an inclusion $[C',D']\leq [C,D]$.}
   \label{fig:red}
 \end{figure}
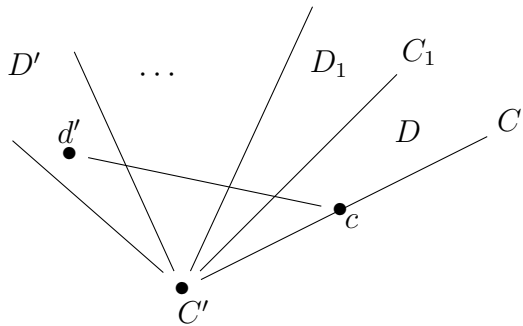

Take generic points $c\in C, d'\in D'$ and form the straight line interval $[c,d']$ oriented towards $d'$.
By Proposition \ref{prop:composition} (b), this interval, after leaving $C$, first hits $D=C\circ D'$. Since $D'$ is a chamber,
 our interval will, after leaving $D$, hit some cell $C_1\geq C'$ of positive codimension, then 
 a chamber $D_1\geq C_1$,   and so on, see Fig. \ref{fig:red}.
 Note, that by choosing $d'$ in a generic enough way, we can ensure that   $C_1$
 is of codimension 1, which we will assume. 
  Note also that $\delta(D_1, D') < \delta(D, D')$. 
 
 We then have a commutative square in (the category corresponding to) the poset
 $(\Sc^{(1)}, \leq)$:
 \[
 \xymatrix{
 [C',D']
 \ar[d]_q
  \ar[r]^p&[C_1, D_1]\ar[d]^r
 \\
 [C,D] \ar[r]_s& [D,D]. 
 }
 \]

 \noindent In this square, $q$ is our inclusion in question, $p$ is an inclusion with smaller $\delta$, while
  $r$ is an inclusion of type (2), and $s$ is an inclusion of type (1). 
  Now, the generalization maps for
 $\Fc$ can be seen as a functor $\gamma^\Fc: (\Sc^{(1)}, \leq)\to\on{Vect}_\k$, so
 the above  square gives a commutative square of vector spaces. By (ii) and our inductive
 assumption,   $p,r,s$ are taken
 by $\gamma^\Fc$ into isomorphisms. 
   It remains to
  deduce that $\gamma^\Fc$ takes $q$ to an isomorphism, invoking the following obvious

\begin{lem}\label{lem:3-4}
 If, in a commutative square of morphisms in any category, 3 out of 4 arrows
 are isomorphism, then the fourth arrow is an isomorphism as well. \qed
\end{lem}

Proposition \ref{prop:S1-to-S0} is proved.

\vskip .3cm

\noindent {\bf C. Complexes of cohomology with support and their linear algebra data.}

Let $\Fc\in D^b_{\Sc^{(0)}}\Sh_{\CC^n}$. Define
\be
\Rc^\bullet_\Fc =  \underline{R\Gamma}_{\RR^n}(\Fc)[n] \,\,\,\in \,\,\, D^b_\Cc(\Sh_{\RR^n}).
\ee
This is a cellular complex on $\RR^n$ (with respect to $\Cc$, the stratification by faces of $\Hc$). 
Therefore it can be defined by the linear algebra data consisting of complexes and generalization maps
\be
E^\bullet_C = E^\bullet_C(\Fc) = R\Gamma(C, \Rc^\bullet_\Fc), \quad \gamma_{C'C}: E^\bullet_{C'}\to E^\bullet_C,
\quad 
C'\leq C\in\Cc
\ee
forming a representation of $(\Cc,\leq)$ in complexes of vector spaces. 

More generally, let $C\in\Cc$ be a face of codimension $d$. Consider the ``tube face"
\be
j_C: \RR^n+iC  \hookrightarrow\CC^n
\ee
and form the complex
\be\label{eq:E-C-complex}
\Ec^\bullet_C \,\,=\,\, \Ec^\bullet_C(\Fc) \,\,=\,\, \underline{R\Gamma}_{\RR^n+iC}(\Fc)[d]
\,\, = \,\, j_{C*}j_C^!(\Fc)[d].
\ee
Thus $\Ec_C^\bullet$ is an $\RR$-constructible complex
on $\CC^n$,  supported on the closure $\RR^n + i \overline C$. 
We denote by $p_C, \overline p_C$ the composite projections
\be
\begin{gathered}
 \RR^n+ iC \buildrel \Re \over\lra \RR^n\lra \RR^n/\LL(C)\\
\RR^n + i \overline C  \buildrel \Re \over\lra \RR^n\lra \RR^n/\LL(C)
\end{gathered}
\ee

\begin{prop}\label{prop:E-C-red}
The complex $\Ec_C^\bullet$, considered as a complex of sheaves on $\RR^n + i\overline C$,
has the form $\overline p_C^*\Ec_C^{\bullet, \on{red}}$, where
$\Ec_C^{\bullet, \on{red}}$ is a  complex on $\RR^n/\LL(C)$ cellular with respect
to the stratification by chambers of the quotient arrangement $\Hc/\LL(C)$. 
\end{prop}

\noindent {\sl Proof:} It is enough to show that $j_C^!\Fc$ has the form $p_C^* \Ec_C^{\bullet, \on{red}}$
for some $\Ec_C^{\bullet, \on{red}}$ as above, since the derived direct image extension from
$\RR^n + iC$ to $\RR^n + i \overline C$ will then proceed along the directions in which 
$p_C^* \Ec_C^{\bullet, \on{red}}$ is constant. 

In the remainder of the proof the word ``manifold" will mean a $C^\infty$-manifold. 
For a complex of  sheaves $\Gc$ on a  manifold $X$ we denote by $\SS(\Gc)\subset T^*X$
the micro-support of $\Gc$, see \cite[Ch. V]{KS}. Our desired result about $j_C^!\Fc$ can be
reformulated by saying that
\be\label{eq:SS-reform}
\begin{gathered}
SS(j_C^!\Fc) \,\,\subset \,\, \bigcup_{\Lambda\in\Lc_{\Hc/\LL(C)}} T^*_{p_C^{-1}(\Lambda)} (\RR^n+iC)  \,\, =
\\
=\,\, \bigcup_{L\in\Lc, \,\, L\supset C} T^*_{L+iC}(\RR^n+iC) \,\,=\,\, \bigcup_{L\in\Lc, \,\, L\supset C}  (L+iC) \times L^\perp,
\end{gathered}
\ee
where $L^\perp\subset\RR^{n*}$ is the orthogonal to $L$. 
We recall the  estimate for the micro-support of the direct image \cite[Cor. 6.4.4]{KS}
specialized to the case of a locally closed embedding.

For a locally closed submanifold $M$ of a manifold $X$ and a subset $S\subset X$ we denote by $K_M(S)$ the
{\em normal cone} to $S$ along $M$, see \cite[Def. 4.1.1]{KS}. It is a closed conic subset of $T_M X$,
 the normal bundle to $M$ in $X$. We note the properties:
 \be\label{eq:normal-add}
 K_M(S_1\cup S_2) \,\,=\,\, K_M(S_1)\cup K_M(S_2), \quad S_\nu\subset X,\,\,\nu=1,2; 
 \ee
 \be
  \begin{gathered}\label{eq:normal-mult}
 K_{M_1\times M_2}(S_1\times S_2) = K_{M_1}(S_1)\times K_{M_2}(S_2) \subset
 \\\subset  T_{M_1}(X_1) \times
 T_{M_2}(X_2) = T_{M_1\times M_2}(X_1\times X_2),
  \quad M_\nu, S_\nu\subset X_\nu. 
 \end{gathered}
 \ee
Given a locally closed embedding of manifolds $f: Y\to X$, the conormal bundle $T^*_YX$ is a Lagrangian submanifold
in $T^*X$, and so its normal bundle there can be written as
\[
T_{T^*_YX} (T^*X)\,\, \simeq \,\, T^*(T^*_YX). 
\]
In particular, the projection $T^*_YX\to Y$ gives a closed embedding
\[
T^*Y\,\, \subset \,\, T^*(T^*_YX) \,\,\simeq \,\, T_{T^*_YX} (T^*X). 
\]
For a complex of sheaves $\Fc$ on $X$, Cor. 6.4.4 and Prop. 6.2.4 of \cite{KS} give:
\be\label{eq:estimate-SS}
\SS(f^!\Fc)\,\, \subset\,\,  K_{T^*_YX}(\SS(\Fc)) \cap T^*Y, 
\ee
the intersection inside $T_{T^*_YX}(T^*X)$.

We now specialize this to
\[
f=j_C: Y=\RR^n+iC \lra \CC^n= X. 
\]
By our assumptions on $\Fc$,
\[
\SS(\Fc)\,\, \subset\,\, \bigcup_{L\in\Lc} T^*_{L_\CC}\CC^n \,\,=\,\, \bigcup_{L\in\Lc} L_\CC \times L_\CC^\perp
\,\,\subset\,\, \CC^n \times\CC^{n*}\,\,
= \,\,T^*\CC^n \,\,=\,\, T^*X. 
\]
We further have
\[
T^*_YX \,\,=\,\, T^*_{\RR^n+iC}\CC^n \,\,=\,\,
(\RR^n+iC) \times i\LL(C)^\perp \,\,\subset \,\, \CC^n \times\CC^{n*}\,\,=\,\, T^*X.
\]
Therefore
\[
T_{T^*_YX}(T^*X) \,\,=\,\, 
\bigl( (\RR^n+iC)\times i\LL(C)^\perp\bigr) \times
\bigl( i(\RR^n/\LL(C)) \times (\RR^{n*} + i \LL(C)^*)\bigr),
\]
in which $T^*Y = T^*(\RR^n+iC)$ is embedded as the product of all the factors except $i\LL(C)^\perp$ and
$i(\RR^n/\LL(C))$, the coordinates corresponding to these factors being put to $0$. 
Applying \eqref{eq:normal-add}, we find
\be\label{eq:K-union}
K_{T^*_YX}(\SS(\Fc))
\,\, \subset\,\,  \bigcup_{L\in\Lc} K_{ (\RR^n+iC) \times i\LL(C)^\perp} (L_\CC \times L_\CC^\perp).
\ee

\begin{lem}\label{lem:dichotomy}
For $L\in \Lc$ we have one of two possibilities:
\begin{enumerate}
\item[(1)] $L_\CC \cap(\RR^n+iC)=\emptyset$; 
\item[(2)] $L\supset C$. 
\end{enumerate}
\end{lem}

Assuming the lemma, we note that in the case (1), $L_\CC \times L_\CC^\perp$ does not meet $(\RR^n+iC)\times i\LL(C)^\perp$
and so will not contribute to the union in \eqref{eq:K-union}. In the case (2) we write, using \eqref{eq:normal-mult}:
\[
\begin{gathered}
K_{(\RR^n+iC)\times i \LL(C)^\perp} (L_\CC\times L_\CC^\perp) \,\,=\,\, K_{(\RR^n+iC) \times (0+i \LL(C)^\perp)}
\bigl( ( L+iL)\times (L^\perp +iL^\perp)\bigr) \,\,=
\\= \,\,
\bigr( K_{\RR^n}(L) +i K_C(L)\bigr)\, \times\, \bigr( K_0(L^\perp) + i K_{\LL(C)^\perp}(L^\perp)\bigr) \,\,= 
\\
= \,\, \bigl( L+ i(C\times (\RR^n/\LL(C)))\bigr) \, \times \bigr( L^\perp + i L^\perp\bigr), 
\end{gathered}
\]
since $L^\perp \subset \LL(C)^\perp$. Now, putting the coordinates from $i\LL(C)^\perp$ and $i(\RR^n/\LL(C))$ to be $0$,
we get $(L+iC)\times L^\perp$, which is the contribution of $L\supset C$ into \eqref{eq:SS-reform}. 
This proves Proposition \ref{prop:E-C-red} modulo Lemma \ref{lem:dichotomy}.

\vskip .2cm

\noindent {\sl Proof of Lemma \ref{lem:dichotomy}:} 
By taking intersections, the lemma reduces to the case when $L=H$ is a hyperplane from $\Hc$,
which we now assume. Let $f_H: \RR^n\to\RR$ be a linear equation of $H$. If (1) does not hold for $H$, then there are
$b\in\RR^n, c\in C$ such that
\[
0\,\,=\,\, f_H(b+ic) \,\,=\,\, f_H(b)+i f_H(c).
\]
Since $f_H(b), f_H(c)\in\RR$, this implies that $f_H(c)=0$ and so $f_H|_C=0$ and $H\supset C$. 
Lemma \ref{lem:dichotomy} and Proposition \ref{prop:E-C-red} are proved.

\vskip .2cm

Since $\Ec_C^{\bullet, \on{red}}$  is a cellular complex, it is determined by linear algebra data
which we denote
\be
E^{C\bullet}_D,\,\,\, D\in\Cc_{\Hc/\LL(C)}, \quad \gamma^C_{D'D}: E^{C\bullet}_{D'}\lra
E^{C\bullet}_D,\,\,\, D'\leq D. 
\ee

 \begin{prop}\label{prop:E-C-D-lift}
  We have canonical identifications (quasi-isomorphisms of complexes of $\k$-vector spaces
 compatible with the maps $\gamma$)
 \[
 E^{C\bullet}_D \simeq E^\bullet_{\overline \pi_C^{-1}(D) }\otimes_\k \OR(C), \quad \gamma^C_{D'D} =
  \gamma_{\overline\pi_C^{-1}(D'), \overline\pi_C^{-1}(D)}\otimes\on{Id}. 
 \]
 Here $\overline \pi_C$ is the isomorphism of posets from
 \eqref{eq:overline-pi-C}. 
    \end{prop}
  
    \noindent {\sl Proof:} 
   Let $K=\overline \pi_C^{-1}(D)$, so $K \geq C$. We start  with recalling the definitions of 
   $E^{C\bullet}_D$ and $E^\bullet_K$ side by side. 
   First,  we remind the notation $d= \on{codim}_{\RR^n}(C)$. 
   Let also $\RR^n_{\geq C}=\bigcup_{K'\geq C} K'$, which is an open subset of $\RR^n$. 
   
      The natural projection $\pi: K\to D$ has contractible
  fibers. The cohomology sheaves of the complex $\underline {R\Gamma}_{\RR^n+iC}(\Fc)[d]$ are
   constant on $K+iC$, while
  the cohomology sheaves of   $\underline{R\Gamma}_{\RR^n}(\Fc)[n]$ 
  are constant on $K$. 
  Denote
  \[
  \Gc'= \underline{R\Gamma}_{\RR^n +  i\LL(C)} (\Fc)[d], \quad \Gc=\Gc'|_{\RR^n_{\geq C}+i\LL(C)}. 
  \]
  Because $C$ is open in $\LL(C)$ and $\RR^n_{\geq C}$ in $\RR^n$, we have
  \[
  E^{C\bullet}_D = \Ec_C^{\bullet, \on{red}} |_D = \underline{R\Gamma}_{\RR^n + iC} (\Fc)[d]
  \bigl|_{K+iC} = 
    \Gc \bigl|_{K+iC}
  \]
  (restriction of a cellular complex to a cell,  considered as a complex of vector spaces). 
  On the other hand,  we have
    \[
    E_K^\bullet \,\, =\,\, \Rc_\Fc|_K \,\,=\,\,
    \underline{R\Gamma}_{\RR^n}(\Fc)[n] \bigl|_{K}\,\, =\,\,  
  \underline{R\Gamma}_{\RR^n}  \Gc'  [n-d]\bigl| _{K} \,\,=\,\,   \underline{R\Gamma}_{\RR^n_{\geq C}}  \Gc  [n-d]\bigl| _{K} . 
  \]
  To compare these, we note that $K =  K+i0$ lies in $K+i \overline C $ but not
  in $K+iC$ proper. We use the following

  \begin{lem}\label{lem:non-char}
  For any $K'\geq C$ the cohomology sheaves of $\Gc$ are
    constant on the whole of $K'+i\LL(C)$ (and hence the restriction of $\Gc$ to $K'+i\LL(C)$
    is quasi-isomorphic to a constant complex of sheaves).
  \end{lem}
  This is a consequence of
    the following analog of
     Lemma 
    \ref{lem:dichotomy}.
   
   \begin{lem}\label{lem:3.14}
  Let $H$ be a hyperplane from $\Hc$ and $K'\geq C$. Then $H_\CC$ meets $K'+i\LL(C)$ if and only
   if $H_\CC$ contains $K'+i\LL(C)$. 
   \end{lem}

   \noindent {\sl Proof of  lemma \ref{lem:3.14}:} 
   If $f_H(k'+iy)=0$ for some $k'\in K$, $y\in \LL(C)$, then by taking the real part, we find
   $f_H(k')=0$ and therefore     $(f_H)|_{K'} =0$. Now, since $C\leq K'$, $(f_H)|_C=0$  as well and therefore
   $(f_H)|_{\LL(C)}=0$. Therefore $f_H(k'+iy) =f_H(k')+if_H(y)=0$ for each
    $k'\in K'$ and  $y\in\LL(C)$. \qed  
    
    \vskip .2cm
    
    Consider now the Cartesian square of  embeddings (transverse intersection)
  \[
  \xymatrix{
  K+i\LL(C) \ar[r]^v & \RR^n_{\geq C}+i\LL(C)
  \\
  K=K+i0
  \ar[u]^u \ar[r]_{\hskip -.5cm t} & \RR^n_{\geq C}=\RR^n_{\geq C}+i0. \ar[u]_s
  }
  \]
  Lemma \ref{lem:non-char} implies that $s$ is non-characteristic for $\Gc$ and therefore we have ``local Poincar\'e duality"
  $s^!\Gc \simeq s^*\Gc \otimes \orr(C)[d-n]$, see \cite[cor. 5.4.11]{KS}. Therefore we identify
  \[
  \begin{gathered}
  E^\bullet_K =\Gamma(K,  t^*s^!\Gc[n-d]) \,\,\simeq \,\,  \Gamma(K, t^*s^*\Gc\otimes\orr(C)) \,\, =\,\, 
   \\
=\,\,   \Gamma(K,  u^*v^*\Gc\otimes\orr(C))\buildrel \eqref{lem:non-char}\over = 
\Gamma(K+i\LL(C), v^*\Gc\otimes\orr(C) ) \,\,= \,\,E_D^{C\bullet}\otimes\orr(C).
   \end{gathered}
  \]
    Identification of the maps $\gamma^C_{D'D}$ is done similarly. 
       Proposition \ref{prop:E-C-D-lift} is proved.

       \begin{rem}  Let us indicate another, perhaps more geometric view on
      Proposition \ref{prop:E-C-D-lift}.  Let $x\in C$ and $T$ be a small transversal slice to $\LL(C)$
       at $x$, i.e., a subset of the form $S+x$ where $S$ is a small open ball in a linear subspace $M\subset\RR^n$
       such that $M\oplus \LL(C)=\RR^n$. Thus $T$ fits into a diagram of embeddings
       \[
       T\buildrel\epsilon\over\lra \RR^n \buildrel j_0\over\lra \CC^n.
       \]
      The composition of $\epsilon$ with the projection $q: \RR^n\to\RR^n/\LL(C)$
       is an embedding of a small ball into $\RR^n/\LL(C)$. Proposition \ref{prop:E-C-D-lift} can be reformulated
       by saying that
       \[
       (q\epsilon)^*\Ec_C^{\bullet, \on{red}}\,\,=\,\, \epsilon^*\Rc_\Fc\otimes_\k \on{or}(C). 
       \]
      Both sides of this proposed isomorphism can be identified as certain complexes of cohomology with support. First, 
      by transversality and Poincar\'e duality,
       \[
       \epsilon^*\Rc_\Fc \,\,=\,\, \epsilon^!\Rc_\Fc \otimes_\k \on{or}(C)[n-d],
       \]
       and so, since $\Rc_\Fc=j_0^!\Fc[n]$, 
       \[
       \epsilon^*\Rc_\Fc \,\,=\,\, (j_0\epsilon)^!\Fc\otimes_\k \on{or}(C)[2n-d].
       \]
       Second, take another point $y\in C$ and let $T'=T+iy\subset\CC^n$ be the translation of $T$ by $iy$, fitting into
       the diagram of embeddings
       \[
       T'\buildrel\epsilon'\over\lra \RR^n+iC \buildrel j_C\over\lra\CC^n. 
       \]
       Then $T'$ is transversal to the fibers of $p_C$ and therefore $(q\epsilon)^*\Ec_C^{\bullet, \on{red}}$ is identified
       with 
        $(\epsilon')^*\Ec_C^\bullet$, after identification of $T$ with $T'$ via the shift. Again, by transversality of $T'$
        to the fibers of $p_C$ (and the fact that they are canonically oriented, being complex manifolds), we find that
        \[
        (\epsilon')^*\Ec_C^\bullet = (\epsilon')^!\Ec_C^\bullet [2n-2d]  = (j_C\epsilon')^! \Fc[2n-d].
        \]   
        Using the canonical identification $\on{or}(C)^{\otimes 2}=\k$, 
       the proposition thus reduces to the claim that  the complexes of cohomology
      with support in $T$ and $T'$ for $\Fc$ are identified. Further,  Lemma \ref{lem:3.14}
      can be seen as identifying the stratifications induced by $\Hc$ on $T$ and $T'$ and so allows 
      one to prove the desired claim by a homotopy argument, deforming $y$ to $0$. 
         \end{rem}

\begin{cor}\label{cor:stalks-E-C}
(a) The complex $\Ec_C^\bullet$ of sheaves on $\CC^n$
is smooth with respect to the stratification $\Sc^{(1)}$ (and therefore to $\Sc^{(2)}$).

(b) The stalks  and generalization maps of $\Ec_C^\bullet$ on cells of $\Sc^{(2)}$
are given by 
\[
\begin{gathered}
 \Ec_C^\bullet \bigl|_{iC_1+D} = \begin{cases}
  E^\bullet_{C\circ D}\otimes\OR(C), & \text{if } C_1\leq C; \\
  0, & \text{otherwise};
  \end{cases}
  \\
  \gamma^{\Ec_C^\bullet}_{iC'_1+D', iC_1+D} = \gamma_{C\circ D', C\circ D}\otimes\on{Id}, \quad 
  C'_1\leq C_1\leq C, \,\, D'\leq D.
  \end{gathered}
  \]

 (c) The stalks  of $\Ec_C^\bullet$ on cells of $\Sc^{(1)}$
are given by 
\[
 \Ec_C^\bullet\bigl|_{[C_1,D]} =\begin{cases}
E^\bullet_{C\circ D} \otimes \OR(C) & \text{if } C_1\leq C,\\
0 & \text {otherwise}.
\end{cases}
\]
Further, if $[C'_1, D'] \leq [C_1, D]$ is an inclusion of two $\Sc^{(1)}$-cells in the support of $\Ec^\bullet_C$,
 then $C\circ D' \leq C\circ D$, and the corresponding generalization map for $\Ec^\bullet_C$ has the form
 \[
\gamma^{\Ec^\bullet_C}_{[C'_1, D'], [C_1, D]} = \gamma_{ C\circ D', C\circ D}\otimes\on{Id}. 
 \]
 
\end{cor}

     \noindent   In particular, putting $C_1=D=\{0\}$ in part (b) of the corollary, we get  quasi-isomorphisms
    \be
    \begin{gathered}
    R\Gamma(\CC^n, \Ec_C^\bullet) \,\,\simeq \,\, R\Gamma(\RR^n/\LL(C), \Ec^{\bullet, \on {red}}_C)
    \,\,\simeq \\
   \simeq \,\, E^{C\bullet}_0   \,\,=\,\, E^\bullet_C\otimes\OR(C). 
    \end{gathered}
    \ee
    At the very right we have the complexes appearing in the description of the cellular complex
    $\Rc^\bullet_\Fc$. 
    
    \vskip .3cm
    
    \noindent {\sl Proof of Corollary \ref{cor:stalks-E-C}:} Part (a) follows  from Proposition
    \ref{prop:E-C-D-lift}, because cells of $\Sc^{(1)}$ are the the lifts of cells of  the quotient arrangements,
    by their definition in \S \ref{sec:back-arr}D. Part (b) is a simple translation of Proposition 
      \ref{prop:E-C-D-lift} by using the interpretation of the $\circ$ operation in Proposition
      \ref{prop:circ-pi}.   The inclusion $C\circ D' \leq C\circ D$ in part (c)  is proved as follows.
      
      The condition $[C'_1, D_1] \leq [C_1, D]$ means, by  Proposition \ref{prop:3-strat}(d), that $C'_1\leq C$
      and $C_1\circ D'\leq D$.   The inclusion $C_1\leq C$ implies that $C\circ C_1 = C$. 
      Since the operation $\circ$ is associative and monotone in the second argument,
       we then have
      \[
      C\circ D' \,\,=\,\,  (C\circ C_1) \circ D' \,\, = \,\, C\circ (C_1\circ D')\,\, \leq \,\, C\circ D. 
      \]
      Once this inclusion is established, 
      part (c) becomes a reformulation of part (b) using the explicit relation between strata of $\Sc^{(2)}$
      and $\Sc^{(1)}$  as given in Proposition \ref{prop:3-strat}(c). \qed

    \section{Perverse sheaves and  double quivers}
    
    \noindent {\bf A. The Cousin complex.}
    We keep the notations of the previous section. 
    Put $X=\CC^n$ and consider the filtration of $X$ by closed subspaces
    \be\label{eq:tube-filtr}
    X_d \,\,=\,\, \bigcup_{\dim(C)\leq d} \RR^n + iC,
    \ee
    so that
    \[
    Y_d\,\,=\,\, X_d\setminus X_{d-1} \,\,\, = \,\,\, \bigsqcup_{\dim(C)=d} \RR^n + iC. 
   \]
    Let $\Fc\in D^b_{\Sc^{(0)}} \Sh_{\CC^n}$ be an $\Sc^{(0)}$-constructible complex. 
    We consider the corresponding homological Postnikov system for $\Fc$.  The associated complex
    in the derived category with total object $\Fc$ has the form
    \be\label{eq:cousin-double-cx}
    \Ec^{\bullet\bullet} = \Ec^{\bullet\bullet}(\Fc) \,\,=\,\,
    \biggl\{ \bigoplus_{\cod(C)=0} \Ec^\bullet_C \buildrel {\underline{\widetilde \delta}}\over\lra
      \bigoplus_{\cod(C)=1} \Ec^\bullet_C \buildrel {\underline{\widetilde\delta}}\over\lra
     \cdots  \buildrel {\underline{\widetilde\delta}}\over\lra \Ec_0^\bullet\biggr\}.  
    \ee
    Here $\Ec_C^\bullet$ has been defined in \eqref{eq:E-C-complex}. 
     In particular, by applying the functor $R\Gamma(\CC^n, -)$, we get a Postnikov system
     in $D^b\Vfd$ whose associated complex with total object $R\Gamma(\CC^n,\Fc)$ has the form
     \be\label{eq:cousin-1}
     \bigoplus_{\cod(C)=0} E_C^\bullet\otimes\OR(C) 
     \buildrel \widetilde\delta\over\lra 
     \bigoplus_{\cod(C)=1} E_C^\bullet\otimes\OR(C) 
     \buildrel\widetilde \delta\over\lra \cdots  \buildrel\widetilde \delta\over\lra E_0^\bullet \otimes\OR(0).
     \ee
     We will call \eqref{eq:cousin-1} the {\em Cousin complex} of $\Fc$. 
     In particular, the differential $\widetilde\delta$ of this complex splits into matrix elements
     \[
     \widetilde\delta_{CC'}: E_C^\bullet\otimes\OR(C) \lra E_{C'}^\bullet\otimes\OR(C'), \quad
      C'<_1 C. 
     \]
     Note that for a codimension 1 face $C'$ of a convex polyhedron $C$ we have a canonical
     coorientation, i.e., a canonical trivialization of $\OR(C)^*\otimes\OR(C')$, which allows
     us to write the matrix elements as maps
     \be\label{eq:delta-cooriented}
     \delta_{CC'} = \delta_{CC'}^\Fc: E^\bullet_C \lra E^\bullet_{C'}, \quad
     C'<_1 C. 
\ee

\begin{prop}\label{prop:delta-commute}
The maps of complexes $\delta_{CC'}$ commute, i.e., extend to a contravariant representation
of the poset $(\Cc, \leq)$ in complexes of $\k$-vector spaces. 

\end{prop}

\noindent{\sl Proof:} Indeed, the condition $\widetilde\delta^2=0$ implies that the $\widetilde\delta_{CC'}$
anticommute, so the $\delta_{CC'}$ commute by the antisymmetry of the orientation isomorphisms. 
\qed

\begin{prop}\label{prop:E-C-dual}
We have canonical isomorphisms $E^\bullet_C(\Fc^\bigstar) \,\,\simeq \,\, E^\bullet_C(\Fc)^*$
(dual complexes and adjoint maps)
 in the derived category of vector spaces. Further, under these isomorphisms the morphisms of complexes
 $ \delta_{CC'}^{\Fc^\bigstar}$
 (considered as morphisms in the derived category) are dual to 
  $(\gamma^\Fc_{C'C})$.

\end{prop}

\noindent {\sl Proof:} Since the stratification given by $\Hc_\CC$ is invariant under the $\CC^*$-action on $\CC^n$,
we can use cohomology with support in $i\RR^n$ instead of $\RR^n$ to define $\Rc_{\Fc}$ and $\Rc_{\Fc^\bigstar}$.
More precisely, for $t\in\RR$ let $j_t: \RR^n\to\CC^n$ be the   embedding given by the  multiplication with $e^{it}\in\CC^*$. 
Then the sheaves $\Rc_{\Fc, t}= j_t^!\Fc[n]$ on $\RR^n$ for $t\in [0, \pi/2]$ are canonically identified with each other
because they unite into a sheaf on $[0, \pi/2]\times\RR^n$ locally constant (and therefore constant)
  on each interval $[0, \pi/2]\times\{x\}$, $x\in\RR^n$. 
  So we denote $j=j_{\pi/2}$ to be the embedding of $i\RR^n$ and use $j^!$ to define $\Rc_{\Fc^\bigstar}$. 
 
  Since Verdier duality interchanges
$j^!$ and $j^*$, we have a canonical quasi-isomorphism
\[
(\Rc_{\Fc^\bigstar})^\bigstar \,\,= \,\,  \bigl( \underline{R\Gamma}_{i\RR^n}(\Fc^\bigstar)[n]\bigr)^\bigstar
\,\,\simeq \,\,  \Fc|_{i\RR^n}. 
\]
Now, $\Fc|_{i\RR^n}$ can be understood by restricting to $i\RR^n$ the homological Postnikov system
corresponding to the filtration \eqref{eq:tube-filtr}. This gives a complex in derived category with
total object  $\Fc|_{i\RR^n}$ obtained by restricting $\Ec^{\bullet\bullet}(\Fc)$ to $i\RR^n$ which,
after identifying $i\RR^n$ back  with $\RR^n$, gives the following:
\be\label{eq:A}
\bigoplus_{\cod(C)=0} \underline {E_C^\bullet\otimes\OR(C)}_{\overline C} 
\buildrel (\widetilde\delta_{CC'})\over\lra 
\bigoplus_{\cod(C)=1} \underline {E_C^\bullet\otimes\OR(C)}_{\overline C} 
\buildrel (\widetilde\delta_{CC'})\over\lra \cdots
\ee
Indeed, for any $C$ we find from Corollary \ref{cor:stalks-E-C}(b) that
\[
\Ec_C^\bullet|_{i\RR^n} \,\,=\,\,\underline {E_C^\bullet\otimes\OR(C)}_{i\overline C}. 
\]
On the other hand,  by its original definition (involving cohomology with support in $\RR^n$), 
$\Rc_{\Fc^\bigstar}$  is given by complexes of vector spaces
$E_C^\bullet(\Fc^\bigstar)$ and generalization maps $\gamma_{C'C}^{\Fc^\bigstar}$. 
So by Proposition \ref{prop:VD-explicit} the shifted Verdier dual $(\Rc_{\Fc^\bigstar})^\bigstar$
is the complex
\be\label{eq:B}
\bigoplus_{\cod(C)=0} \underline{
E_C(\Fc^\bigstar)^*\otimes\OR(C)}_{\overline C} 
\buildrel 
(\gamma_{C'C"}^{\Fc^\bigstar})^*\over\lra
\bigoplus_{\cod(C)=1} \underline{
E_C(\Fc^\bigstar)^*\otimes\OR(C)}_{\overline C} 
\buildrel 
(\gamma_{C'C"}^{\Fc^\bigstar})^*\over\lra\cdots
\ee
Comparing \eqref{eq:A} and \eqref{eq:B}, we get our statement. \qed

\vskip .3cm

\noindent {\bf B. Perverse sheaves.}
Let 
$\Perv(\CC^n, \Hc)\subset D^b_{\Sc^{(0)}}\Sh_{\CC^n}$
be the full subcategory of perverse sheaves. 
 We  choose the following normalization of the perversity conditions for a complex $\Fc$
(differing by a shift from that of \cite{KS} \S 10.3): 
\begin{itemize}
\item[($P^-$)] For each $p$ the sheaf $\underline H^p(\Fc)$ is supported on a closed
complex subspace of codimension $\geq p$.
\item[($P^+$)] If $l: Z\to X$ is a locally closed embedding of a smooth analytic sub manifold
of codimension $p$, then the sheaf  $\underline H^q(l^!\Fc) = \underline \HH^q_Z(\Fc)|_Z$ is zero for $q<p$. 
\end{itemize}
With respect to this definition, a constant sheaf on $\CC^n$ is perverse, if put in degree 0.
The normalized Verdier duality functor $\bigstar$ interchanges $(P^-)$ and $(P^+)$ and
preserves $\Perv(\CC^n, \Hc)$.

\begin{prop}\label{prop:E-C-perverse}
(a) If $\Fc$ is perverse, then each $E_C^\bullet(\Fc)$ is quasi-isomorphic to one
vector space $E_C(\Fc)$ in degree 0. 

(b) For any $C\in \Cc$ the functor
\[
    E_C: \Perv(\CC^n, \Hc)\lra\Vfd, \quad \Fc\mapsto E_C(\Fc)
    \]
    is an exact functor of abelian categories. 
\end{prop}

\noindent {\sl Proof:} The functor
\[
E^\bullet_C: D^b_{\Hc_C}\Sh_{\CC^n} \lra D^b\Vfd, \quad \Fc\mapsto E_C^\bullet (\Fc)
\]
is an exact functor of triangulated categories. So (b) will follow from (a), and we concentrate
on (a). 

 Because of Proposition \ref{prop:E-C-dual}, it is enough to show that
$H^iE_C^\bullet(\Fc)$ vanishes for $i>0$, since vanishing for $i<0$ will then follow
by duality. Since the $E_C^\bullet(\Fc)$ are the linear algebra data describing the
cellular complex $\Rc_\Fc=\underline{R\Gamma}_{\RR^n}(\Fc)[n]$, we need to show that 
\[
\underline{\HH}^p_{\RR^n}(\Fc) = 0, \quad p>n. 
\]
From the spectral sequence
\[
\underline{H}^i_{\RR^n}(\underline H^j(\Fc))\,\,\Rightarrow \,\, \underline\HH^{i+j}_{\RR^n}(\Fc)
\]
we see that it is enough to show that
\[
\underline H^i_{\RR^n}(\underline H^j(\Fc))=0 \quad \text{for} \quad i>n-j.
\]
But by  $(P^+)$, the sheaf $\underline H^j(\Fc)$ is supported on the union
of complex flats $L_\CC$ of $\Hc$ of codimension $\geq j$, i.e., of dimension $\leq n-j$.
So our statement follows from the next lemma.

\begin{lem}
Let $\Gc$ be an $\Sc^{(0)}$-smooth sheaf on $\CC^n$, with $\dim_\CC (\on{Supp}(\Gc))\leq q$.
Then $\underline H^i_{\RR^n}(\Gc)=0$ for $i>q$. 

\end{lem}

\noindent {\sl Proof of the lemma:} We first consider the case when $\on{Supp}(\Gc)\subset L_\CC$ where $L\subset \RR^n$ is
a  flat of $\Hc$ of dimension $d\leq q$. In this case
$L_\CC\cap \RR^n=L$ has real codimension $d=\dim_\RR(L)$ in $L_\CC$.
In other words, the embedding $L\subset L_\CC$ is isomorphic to $\RR^d\subset\CC^d$.
Now, it is a general  property that for any sheaf $\Gc$   on  $\CC^d = \RR^{2d}$
the local cohomology sheaves $\underline H^i_{\RR^d}(\Gc)$ vanish for $i>d$. 

We now consider the general case and use induction on $d=\dim(\on{Supp}(\Gc))$ and then on the
number $\nu$ of $d$-dimensional flats $L$ of $\Hc$ such that $L_\CC\subset\on{Supp}(\Gc)$. 
The case $\nu=1$ is covered by the previous paragraph. In general, for a fixed $d$ and $\nu>1$, 
we choose one (out of the $\nu$) flat $L$
such that $L_\CC\subset\on{Supp}(\Gc)$, and denote by $i: L_\CC\to\CC^n$ the embedding. 
Then we have a short exact sequence
\[
0\to \Gc'\lra \Gc\buildrel c\over\lra i_* i^{-1}\Gc\to 0,
\]
with $c$ be the canonical map and $\Gc'=\on{Ker}(c)$. Here
 $ i_* i^{-1}\Gc$ has support on one complex flat (and so  its $\underline H^i_{\RR^n}=0$ for $i>q$),
while $\Gc'$ is supported on the union of $\nu-1$ complex flats and so 
its $\underline H^i_{\RR^n}=0$ for $i>q$ by induction. Writing the long exact sequence for sheaves of cohomology
with support,  we deduce that $\underline H^i_{\RR^n}(\Gc)=0$ for $i>q$.
This finishes the proof of the lemma as well as of Proposition \ref{prop:E-C-perverse}. 

\qed

\begin{cor}\label{cor:cousin-complex}
 (a)
If $\Fc\in\Perv(\CC^n, \Hc)$, then each complex $\Ec_C^\bullet(\Fc)$
is quasi-isomorphic to a single sheaf $\Ec_C(\Fc)$ in degree 0, and
  \[
   \Ec^{\bullet}(\Fc) \,\,=\,\,
    \biggl\{ \bigoplus_{\cod(C)=0} \Ec_C \buildrel {\underline{\widetilde \delta}}\over\lra
      \bigoplus_{\cod(C)=1} \Ec_C \buildrel {\underline{\widetilde\delta}}\over\lra
     \cdots  \buildrel {\underline{\widetilde\delta}}\over\lra \Ec_0\biggr\}
    \]
    is a complex of sheaves on $\CC^n$ in the usual sense, quasi-isomorphic to $\Fc$. 
    
    (b) For any $C\in \Cc$ the functor
    \[
    \Ec_C: \Fc\in\Perv(\CC^n, \Hc)\lra\Sh_{\CC^n}, \quad \Fc\mapsto \Ec_C(\Fc)
    \]
    is an exact functor of abelian categories.

\end{cor}

We will call $\Ec^\bullet(\Fc)$ the {\em Cousin resolution} of $\Fc$. 

\vskip .2cm

\noindent{\sl Proof:} (a) follows from  Propositions  \ref{prop:E-C-perverse}(a) and \ref{prop:E-C-D-lift}.
Part (b) follows, similarly to Proposition  \ref{prop:E-C-perverse}(b), from the exactness
of $\Fc\mapsto \Ec_\Fc^\bullet$ on the derived category. \qed

\begin{Defi}
By a {\em double representation} of the poset $(\Cc, \leq)$ we mean a datum
$
Q = (E_C, \gamma_{C'C}, \delta_{CC'})
$
consisting of finite-dimensional $\k$-vector spaces $E_C, C\in C$ and linear maps
\[
\gamma_{C'C}: E_{C'}\to E_C, \quad \delta_{CC'}: E_C\to E_{C'}, \quad C'\leq C, 
\]
so that $(\gamma_{C'C})$ is a covariant representation and $(\delta_{CC'})$ is
a contravariant representation of $(\Cc, \leq)$. 
\end{Defi}

Double representations of $\Cc$
form an abelian category which we denote $\on{Rep}^{(2)}(\Cc)$. 
This category has a perfect duality
\[
Q \mapsto Q^* = (E_C^*, \delta_{CC'}^*, \gamma_{C'C}^*). 
\]

The results of this sections imply that we have an exact functor
\be\label{eq:functor-Q}
\begin{gathered}
\Qc: \Perv(\CC^n, \Hc_\CC) \lra \on{Rep}^{(2)}(\Cc),\\
\Fc \mapsto \Qc(\Fc) = \bigl( E_C(\Fc), \gamma^\Fc_{C'C"}, \delta^\Fc_{CC'}\bigr), 
\end{gathered}
\ee
commuting with duality. 
We will call $\Qc(\Fc)$ the {\em double quiver} associated to the perverse sheaf
$\Fc$. 

Let us note the following converse to Proposition \ref{prop:E-C-perverse}. 

\begin{prop}\label{prop:perversity-char}
Let $\Fc\in D^b_{\Sc^{(0)}}\Sh_{\CC^n}$ be an $\Sc^{(0)}$-constructible complex such that
each $E_C^\bullet(\Fc)$, $C\in\Cc$, is quasi-isomorphic to a single vector space in degree 0.
Then $\Fc$ is perverse. 
\end{prop}

\noindent {\sl Proof:} Our assumptions imply that we have a Cousin resolution $\Ec^\bullet(\Fc)$ of
$\Fc$ as in Corollary \ref{cor:cousin-complex}. So it is enough to show that $\Ec^\bullet(\Fc)$
satisfies $(P^-)$ and $(P^+)$. By construction, $\Ec^p(\Fc)$ is supported on the union of the
$\RR^n+i\overline C$ for $C\in\Cc$, $\on{codim}(C)=p$. Therefore $\underline H^p(\Ec^\bullet(\Fc))$
is supported on the union of complex flats $L_\CC$, $L\in\Lc$ which are contained in the above
union  of the $\RR^n+i\overline C$. Each such $L$ must have codimension $\geq p$. So
$\Ec^\bullet(\Fc)\sim \Fc$ satisfies $(P^-)$. Now, look at $\Fc^\bigstar$. The double quiver corresponding
to $\Fc^\bigstar$, being, by Proposition \ref{prop:E-C-dual}, identified with $\Qc(\Fc)^*$, also consists
of single vector spaces in degree 0. Therefore the above reasoning shows that $\Fc^\star$ satisfies
$(P^-)$ and so $\Fc$ satisfies $(P^+)$ and so is perverse. \qed

\vskip .3cm

\noindent {\bf C. Relation to earlier works.} 
\noindent {\bf (1)} 
 Suppose that   $\k$ (our coefficient field for perverse sheaves) is equal to $\CC$.
By the Riemann-Hilbert correspondence, any 
$\Fc\in\Perv(\CC^n, \Hc)$ can be represented as the solution sheaf of
a  holonomic $\Dc$-module $\Mc$  on $\CC^n$:
\[
\Fc=\underline{R\Hom}_{\Dc_{\CC^n}} (\Mc, \Oc_{\CC^n}).
\]
In this case one can give another, analytic proof of 
Proposition
 \ref{prop:E-C-perverse}. Indeed, the complex $\Rc^\bullet_\Fc$ can be written as
 \[
 \Rc^\bullet_\Fc = \underline{R\Hom}_{\Dc_{\CC^n}}(\Mc, \underline{R\Gamma}_{\RR^n}(\Oc_{\CC^n}))[n].
 \]
 By the classical result of Sato \cite{SKK}
 \[
 \underline{R\Gamma}_{\RR^n}(\Oc_{\CC^n})[n] \,\,\sim\,\, \underline H^n_{\RR^n}(\Oc_{\CC^n}) = 
 j_*\Bc_{\RR^n}\otimes_\CC \OR(\RR^n), 
 \quad
 j: \RR^n\hookrightarrow\CC^n
\]
 reduces to the sheaf of hyperfunctions $\Bc_{\RR^n}$, so
 \[
 \Rc^\bullet_\Fc = \underline{R\Hom}_{\Dc_{\CC^n}}(\Mc, j_*\Bc_{\RR^n}) \otimes_\CC \OR(\RR^n)
   \]
 is  the complex of hyperfunction solutions of $\Mc$. 
 The fact that this complex reduces to a single sheaf, follows from the 
 result of Lebeau \cite{lebeau} 
 (see \cite{honda-schapira} for  the proof of a more general statement)
 which implies that under our assumptions, 
$\underline{\on{Ext}}^q_{\Dc_{\CC^n}}(\Mc,  j_* \Bc_{\RR^n})=0$ for $q>0$.
So $\Rc_\Fc^\bullet$
 is  quasi-isomorphic to the  sheaf of hyperfunction solutions of $\Mc$ in the  non-derived sense. 
 
  In particular, $E_0(\Fc)= (\Rc_\Fc)_0 = \Gamma(\RR^n, \Rc_\Fc)$
  is the space of global hyperfunction solutions of $\Mc$, and
 its dimension was found by Takeuchi \cite{takeuchi} to be
 the sum of multiplicities of $\Fc$ (or $\Mc$) along all
 the  possible components of the
 characteristic variety: 
 \[
\dim_\CC E_0(\Fc) \,\,=\,\, \sum_{L\in \Lc} \on{mult}_{T^*_{L_\CC}\CC^n} \Fc.  
 \]
 This generalizes the classical index formula of Kashiwara-Komatsu in dimension 1,
 see \cite[Th.4.2.7]{kashiwara} \cite{komatsu}. See also \cite{schurmann}, Th. 1.2, for a generalization
 to arrangements of non-linear analytic subvarieties. 
 
 More generally, for any $C\in\Cc$ we deduce, by passing to the transverse slice to $\LL(C)$, that
 \[
 \dim\,\,  E_C(\Fc) \,\,\,=\,\,\,\sum_{L\in\Lc \text{ s.t. } C\subset L} 
 \on{mult}_{T^*_{L_\CC} \CC^n} \Fc. 
 \]
 See \S \ref{sec:excom} for a discussion of low-dimensional cases and identification of the $E_C(\Fc)$
 in these cases, in terms of the standard functors of nearby and vanishing cycles. 
   \vskip .2cm

\noindent {\bf (2)} The vector spaces $E_C(\Fc)$ are the same as  the spaces of
``generalized vanishing cycles"
introduced in \cite{BFS}, Part I, \S 3.3.  Our Proposition
 \ref{prop:E-C-perverse} corresponds to Theorem 3.9 of \cite{BFS}, Part I, 
 which says that the complexes of generalized vanishing cycles reduce
 to single vector spaces, while the part of  Proposition
 \ref {prop:E-C-dual} pertaining to the spaces $E_C$,  
 corresponds  to Theorem 3.5 of \cite{BFS}, Part I.  
 Note that in our approach the more immediate maps
 among the $E_C$ are the $\gamma_{C'C}$, while
 in \cite{BFS}, Part I,  \S 3.11 it is the $\delta_{CC'}$ (the ``variation maps"). This is due to the fact that the definition of 
 {\it op. cit.} is ''Verdier dual'' to ours: the spaces there 
 are ours $E_C(\Fc^*)^*$.

 \section{Functorialities of the double quiver}
  
 \noindent {\bf A. Hyperbolic restriction.} Let $L\in\Lc$ be a flat of $\Hc$ of codimension $d$.
 Consider the embeddings
 \[
 L_\CC = L+iL \buildrel k\over\lra \RR^n+iL \buildrel j\over\lra\CC^n. 
 \]
 Recall that $L$ carries the induced arrangement $\Hc\cap L$ with poset of faces $\Cc_{\Hc\cap L} \simeq \Cc ^{\leq L}$
 naturally a subposet of $\Cc=\Cc_\Hc$. We have therefore the restriction functor
 \[
 \Rep^{(2)}(\Cc) \lra \Rep^{(2)} (\Cc^{\leq L}), \quad Q \mapsto Q^{\leq L},
 \]
 which takes a double quiver $Q=(E_C, \gamma_{C'C}, \delta_{CC'})$ to its subdiagram involving only $C', C$ which
 are contained in $L$. 
 
 \begin{prop}\label{prop:hyp-res}
 Let $\Fc\in\Perv(\CC^n, \Hc)$ be an $\Sc^{(0)}$-smooth perverse sheaf with
 double quiver $Q$. Then $k^* j^! \Fc[d]$ is an object of $\Perv(L_\CC, \Hc\cap L)$, and
 its associated double quiver is $Q^{\leq L}$. 
 \end{prop}
 
 \begin{rem}
 The statement about perversity of $k^* j^! \Fc[d]$ can be seen as a real analytic  analog
 of the main result of Braden \cite{braden} (see also \cite{drinfeld-gaitsgory} for a more in-depth treatment).
  Extending the terminology of \cite{braden}, we will
 call the perverse sheaf $k^* j^! \Fc[d]$ the {\em hyperbolic restriction} of $\Fc$ to $L_\CC$.
 \end{rem}
 
 \noindent {\sl Proof of Proposition \ref{prop:hyp-res}:} First, we notice that $k^*j^!\Fc$ is an $\Sc^{(0)}$-smooth
 constructible complex on $L_\CC$. This follows at once from the    estimate  \eqref{eq:estimate-SS} for the
 singular support of the inverse image.  
   Thus $k^* j^! \Fc[d]$
 can be described by a double quiver of complexes of vector spaces, and we analyze this double quiver.
 
 Recall that $\Fc$ is represented by its Cousin resolution $\Ec^\bullet = \Ec^\bullet(\Fc)$ with $\Ec^p$ being the direct
 sum of cellular sheaves $\Ec_C$ for $C$ running over faces of $\Hc$ of codimension $p$. 
 
 If $C\subset L$, then $\RR^n+i\overline C$, the support of $\Ec_C$, is contained in $\RR^n+iL$, the source of $j$. 
 Therefore $j^!\Ec_C=j^*\Ec_C$ and so 
 \[
 k^* j^!\Ec_C\,\, =\,\, k^*j^*\Ec_C \,\,=\,\, (\Ec_C)|_{L_\CC}.
 \]
 Recall that $\Ec_C= \overline p_C^*\Ec_C^{\on{red}}$, where $\Ec_C^{\on{red}}$ is a cellular sheaf on
 $\RR^n/\LL(C)$ given by the spaces $E_K, K\geq C$ from $Q$ and their $\gamma$-maps.
 So the restriction $(\Ec_C)|_{L_\CC}$ is a similar pullback of the restriction of $\Ec_C^{\on{red}}$
 to $L/\LL(C)$, which is given by the $E_K$ for $C\leq K\subset L$, so the proposition is true in the case $C\subset L$.
 
 If $C\not\subset L$, then $C\cap L=\emptyset$, and so
 \[
 \xymatrix{
 \emptyset\ar[r]\ar[d]& \RR^n+iC
 \ar[d]^{j_C}
 \\
 \RR^n+iL \ar[r]_j&\CC^n
 }
 \]
 is a Cartesian square. Therefore, $j^!\Ec_C = j^! j_{C*}(j^!\Fc)=0$ by
  base change \cite[Prop. 3.1.9]{KS}.
 
  We conclude that $k^* j^! \Ec_C^\bullet[d]$ is an $\Sc^{(0)}$-constructible complex given by the double quiver $Q^{\leq L}$.
 Since $Q^{\leq L}$ consists of single vector spaces (not just complexes), Proposition 
 \ref{prop:perversity-char} implies that $k^*j^!\Fc[d]$ is perverse.\qed
 
 \vskip .3cm
 
 \noindent {\bf B. Transversal slice.} Let $L\in\Lc$ be as before. The normal bundle $N_{L/\RR^n}$ is canonically
 trivialized, with fiber $\RR^n/L$. Recall that $\RR^n/L$ carries the quotient arrangement $\Hc/L$. Choose
 a face $C\subset L$ open in $L$. Then the projection $\overline \pi_C: \Cc^{\geq C}\to \Cc_{\Hc/L}$
 identifies $\Cc_{\Hc/L}$ with the subposet $\Cc^{\geq C}\subset\Cc$. This leads to another 
 restriction functor
 \[
 \Rep^{(2)}(\Cc) \lra \Rep^{(2)} (\Cc^{\geq C}), \quad Q\mapsto Q^{\geq C},
 \]
 defined similarly to that in \S A.

 As before, suppose given $\Fc\in\Perv(\CC^n, \Hc)$ with double quiver $Q$. We have 
 then the {\em Verdier specialization} $\on{Sp}_{L_\CC}(\Fc)$, see \cite{KS} for background.
 It is a perverse sheaf on (the total space of) the normal bundle $N_{L_\CC / \CC^n} = L_\CC \times (\CC^n/L_\CC)$. 
 Choose some point $c\in C$ (this choice will be immaterial) and let $N_c\simeq \CC^n/L_\CC$
 be the fiber of the above normal bundle over $c$. We can think of $N_c$ as a transversal slice to $L_\CC$
 at $c$. Note that $N_c$ is  transverse to the characteristic variety of $\on{Sp}_{L_\CC}(\Fc)$, and so
 we have the perverse sheaf
 \[
 \Fc|_{N_c} \,\, := \,\, \on{Sp}_{L_\CC}(\Fc)|_{N_c} \,\, \in \,\, \Perv(\CC^n/L_\CC, \Hc/L). 
 \]
 
 \begin{prop}
 The double quiver of $\Fc|_{N_c}$ is identified with $Q^{\geq C}$. 
 \end{prop}
 
 \noindent {\sl Proof:} Let $T\subset \RR^n$ be an affine subspace forming
 a transversal slice to $L$ at $c$. We consider $T$ as an $\RR$-vector space with 
 origin $c$. Then the composition $T\hookrightarrow \RR^n\to \RR^n/L$ is an
 isomorphism of vector spaces. To understand the specialization 
 $\on{Sp}_{L_\CC}(\Fc)$ and the perverse sheaf $\Fc|_{N_c}$, we use the Cousin resolution
 $\Ec^\bullet = \Ec^\bullet(\Fc)$ and analyze $\on{Sp}_{L_\CC}(\Ec_K)$ for each summand $\Ec_K$ of
 each term of $\Ec^\bullet$. We claim that
 \be\label{eq:specializ}
 \on{Sp}_{L_\CC}(\Ec_K) \,\,=\,\, \begin{cases}
 \Ec_{\overline \pi_C(K)} (\on{Sp}_{L_\CC}(\Fc)),& \text {if } K\geq C, \\
 0,& \text{otherwise}. 
 \end{cases}
 \ee
 Indeed, the small neighborhood of the origin in $T_\CC$ meets only those closures of tube cells $\RR^n+i\overline K$,
 for which $K\geq C$. If $K\geq C$, the statement follows from the fact that $T_\CC$
 is transversal to the characteristic varieties of all the sheaves involved and therefore taking cohomology
 with support commutes with restriction to $T_\CC$. 
 
 Now, the statement about the double quiver of $\Fc|_{N_c}$ follows from \eqref{eq:specializ}
 immediately. \qed

\section{The double quiver determines a perverse sheaf}\label{sec:determines}

The goal of this section is to prove the following preliminary result.

\begin{thm}\label{thm:faithful}
The functor $\Qc: \Perv(\CC^n, \Hc)\to \on{Rep}^{(2)}(\Cc)$
from \eqref{eq:functor-Q}
 is fully faithful. 
\end{thm}

\noindent {\bf A. Orthogonality relations.} 
For an abelian category $\Ac$ let $C^b\Ac$ be the abelian category
formed by bounded complexes over $\Ac$ and morphisms of complexes
(not homotopy classes) in the usual sense. We start with

\begin{prop}\label{prop:cousin-faithful}
The Cousin resolution functor
\[
\Ec^\bullet:  \Perv(\CC^n, \Hc)\lra C^b\Sh_{\CC^n}, \quad \Fc \mapsto \Ec^\bullet(\Fc), 
\]
 is fully faithful. 
\end{prop}

\noindent {\sl Proof:} 
Let $\Ic$ be the  image of the functor $\Ec^\bullet$, i.e., the full subcategory in $C^b\Sh_{\CC^n}$ consisting
of complexes of sheaves  of the form $\Ec^\bullet(\Fc)$ for
$\Fc\in\Perv(\CC^n, \Hc)$. Thus we need to show that
$\Ec^\bullet: \Perv(\CC^n, \Hc)\to\Ic$ is an equivalence of categories. 
Define the functor $\Xi: \Ic \to \Perv(\CC^n, \Hc)$ to fit into the commutative diagram of functors
\[
\xymatrix{
\Ic \ar[d]_{\on{emb}'} \ar[r]^{\hskip -1cm \Xi} & \Perv(\CC^n, \Hc) \ar[d]^{\on{emb}''}
\\
C^b\Sh_{\CC^n} \ar[r]_{\on{can}}& D^b\Sh_{\CC^n}.
}
\]
Here $\on{emb}'$ is the embedding of $\Ic$ into the abelian category of complexes, $\on{emb}''$
is the embedding of $\Perv(\CC^n, \Hc)$ into the derived category, and $\on{can}$ is the canonical functor
from the abelian category of complexes to the derived category. Thus $\Xi(\Gc)$  for $\Gc\in\Ic$
is defined as the image of $\Gc$ is the derived category.   We know this image to be a perverse sheaf, i.e.,
an object of  $\Perv(\CC^n, \Hc)$.  Indeed, since $\Gc\in\Ic$,  we have that $\Gc = \Ec^\bullet(\Fc)$
for some $\Fc\in\Perv(\CC^n, \Hc)$, and we know that $\Ec^\bullet(\Fc)$ is quasi-isomorphic to $\Fc$ by
Corollary \ref{cor:cousin-complex}. 

We now prove that the functors $\xymatrix{    \Perv(\CC^n, \Hc) \ar@<.5ex>[r]^{\hskip 1cm \Ec^\bullet} & \Ic 
 \ar@<.5ex>[l]^{\hskip 1cm \Xi}}$ are quasi-inverse to each other. First for $\Fc\in\Perv(\CC^n, \Hc)$ we notice
 that $\Xi(\Ec^\bullet(\Gc)$ is canonically identified with $\Fc$ by Corollary \ref{cor:cousin-complex}.
 Further, for $\Gc\in\Ic$ we prove that $\Ec^\bullet(\Xi(\Gc))$ is canonically isomorphic to $\Gc$ in the abelian
 category of complexes. 
 This is an immediate consequence of the following ``orthogonality relations".

\begin{lem}\label{lem:orthogonality}
 Let $\Fc\in \Perv(\CC^n, \Hc_\CC)$. 
For any $C, C'\in\Cc$ we have
\[
\underline{R\Gamma}_{\RR^n + iC'}\Ec_C(\Fc) = 
\begin{cases}
\Ec_C(\Fc), & \mbox{if }C'=C;\\
0, & \mbox{if } C'\neq C.
\end{cases}
\]
\end{lem}

 \noindent {\sl Proof of the lemma:} If $C'=C$, the statement of the lemma
  is obvious. If $C'\neq C$, then $C'\cap C=\emptyset$ and so
  \[
  \xymatrix{
  \emptyset \ar[r]\ar[d] & \RR^n+iC\ar[d]^{j_C}
  \\
  \RR^n+iC'\ar[r]_{j_{C'}}&\CC^n
  }
  \]
  is a Cartesian square. So base change \cite[Prop. 3.1.9]{KS} implies that $j_{C'}^! j_{C*}\Gc=0$ for any $\Gc\in D^b\Sh_{\RR^n+iC}$.
  In particular,
  \[
  \underline{R\Gamma}_{\RR^n+iC'}\Ec_C(\Fc) \,\,=\,\, j_{C'*} j_{C'}^! j_{C*} j_C^!\Fc \,\,=\,\, 0. 
  \]
 Lemma   \ref{lem:orthogonality} and Proposition \ref{prop:cousin-faithful} are proved. \qed
  
  \vskip .3cm
  
  \noindent {\bf B. Recovery of $\Ec^\bullet(\Fc)$ from $\Qc(\Fc)$.}
  We prove
   Theorem \ref{thm:faithful} by the argument similar to that in the proof of Proposition \ref {prop:cousin-faithful}.
   That is, as a first step,  we explain how to 
     ``recover" the entire complex
  $\Ec^\bullet(\Fc)$ from the double quiver $\Qc(\Fc)$. Next, the second step is to consider the image $\Jc\subset \Rep^{(2)} (\Cc)$
  of the functor $\Qc$ and to interpret the ``recovery" procedure  by constructing a functor $\Theta: \Jc \to \Perv(\CC^n, \Hc)$
  quasi-inverse to $\Qc: \Perv(\CC^n,\Hc)\to\Jc$. 
  
  \vskip .2cm
  
  The first step proceeds as follows. 
  By Proposition \ref{prop:E-C-D-lift},
  each sheaf $\Ec_C(\Fc)$ is recovered from the data of $E_D(\Fc)$, $D\geq C$
  and of $\gamma_{D'D}$, $C\leq D'\leq D$. So the only remaining data are the
 matrix elements of the differentials. Eliminating the orientation torsors
 as in \eqref{eq:delta-cooriented}, we write these matrix elements as morphisms of sheaves
 \be
\underline \delta_{CC'}: \Ec_C(\Fc) \lra \Ec_{C'}(\Fc), \quad C'<_1 C.
 \ee
 Similarly to Proposition \ref{prop:delta-commute}, we have the following fact.
 
 \begin{prop}
 The morphisms of sheaves $\underline\delta_{CC'}$  commute, i.e., extend to 
  a contravariant representation of $(\Cc, \leq)$ in $\Sh_{\CC^n}$. 
 \qed
 \end{prop}
 
 In other words, we have a well-defined map $\underline\delta_{CC'}$ for any inclusion $C'<C$,
 not necessarily of codimension 1 obtained by composing the ``elementary" maps corresponding
 to codimension 1 inclusions.

  The data contained in $\underline\delta_{CC'}$ are precisely the induced morphisms on
 stalks of the sheaves $\Ec_C, \Ec_{C'}$  over all the cells $D+iC', D\in \Cc$, which are linear maps
 \be
 \delta_{CC'|D}: E_{C\circ D}\lra E_{C'\circ D}. 
 \ee
 So it is enough to express each $\delta_{CC'|D}$ through the $\gamma$ and $\delta$ maps of the
 double quiver $Q=\Qc(\Ec)$.
 We start with the following statement. 
 
 \begin{prop}\label{prop:gammadelta=id}
 For each $C'\leq C$ we have $\gamma_{C'C}\circ\delta_{CC'}=\on{Id}_{E_{C}}$. 
 In particular, each $\gamma_{C'C}$ is surjective, each $\delta_{CC'}$ is injective and
 $\dim(E_{C'})\geq\dim(E_C)$. 
 \end{prop}
 
 \noindent {\sl Proof:} It is enough to prove the statement for $C' <_1 C$, which we assume. 
   We use the fact that the maps of stalks induced by the morphism of sheaves
 $\underline\delta_{CC'}$ commute with the generalization maps from $C'+iC'$ to $C+iC'$. 
 This translates into the commutativity of
 \[
 \xymatrix{
 \Ec_C|_{C'+iC'} = E_C
 \ar[d]_{\on{Id}}
 \ar[rr]^{\delta_{CC'}=\delta_{CC'|C'}}&& E_{C'} = \Ec_{C'}|_{C'+iC'}
 \ar[d]^{\gamma_{C'C}}
 \\
 \Ec_C|_{C+iC'} = E_C \ar[rr]_{\delta_{CC'|C}} && E_C = \Ec_{C'}|_{C+iC'}.
 }
 \]
  Here the vertical $\on{Id}$ is the generalization map for $\Ec_C$ and $\gamma_{C'C}$
 is the generalization map for $\Ec_{C'}$. 
 Our statement follows therefore from the next lemma.
 
 \begin{lem}\label{lem:delta=id}
 The map $\delta_{CC'|C}$ is equal to $\on{Id}_{E_C}$. 
  \end{lem}
  
  \noindent {\sl Proof of the lemma:} 
 Recall the identifications of stalks
  \[
  \Ec_{C}|_{C+iC'}\buildrel\alpha\over \lra E_C, \quad \Ec_{C'}|_{C+iC'} 
  \buildrel\alpha'\over \lra E_C
  \]
 given in the proof of Proposition  \ref{prop:E-C-D-lift}.
   Let $d=\on{codim}(C)$. 
 In the notation of   \ref{prop:E-C-D-lift}, we take $K=C$ and form the complex $\Gc = \ul{R\Gamma}_{\RR^n+i\LL(C)}(\Fc)[d]$,
 so we have:
 \begin{itemize}
 \item[(1)] The restriction of  $\Gc$  on $C+i\LL(C)$ reduces to a  constant  sheaf in degree 0, denote it $\Nc$.
 
 \item[(2)] $\Ec_C |_{C+iC'} = \Gc|_{C+iC'}$, while  $E_C = \ul{R\Gamma}_{\RR^n+i0}(\Gc)[n-d]|_{C+i0}$, and
 $\alpha$ is induced by the local Poincar\'e duality \eqref{eq:LPD} for $\Nc$. 
 \end{itemize}
 
 \noindent In the same way, $\alpha'$ is induced by the local Poincar\'e duality for the constant sheaf $\Nc'$ obtained as
 the restriction to $C+i\LL(C')$ of
  \[
  \Gc' \,\,=\,\,  \underline{R\Gamma}_{\RR^n+i \LL(C')}(\Fc) [d+1] \,\,=\,\, 
  \ul{R\Gamma}_{\RR^n+i\LL(C')}(\Gc)[1].
  \]
  The last equality above means that
  \[
  \Nc' = \ul H^1_{C+i \LL(C')}(\Nc), 
  \]
  which we interpret, by codimension 1 local Poincar\'e duality,  as
  \be\label{eq:ororr}
  \Nc'\,\, \simeq \,\, k^*\Nc \otimes \orr(C'/C) \buildrel  \sim \over\lra  k^*\Nc. 
  \ee
  Here $k: C+i\LL(C')\hookrightarrow C+i\LL(C)$ is the embedding, and  the last isomorphism  comes from the identification of
  orientation torsors given by the codimension 1 inclusion $C' <_1 C$. 
  Now,  $\underline\delta_{CC'}$ is the coboundary map
 \[
\ul \HH^d_{\RR^n+iC}(\Fc) \lra \ul\HH^{d+1}_{\RR^n+i C'}(\Fc). 
 \]
 This means that, under our identifications, its stalk at $C+iC'$ becomes equal to the stalk, also at $C+iC'$, of the coboundary map
 for $\Nc$ which we write as:
 \[
  j_* j^*\Nc \buildrel \ul\delta\over\lra \ul H^1_{C+iC'}(\Nc) = (j')^* \Nc' \buildrel \eqref{eq:ororr}\over = l^*\Nc,
  \]
 Here
\[ 
 j: C+iC \hookrightarrow C+i\LL(C), \quad  j': C+iC'\hookrightarrow C+i\LL(C'), \quad l: C+iC'\hookrightarrow C+iC
  \]
  are the embeddings. So our statement reduces to the following elementary fact (``codimension  1 Poincar\'e duality is
  given by the coboundary map").  
  
   \begin{lem}
   Let $\epsilon: M_0\to M$ be a  closed codimension 1 embedding of $C^\infty$-manifolds, 
    $J: M_+\to M$ an open embedding such that the closure $\overline {J(M_+)}$ is a manifold with boundary $M_0$
    (and so gives a trivialization of the orientation torsor $\orr_{M_0/M}$). For any locally constant sheaf $\Kc$ on $M$
    the coboundary map
    \[
   \epsilon^* J_* J^* \Kc\buildrel\ul\delta\over\lra \ul H^1_{M_0}(\Kc)
    \]
    corresponds, after the identification $  \epsilon^* J_* J^* \Kc = \epsilon^*\Kc$ and the Poincar\'e duality
    $\ul H^1_{M_0}(\Kc) \simeq \epsilon^*\Kc$, to the identity of $\epsilon^*\Kc$. \qed
   \end{lem}
   This finishes the proof of Lemma \ref{lem:delta=id} and Proposition \ref{prop:gammadelta=id}.

 \vskip .2cm
 
 We now let
 $C, C', D$   be three  arbitrary faces such that $C'\leq C$. 
 Put 
 \be\label{eq:lifts-of-C-C'}
K= C\circ D\geq C, \quad K'= C'\circ D\geq C'. 
\ee
By the associativity of the operation $\circ$, we have $C\circ K'=K$ and $C'\circ K'=K'$. 
Note that because $\Ec_C, \Ec_{C'}$ are pullbacks of sheaves on $\RR^n/\LL(C)$,
resp. $\RR^n/\LL(C')$, we have
\be\label{eq:delta-c-c'-delta'}
\delta_{CC'|D} = \delta_{CC'|K'}.
\ee
That is, 
\[
\delta_{CC'|D}: E_{C\circ D} = E_K \lra C_{C'\circ D}=E_{K'}
\]
is equal to
\[
\delta_{CC'|K'}: E_{C\circ K'} = E_K \lra E_{C'\circ K'}=E_{K'}. 
\]
To complete the recovery procedure of $\Fc$ from $\Qc(\Fc)$, we prove:

  \begin{prop}\label{prop:short-relation}
 In the described situation, we have
 \[
 \delta_{CC'|D}\,\, = \,\,  \gamma_{C'K'}\, \delta_{KC'}: E_K\lra E_{K'}.
\]
 \end{prop}

 \noindent {\sl Proof:} We first use that 
 the maps of stalks induced by the morphism of sheaves
 $\underline\delta_{CC'}$ commute with the generalization maps from $C'+iC'$ to $K'+iC'$.
 This gives the commutativity of 
 \[
 \xymatrix{
 E_C 
 \ar[d]_{\gamma_{CK}}
 \ar[r]^{\delta_{CC'}}& E_{C'}
 \ar[d]^{\gamma_{C'K'}}
 \\
 E_K \ar[r]_{\delta_{CC'|K'}}& E_{K'},
 }
 \]
 i.e., the equality
 \[
\delta_{CC'|K'}\, \gamma_{CK}  =\gamma_{C'K'}\, \delta_{CC'}. 
 \]
 We precompose this equality with $\delta_{KC}$:
  \[
  \delta_{CC'|K'}\, \gamma_{CK}\, \delta_{KC} = \gamma_{C'K'}\, \delta_{CC'}
  \, \delta_{KC}.
  \]
  Now, using on the left, Proposition 
 \ref{prop:gammadelta=id}, and on the right, the fact that 
  the $\delta$ maps form a contravariant representation of $(\Cc, \leq)$, 
  and  also invoking \eqref{eq:delta-c-c'-delta'}, 
  we get the desired statement. \qed
  
  \begin{cor}
  In the above situation, we also have
  \[
  \delta_{CC'|D} = \gamma_{0K'}\, \delta_{K0}. 
  \]
  \end{cor}
  
  \vskip .3cm

  \noindent {\bf C. End of proof of Theorem \ref{thm:faithful}.} We now perform the second step
  outlined at the beginning of \S B by interpreting the above in a more categorical language. 
  So we denote by $\Jc\subset\Rep^{(2)}(\Cc)$ the image of the functor $\Qc$
  and construct a functor $\Theta$ as in the diagram  
   \[
\xymatrix{
  \Perv(\CC^n, \Hc) \ar@<.5ex>[r]^{\hskip .7cm \Qc} & \Jc
  \ar@<.5ex>[l]^{\hskip 1cm \Theta} 
  }
  \]
  so that the two functors are quasi-inverse. Explicitly, let 
  $Q=(E_C, \gamma_{C'C}, \delta_{CC'})$ be a double quiver from $\Jc$. 
   For each $C\in\Cc$ we {\em define} the sheaf $\Ec_C(Q)$ by postulating the formulas of Corollary 
   \ref {cor:stalks-E-C}(b), i.e., by  setting
  
\be\label{eq:stalks-Q}
\begin{gathered}
 \Ec_C(Q)\bigl|_{iC_1+D} = \begin{cases}
  E_{C\circ D}, & \text{if } C_1\leq C; \\
  0, & \text{otherwise};
  \end{cases}
  \\
  \gamma^{\Ec_C(Q)}_{iC'_1+D', iC_1+D} = \gamma_{C\circ D', C\circ D}\otimes\on{Id}, \quad 
  C'_1\leq C_1\leq C, \,\, D'\leq D.
  \end{gathered}
  \ee
  Next, for any $C' <_1  C$ we {\em define} a morphism of sheaves 
  \[
  \underline\delta_{CC'}: \Ec_C(Q) \lra \Ec_{C'}(Q)
  \]
  by postulating the formulas of Proposition \ref{prop:short-relation}, i.e., by defining the action of on the stalk
  over $D+iC'$, $D\in\Cc$, to be
 \be\label{eq:delta-Q}
 \delta_{CC'|D} \,\, =\,\,  \gamma_{C'K'} \, \delta_{KC'}:  E_{C\circ D} =  E_K\lra E_{K'} =  E_{C'\circ D}, 
\ee
where $K$ and $K'$ are defined by \eqref{eq:lifts-of-C-C'}.  Since we know that $Q\in\Jc$, i.e., 
$Q=\Qc(\Fc)$ for some $\Fc\in\Perv(\CC^n, \Hc)$, the results of n$^\circ$ B imply that these maps of stalks commute
with generalization maps and so indeed define morphisms of sheaves $\underline\delta_{CC'}$. Further, for the same reason
($Q\in\Jc$),  these
morphisms commute  and so assemble into  a complex of sheaves
\[
\Theta(Q) \,\,=\,\,\biggl\{ \bigoplus_{\on{codim}(C)=0} \Ec_C(Q)\otimes\OR(C) \buildrel \underline\delta \over\to 
\bigoplus_{\on{codim}(C)=1} \Ec_C(Q)\otimes\OR(C)\buildrel \underline\delta \over\to\cdots \buildrel \underline\delta \over\to
\Ec_0(Q)
\biggr\}, 
\]
which, moreover, lies in $\Perv(\CC^n, \Hc)$. This defines the functor $\Theta$, and the results of 
 n$^\circ$ B mean that $\Qc\Theta$ is naturally isomorphic to the identity functor of $\Jc$, while
 $\Theta \Qc$ is naturally isomorphic to the identity functor of  $\Perv(\CC^n, \Hc)$. \qed

\section{Algebraic relations in the double quiver }

\noindent{\bf A. The transitivity relations.}
Let
\be
Q = \bigl( (E_C)_{C\in \Cc}, (\gamma_{C'C}, \delta_{CC'})_{C'\leq C}\bigr)\,\,\in\,\, \Rep^{(2)}(\Cc)
   \ee
 be a  double representation of $\Cc$.
 In this    section we find algebraic relations 
 among the maps $\gamma_{C' C}, \delta_{CC'}$ which are necessary for $Q$ to have the form $Q=\Qc(\Fc)$ for
 some   $\Fc\in\Perv(\CC^n, \Hc)$. 
 
 \vskip .2cm

 Call $Q$ {\em monotone}, if $\gamma_{C'C}\, \delta_{CC'}=\on{Id}$ for any $C'\leq C$. 
 By  Proposition \ref{prop:gammadelta=id}, any $\Qc(\Fc)$ is monotone. 
 
  Given a monotone $Q$, for any $A,B\in \Cc$ we define
  the {\em transition map}
  \be
  \phi_{AB} = \gamma_{MB}\,  \delta_{AM}: E_A\lra E_{B}, \quad  M\leq A,B. 
  \ee
  By monotonicity, the choice of $M$ is immaterial, for example one can
  take $M=0$. Note that $\phi_{AB}$ is equal to $\gamma_{AB}$, if $A\leq B$,
  and to $\delta_{AB}$, if $A\geq B$.

  \begin{Defi} 
  Three faces $A,B,C\in\Cc$ are called {\em collinear}, if there are $a\in A, b\in B, c\in C$
  such that $b\in[a,c]$, i.e., $b$ lies on the straight line segment between $a$ and $c$. 
  \end{Defi}
  
  Collinearity is recovered from the oriented matroid corresponding to $\Hc$. More precisely:
  
  \begin{prop}
  Let us introduce a total order $\preceq$ on $\{+,-,0\}$ which is induced by the standard order on $\RR$,
  i.e., by $-\preceq 0 \preceq +$. Then the following are equivalent:
  \begin{enumerate}
  \item[(i)] The faces $A,B,C$ are collinear.
  
  \item[(ii)] For each $H\in\Hc$ the sequence of sign vectors $(A_H, B_H, C_H)$ is monotone increasing
  or decreasing: either $A_H\preceq B_H \preceq C_H$, or  $C_H\preceq B_H \preceq A_H$.
  \end{enumerate}
  
  \end{prop}
  
  \noindent {\sl Proof:} (i)$\Rightarrow$(ii) is obvious. To see the converse, suppose that three cells $A,B,C$ are
  {\em not} collinear. Then  there is a hyperplane $H\in\Hc$ such that $B$ lies on one side of $H$ and
  $A,C$ lie on the other side. This contradicts the monotonicity of $(A_H,B_H,C_H)$. \qed

  \begin{thm}[Transitivity relations]
  \label{thm:transitivity}
  Let $Q=\Qc(\Fc)$ for some $\Fc\in\Perv(\CC^n, \Hc)$. Then for any collinear faces
  $A,B,C$ we have
  \[
  \phi_{AC}=\phi_{BC}\, \phi_{AB}: E_A\lra E_C. 
  \]
  
  \end{thm}
  The proof will be given after Example \ref{ex:zifferblatt}.

  \begin{rem}[(Transitivity: long form)]
  \label{rem:longform}
  By marking all the faces meeting $[a,c]$, we get 
  a face path (alternating sequence of
  inclusions)
  \[
  A = B_1   \geq B'_1 < B_2 <  B'_2 >  \cdots  < B_{m-1}> B'_{m-1} \leq  B_m=C
  \]
  (there are 4 possibilities as to whether $A=B'_1$ or $B'_m=C$). 
  By iterating Theorem \ref{thm:transitivity}, we can reformulate in in the
  equivalent form:
  \[
  \phi_{AC}= \gamma_{B'_{m-1} B_m} \,  \delta_{B_{m-1} B'_{m-1}} \,  \cdots 
  \,  \delta_{B'_1 B_2}\,  \delta_{B_1 B'_1}. 
  \]
  Note that it may be possible that $A$ and $C$ can be connected by a straight line segment in
  more than one inequivalent way,  in which case $\phi_{AC}$ can be
  expressible through the $\gamma$ and $\delta$ maps in more than one way,
  producing additional algebraic relations in $\Qc$, see Example \ref{ex:zifferblatt} below. 
  
  \end{rem}

  \begin{ex}[(Base change)]
  \label{ex:base-change}
  It is sometimes convenient to view the poset $\Cc$ as a category, that is, to write an inclusion
  $C'\leq C$ as a morphism $u: C'\to C$. Then a double representation $Q$ can be viewed
  as a ``bivariant theory" on $\Cc$: for a morphism $u$ as above, we write 
  \[
  u_*=\gamma_{C'C}, \,\, u^*=\delta_{CC'}. 
  \]
  In this language,  a simplest instance of Theorem \ref{thm:transitivity} can be viewed as a ``base change
  property": we consider a  (necessarily commutative) square of face inclusions
  \be\label{eq;cartesian}
  \xymatrix{
  D
  \ar[d]_{v_2}
  \ar[r]^{v_1}& A
  \ar[d]^{u_1}
  \\
  C\ar[r]_{u_2} & B.
  }
  \ee
  The condition that $A,B,C$ are collinear, means that the square is coCartesian in the
  categorical sense:  $B$ is the minimal face containing $A$ and $C$ in its closure.
  In this case, Theorem \ref{thm:transitivity} says that
  \[
  v_{2*} v_1^* = u_2^* u_{1*}: E_{A} \lra E_{C}.
  \]
  Indeed, the LHS is $\phi_{AC}$, while $u_{1*}= \gamma_{AB}=\phi_{AB}$ and $u_2^*=\delta_{BC}=\phi_{BC}$. 
  \end{ex}
  
     \begin{ex}[(Zifferblatt relations)]\label{ex:zifferblatt}
     Another extreme instance of Theorem \ref{thm:transitivity}
     corresponds to the case  when $C=-A$ is the opposite cell to $A$.  Suppose $\dim(A)\geq 2$,
     and let $L$ be a 2-dimensional subspace such that $\dim(L\cap A)=2$. 
     The arrangement $\Hc\cap L$ then cuts $L\setminus\{0\}$ into some even number $2m$ of 2-dimensional
     open cones, which we number cyclically $\overline B_1, \cdots, \overline B_{2m}$
     and the same number of 1-dimensional  open  rays $\overline B'_1, \cdots , \overline B'_{2m}$. 
     Let $B_\nu, B'_\nu$ be the faces of $\Cc$ which intersect $L$ in $\overline B_\nu, \overline B'_\nu$,
     see Fig.\ref{fig:Zifferblatt}. 
     
      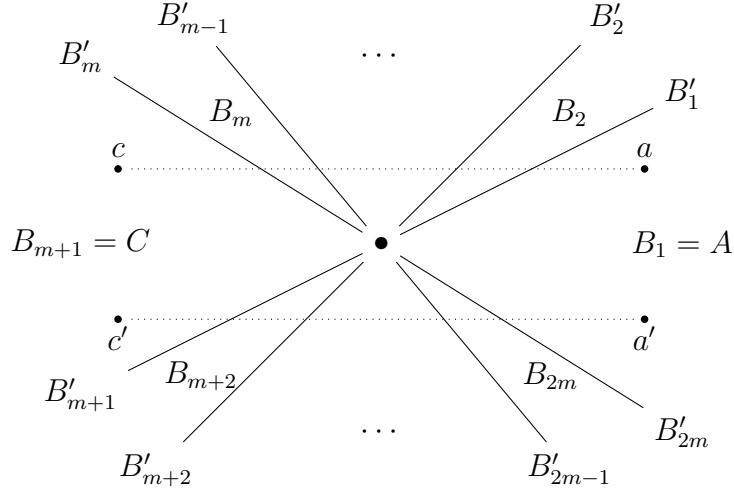
\begin{figure}
      \centering
    \begin{tikzpicture}[scale=0.5]
      \node (0) at (0,0){$\bullet$}; 
      \node(C'1) at (8,4){$B'_1$}; 
      \draw (0) --( C'1); 
      \node (C'2) at (6,6){$B'_2$};
      \draw (0) -- (C'2);
      \node (C'm-1) at (-8,5){$B'_{m}$};
      \draw (0) -- (C'm-1); 
      \node (C'm) at (-8,-4){$B'_{m+1}$};
      \draw (0) -- (C'm);
      \node (C'm+1) at (-6,-6){$B'_{m+2}$};
      \draw (0) -- (C'm+1);
      \node (C'2m-1) at (5,-6){$B'_{2m-1}$};
      \draw (0) -- (C'2m-1);
      \node (C'2m) at (8,-5){$B'_{2m}$};
      \draw (0) -- (C'2m);
      
    \node at (0,5) {$\cdots$}; 
     \node at (0,-5) {$\cdots$}; 
    
    \node (C'm-2) at (-5,6){$B'_{m-1}$}; 
    \draw (0) -- (C'm-2); 
    
    \node at (8,0){$B_1=A$}; 
    \node at (-8,0) {$B_{m+1}=C$}; 
    \node (a) at (7,2){}; 
     \fill (a) circle (0.1); 
     \node (c) at (-7,2){};
      \fill (c) circle (0.1); 
\draw[dotted] (a) -- (c); 
\node at (7,2.5){$a$};
  \node at (-7,2.5){$c$};   
  
   \node (a') at (7,-2){}; 
     \fill (a') circle (0.1); 
     \node (c') at (-7,-2){};
      \fill (c') circle (0.1); 
\draw[dotted] (a') -- (c'); 
\node at (7,-2.5){$a'$};
  \node at (-7,-2.5){$c'$};   
  
  \node at (5,3.5){$B_2$};
  
  \node at (-4,3.5){$B_{m}$};
  
   \node at (-4.8,-3.5){$B_{m+2}$};
   
    \node at (4.5,-3.5){$B_{2m}$};
    \end{tikzpicture}
     \caption{The Zifferblatt. }\label{fig:Zifferblatt}

    \end{figure}

    Note that there are two inequivalent ways to join a point of $A$ with a point of $C$ by
    a straight line segment inside $L$ not passing through 0, represented by the segments $[a,c]$ and $[a',c']$
    in Fig. \ref{fig:Zifferblatt}. 
So the long form of the transitivity relations (Remark \ref {rem:longform}) gives the
{\em Zifferblatt relation} (we borrow the term from \cite{manin-schechtman}):
\[
\begin{gathered}
\gamma_{B'_{m}, C}\,  \delta_{B_{m},B'_{m}} \,  \gamma_{B'_{m-1}, B_{m} }\,  \cdots \,  \delta_{B_2 ,B'_2}
\,  \gamma_{B'_1, B_2} \,  \delta_{A,B'_1} = \\
 = \gamma_{B'_{m+1}, C}\,  \delta_{B_{m+2}, B'_{m+1}} \,  \gamma_{B'_{m+2} ,B_{m+2}}\,  \cdots
 \,  \delta_{B_{2m}, B'_{2m-1}}\,  \gamma_{B'_{2m}, B_{2m}}\,  \delta_{A, B'_{2m}}
\end{gathered} 
\]
(both sides of this equality are equal to $\phi_{A,C}=\gamma_{0,C}\circ\delta_{A,0}$). 

\end{ex}
  
  \vskip .3cm       
     
 \noindent {\bf B. Proof of Theorem \ref{thm:transitivity}. }

     \vskip .2cm
     
    \noindent {\sl Step 1: Base change.}
   We first consider the situation of a coCartesian square  from
  Example \ref{ex:base-change}, as it will serve as an inductive step
  in  treating more general cases.
    The assumption that $A,B,C$ are collinear implies that
  $C\circ A=B$.
  
  \begin{lem}
  In the situation of  Example \ref{ex:base-change}, the map 
  \[
  \delta_{AD|C}: E_{C\circ A}=E_B \lra E_C = E_{C\circ D}
  \]
  is equal to $\delta_{BC}$. 
  \end{lem}
  
  \noindent {\sl Proof of the lemma:} We consider the commutative square of the $\Ec$-sheaves corresponding to 
  \eqref{eq;cartesian} and the corresponding commutative square of stalks over $C+iD$. These squares have the form
  \[
  \xymatrix{\Ec_D & \ar[l]_{\underline\delta_{AD}} \Ec_A
  \\
  \Ec_C
  \ar[u]^{\underline\delta_{CD}}
  & \ar[l]^{\underline\delta_{BC}}\Ec_B, 
  \ar[u]_{\underline\delta_{BA}}
  }
  \quad\quad 
  \xymatrix{
  E_C & \ar[l]_{\delta_{AD|C}} E_B
  \\
  E_C
  \ar[u]^{\on{Id}}
  &\ar[l]^{\delta_{BC}} E_B,
  \ar[u]_{\on{Id}}
  }
  \]
  whence the lemma.\qed
  
  To deduce our particular case of Theorem \ref{thm:transitivity} from the lemma, we 
  spell out the condition that the maps of stalks induced by $\underline\delta_{AD}$ commute
  with generalization from $D+iD$ to $C+iD$. This gives a commutative square of vector spaces
  \[
  \xymatrix{E_B \ar[rr]^{\delta_{AD|C}=\delta_{BC}}& &E_C
  \\
  E_A
  \ar[u]^{\gamma_{AB}}
   \ar[rr]_{\delta_{AD|D}=\delta_{AD}}& &E_D
   \ar[u]_{\gamma_{DC}}
  }
  \]
  in which the path through $E_D$ gives, as the composite map, $\phi_{AC}$ while the path through $E_B$
  consists of $\delta_{BC}=\phi_{BC}$ and $\gamma_{AB}=\phi_{AB}$.

  \vskip .2cm
  
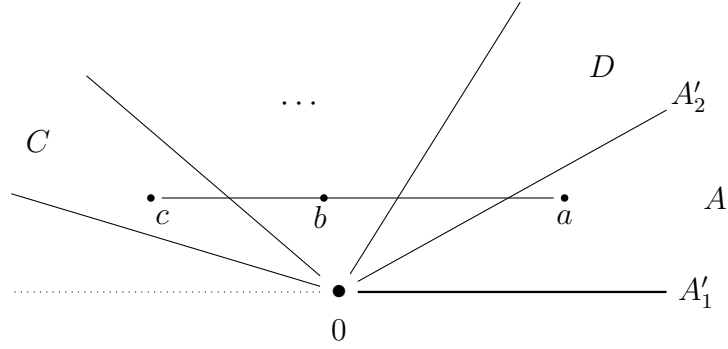
\begin{figure}
    \centering
    \begin{tikzpicture}[scale=0.5]
     
    \node (C') at (0,0){$\bullet$};

     \node (C) at (9.5,0){$ A'_1$};
     \node (preC) at (9,0){};
     \draw[line width=0.3mm] (C') -- (preC){};
     
     \node (pre-d'0) at (9,5){};
     
     \draw (C') -- (pre-d'0);
     
     \node (pred'1) at (5,8){}; 
     
     \draw (C') -- (pred'1); 
     
     \node (pre-d'm+1) at (-9,2.7){};
     
     \draw (C') -- (pre-d'm+1); 
     
     \node (pre-d'm) at (-7, 6){};
     
     \draw (C') -- (pre-d'm);
     
     \node at (0, -1){$0$}; 
     
     \node at (10, 2.5){$ A$};
     
     \node at (9.3,5.2) {$ A'_2$}; 
     
     \node at (7,6) {$D$};
     
     \node at (5.1,8.2) {$$};
     
     \node at (-7.2,6.3) {$$}; 
     
     \node at (-8, 4) {$C$};
     
     \node (left)  at (-9,0){};
     
     \draw[dotted] (C') -- (left);
     
     \node at (-1, 5) {$\cdots$}; 
     
     \node (x) at (6, 2.5) {};
     \fill (x) circle (0.1); 
     
     \node (x') at (-5, 2.5){}; 
     \fill (x') circle (0.1); 
     
      \node (b) at (-.4, 2.5){}; 
     \fill (b) circle (0.1); 
     
     \draw (x) -- (x'); 
     
     \node  at (6,2){$a$};
     \node at (-4.7,2) {$c$}; 
     \node at (-0.5,2) {$b$};

    \end{tikzpicture} 
    \caption{Case $C\neq -A$. }\label{fig:chain}
    \end{figure}

  \noindent {\sl Step 2:  Case when $C\neq -A$.} In this case the segment $[a,c]$ does not pass through 0
  and so its $\RR$-linear span is a 2-dimensional subspace $L\subset\RR^n$. As in Example 
  \ref{ex:zifferblatt}, we then have an induced arrangement $\Hc\cap L$ of lines in $L$, and there are
  various possibilities as to whether  the faces
  \[
  \overline A= A\cap L,\,\,
  \overline B= B\cap L, \,\,
  \overline C=C\cap L
  \]   have dimension 1 or 2. 
  We first consider the case $\dim(\overline A)=2$,
  depicted in 
   Fig. \ref{fig:chain}.
  Denote by $\overline A'_1, \overline A'_2$ the two rays bordering $\overline A$
   so that $\overline A'_2$ meets $[a,c]$. Further, denote $\overline D$ the next 2-dimensional face
   in the direction from $a$ to $c$. Let $D, A'_1, A'_2\in\Cc$ be the face whose intersections with $L$
   are $\overline D, \overline A'_1, \overline A'_2$. Our assumption that $C\neq -A$ implies
    that the intersection with $L$ of  all
   the faces meeting $[a,c]$ lie on the same side of the line $\LL(A'_1)$.

  Consider the square of inclusions $0\leq A'_1, A'_2\leq A$.
   Let us write the corresponding commutative squares of $\Ec$-sheaves and of stalks of these sheaves
   at $C= C+i0$. The latter square has the form
   \[
   \xymatrix{
   \Ec_A\bigl|_{C}= E_A
   \ar[d]_{\phi_{AD}}
    \ar[rr]^{\hskip -1.5cm \on{Id}}&& E_A = E_{A'_1\circ C}= \Ec_{A'_1}\bigl|_{C}
    \ar[d]^{\phi_{AC}}
   \\
   \Ec_{A'_2}\bigl|_{C}=E_{A'_2\circ C}=E_D\ar[rr]_{\hskip 1cm \phi_{DC}} &&
   \Ec_0\bigl|_{C}=E_C,
   }
   \]
   so $\phi_{AC}=\phi_{DC}\phi_{AD}$. Similarly, by considering the commutative square of stalks over $B+i0$,
   we find that $\phi_{AB}=\phi_{DB}\phi_{AD}$. Therefore we reduce to proving that $\phi_{DC}=\phi_{BC}\phi_{DB}$
   which can be further reduced in a similar way.  So  proceeding by induction, we reduce the situation to the
     case $B=A$ or  $B=A'_2$. If $B=A$, there is nothing to prove. If $B=A'_2$, then 
     \[
     \phi_{AC}=\gamma_{0C}\delta_{A0} = \gamma_{0C}\delta_{B0}\delta_{AB}= \phi_{BC}\delta_{AB} 
     =\phi_{BC}\phi_{AB}. 
     \]
     This concludes the treatment of the case $\dim(\overline A)=2$. 
     
     The case when $\dim(\overline A)=1$, is analyzed inductively in a similar way, reducing it to the
     situation of Example  \ref{ex:base-change} which has been analyzed in Step 1.
       This concludes the analysis of Step 2.
     
     \vskip .2cm
     
     \noindent {\sl Step 3: Case $C=-A$.} In this case we use the notation of Example 
     \ref{ex:zifferblatt}. We consider the segment $[a,c]$. Note that by Step 2 we have
     \[
     \phi_{B'_1 B'_{m}}\,\,  :=  \,\, \gamma_{0 B'_{m}}\circ\delta_{B'_1 0} \,\, = \,\,
     \delta_{B_{m} B'_{m}}\circ \gamma_{B'_{m-1} B_{m}} \circ \cdots 
     \circ \delta_{B_2 B'_2}\circ \gamma_{B'_1 B_2}. 
     \]
     Therefore, using the fact that the $\gamma$ and $\delta$ maps commute with compositions of inclusions, we have
     \[
     \begin{gathered}
     \phi_{A,C} \,\, := \,\,\gamma_{0,C}\, \delta_{A,0} \,\,=\,\,
     \gamma_{B'_{m}, C}\, \gamma_{0, B'_{m}} \, \delta _{B'_1, 0}\, \delta_{A, B'_1} = \\
     =  \gamma_{B'_{m}, C} \, \bigl(  \delta_{B_{m}, B'_{m}}\, \gamma_{B'_{m-2}, B_{m-1}} \, \cdots 
     \, \delta_{B_2, B'_2}\, \gamma_{B'_1, B_2}\bigr) \, \delta_{A, B'_1},
     \end{gathered}
     \]
     which is the long form of the transitivity relation. 
     The case of the segment $[a',c']$ is treated in the same way.
     
     This concludes the proof of Theorem \ref{thm:transitivity}.

\section{Equivalence of categories }

\noindent {\bf A. The main result.} 
Let $\Hc$ be an arrangement of hyperplanes in $\RR^n$, with the poset of faces $\Cc$,
and $\Rep^{(2)}(\Cc)$ be the corresponding category of double representations
$Q=(E_C, \gamma_{C'C}, \delta_{CC'})$. 
Let $\Ac=\Ac_\Hc\subset\Rep^{(2)}(\Cc)$ be the full subcategory of 
   double representations satisfying the following  three conditions:
   
   \begin{enumerate}
\item[(Mon)] Monotonicity: $\gamma_{C'C}\, \delta_{CC'}=\on{Id}_{E_C}$, $C'\leq C$. This allows us to define
transition maps $\phi_{AB}=\gamma_{CB}\, \delta_{AC}: E_A\to E_B$ where $C$ is an arbitrary face
$\leq A,B$.

\item[(Tran)] Transitivity: $\phi_{AC}=\phi_{BC}\, \phi_{AB}$ for any three collinear faces $A,B,C$. 

\item[(Inv)] Invertibility: Let $C_1, C_2$ be two faces of the same dimension $d$ with the same linear span
$\LL(C_1)=\LL(C_2)$, which lie on opposite sides of a face $D_1$ of dimension $d-1$, so
$C_1 >_1 D_1 <_1 C_2$. Then $\phi_{C_1C_2} = \gamma_{D_1C_2}\, \delta_{C_1 D_1}: E_{C_1}\to E_{C_2}$
is an isomorphism. 

\end{enumerate}
 
The following is the main result of this paper.

\begin{thm}\label{thm:main}
The functor $\Fc\mapsto \Qc(\Fc)$ defines an equivalence of categories $\Perv(\CC^n, \Hc)\to\Ac$. 
 
\end{thm}

In view of Theorem \ref{thm:faithful},  it suffices to prove the following statement.

\begin{refo}\label{refo:mainref} 
In order for $\Qc\in\Rep^{(2)}(\Cc)$ to have the form $\Qc(\Fc)$ for some $\Fc\in\Perv(\CC^n, \Hc)$,
it is necessary and sufficient that $\Qc$ satisfies (Mon), (Tran) and (Inv). 

\end{refo}

The proof will occupy the remainder of this section. 

\vskip .3cm

\noindent {\bf B. Necessity.} The necessity of (Mon) and (Tran) has already been proved
in Proposition \ref{prop:gammadelta=id} and Theorem \ref{thm:transitivity}. Let us prove the necessity of
(Inv). Because of Proposition \ref{prop:hyp-res}, it is enough to consider the case when $C_1$ and $C_2$
are  faces open in $\RR^n$, because the general case will then follow by considering the hyperbolic
restriction to $L_\CC$, where $L=\LL(C_1)=\LL(C_2)$. 

Assuming the $C_1$ and $C_2$ are open, consider the Cousin resolution $\Ec^\bullet$ of $\Fc$.
Note that $\Fc$ is $\Sc^{(0)}$-smooth, while the summands $\Ec_C$ of $\Ec^\bullet$ are only
$\Sc^{(1)}$-smooth. Look at the $\Sc^{(0)}$-smooth sheaf
\[
\underline H^0(\Fc) \,\,\,=\,\,\, \on{Ker} \,\biggl\{ \bigoplus_{C \text{ open}} \Ec_C \buildrel \widetilde{\underline\delta}
\over\lra \bigoplus_{\on{codim}(D)=1}\Ec_D\biggr\}.
\]
The $\Sc^{(1)}$-cells $[D_1, C_2]\leq [C_2, C_2]$ lie in the same stratum $(\CC^n)^\circ$ of
$\Sc^{(0)}$. Therefore the generalization map for $\underline H^0(\Fc)$ from $[D_1, C_2]$
to $[C_2, C_2]$ must be an isomorphism. But, applying Corollary 
\ref{cor:stalks-E-C}(c), we find that up to tensoring with $\OR(\RR^n)$, we have
\[
\begin{aligned}
\underline H^0(\Fc)|_{[D_1, C_2]} \,\,=\,\, \on{Ker} \bigl\{ E_{C_1}\oplus E_{C_2} \buildrel
\phi_{C_1 C_2}-\on{Id}\over\lra E_{C_2}\bigr\}, 
\\
\underline H^0(\Fc)_{[C_2, C_2]} \,\,=\,\, E_{C_2}, 
\end{aligned}
\]
and the generalization map is induced by the projection $E_{C_1}\oplus E_{C_2}\to E_{C_2}$. In order for this projection
to restrict to an isomorphism $\on{Ker}(\phi_{C_1 C_2}-\on{Id}) \to E_{C_2}$, the map
$\phi_{C_1 C_2}$ must be an isomorphism. 

\vskip .3cm

\noindent {\bf C. Sufficiency: construction of the complex $\Ec^\bullet$.} 
To prove the sufficiency of the three conditions in Reformulation \ref{refo:mainref}, we start with a double 
representation $Q$ satisfying them, construct a complex $\Ec^\bullet =\Ec^\bullet(Q)$
and then prove that it is a perverse sheaf with double quiver $Q$. 

In this procedure, the {\em construction} of the complex $\Ec^\bullet(Q)$ will rely  only on
the properties
(Mon) and (Tran), while (Inv) will be needed to ensure that this complex is an object of $\Perv(\CC^n, \Hc)$. 

The construction will be done by the same procedure as in \S  \ref{sec:determines}C. 
That is, we    define
 the  $\Sc^{(2)}$-smooth sheaves $\Ec_C=\Ec_C(Q)$  by the formulas
 \eqref{eq:stalks-Q}. In fact, it follows from the definition that each $\Ec_C(Q)$
 is $\Sc^{(1)}$-smooth. We next define,  
    for any faces $C'\leq C$ and $D$, the map of the stalks
 \[
 \delta_{CC'|D} = \delta_{CC'|D}^Q: E_{C\circ D}\lra E_{C'\circ D}
 \]
 by the formulas \eqref{eq:delta-Q}, that is,
 \[
 \delta_{CC'|D}= \gamma_{C'K'}\, \delta_{KC'}:  E_{C\circ D} =  E_K\lra E_{K'} =  E_{C'\circ D}, 
\]
where $K=C\circ D$ and $K'=C'\circ D$.  Note that
\be\label{eq:deltaCCD-phi}
\delta_{CC'|D} = \phi_{C\circ D, C'\circ D}. 
\ee

    \begin{prop}\label{prop:delta-comm-gen}
   The $ \delta_{CC'|D}$ commute with the generalization maps and so
   define morphisms of sheaves
   \[
   \underline\delta_{CC'}: \Ec_C \lra \Ec_{C'}, \quad C'\leq C. 
   \]
   \end{prop}

\noindent {\sl Proof:} For any  faces $C'\leq C$ and $D'\leq D$ we need to prove the commutativity of the
   square 
   \be\label{eq:hard-square}
   \xymatrix{
   E_{C'\circ D'}
   \ar[d]_{\gamma_{C'\circ D', C'\circ D}}
    & \ar[l]_{\delta_{CC'|D'}} E_{C\circ D'}
    \ar[d]^{\gamma_{C\circ D', C\circ D}}
   \\
   E_{C'\circ D} & \ar[l]^{\delta_{CC'|D}} E_{C\circ D}. 
   }
   \ee
   In order to do this, we use the bivariant notation of Example \ref{ex:base-change}, and consider the diagram
   of inclusions depicted by arrows
   \[
   \xymatrix{
   C'\circ D'
   \ar[dd]_{v'}
    &&& C\circ D' 
    \ar[dd]^v
   \\
   & C'\ar[ul]_{w'}
   \ar[r]^s 
   \ar[dl]^{u'}
   & C \ar[ur]^w
   \ar[dr]_u
   &
   \\
   C'\circ D &&& C\circ D.
    }
   \]
   Then
   \[
   \begin{gathered}
   \gamma_{C'\circ D', C'\circ D}\, \delta_{CC'|D'} \,\,=\,\, v'_* (w'_* s^* w^*) \,\, = \,\,
   u'_* s^* w^* \,\,=\,\, \gamma_{C', C'\circ D} \, \delta_{C\circ D', C'} \,\,=\,\, \phi_{C\circ D', C'\circ D},
   \\
   \delta_{CC'|D}\, \gamma_{C\circ D', C\circ D} \,\,=\,\, (u'_* s^* u^*) v_* \,\,=\,\, \phi_{C\circ D, C'\circ D} \, \phi_{C\circ D', C\circ D}.
   \end{gathered}
   \]
   In the last identification we used that $v_* = \gamma_{C\circ D', C\circ D}$ is the same as $\phi_{C\circ D', C\circ D}$. 
   So the commutativity of \eqref{eq:hard-square} would follow from (Tran) if we estabish the next lemma.
   
   \begin{lem}\label{lem:circ-collinear}
   For any faces $C'\leq C$, and $D'\leq D$,  the cells $C\circ D', C\circ D, C'\circ D$ form a collinear triple. 
   
   \end{lem}
   
   \noindent {\sl Proof of the lemma:} Choose points $c'\in C', c\in C, d'\in D', d\in D$. Assume that they are in  general position,
   i.e., $T=\on{Conv}\{c',c,d',d\}$ is a tetrahedron in an affine 3-space (the case $\dim(T)\leq 2$ is analyzed easily). 
   We need to find points $x\in C'\circ D$, $y\in C\circ D$, $z\in C\circ D'$ such that $y\in [x,z]$.

\begin{figure}
      \centering
    \begin{tikzpicture}[scale=0.3]
    
    \draw[thick] (0,0) -- (9,9) -- (-4,15) -- (-9,5) -- (0,0); 
    \draw[thick] (-4,15) -- (0,0); 
    \draw[dotted, line width =.4mm] (-9,5) -- (9,9); 
    
    \node (d2) at (-8,7){};
    \fill (d2) circle (0.15); 
    
     \node (d1) at (-7,9){};
    \fill (d1) circle (0.15); 
    
     \node (c2) at (3,3){};
    \fill (c2) circle (0.15); 
    
     \node (c3) at (5,5){};
    \fill (c3) circle (0.15); 
    
    \draw[thick] (0,0) -- (-7,9); 
    
    \draw[dotted, line width =0.3mm] (-8,7) -- (3,3); 
    
    \draw(-9,5) -- (-4.5,5); 
    \draw (-3.6,5) -- (-3.1,5); 
    \draw (-2.1,5) -- (-1.7,5); 
    \draw (-1.1,5) -- (5,5); 
    
    \node at (0,-1){$c'$}; 
    \node at (-10, 4.5){$d'$}; 
    \node at (10,9){$c$}; 
    \node at (-3.5, 15.5){$d$}; 
    \node at (4,2.5){$c_2$}; 
    \node at (6,4.5){$c_3$}; 
    \node at(-9,7){$d_2$}; 
    \node at (-8,9.5) {$d_1$}; 
      
    \end{tikzpicture}
     \caption{Collinearity of $C\circ D'$, $C\circ D$ and $C'\circ D$. }\label{fig:collinearity}

    \end{figure}
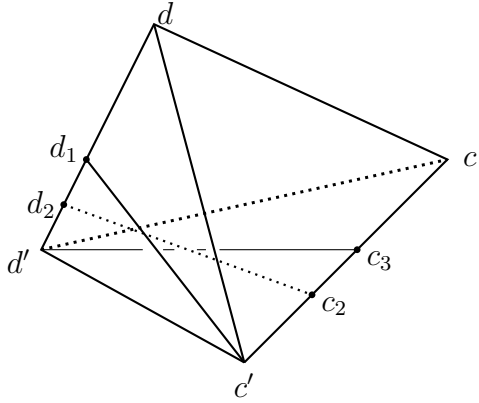

   To choose a possible $x$, we can draw the interval from $c'$ to any point $d_1\in (d',d]$ (any such
   $d_1$ satisfies $d_1\in D$)
   and take any point on this interval sufficietly close to $c'$. The supply of $x$ thus obtained contains a neighborhood of the vertex $c'$
   in the triangle $\triangle (c'd'd)$, with the edge $[c',d']$ removed. See Fig. \ref{fig:collinearity}. 
   
   To choose a possible $y$, we can draw the interval from any $c_2\in (c',c]$ to any $d_2\in (d',d]$ (as any such
   $c_2, d_2$ satisfy  $c_2\in C, d_2\in D$)
   and take any point on this interval sufficiently close to $c_2$. The supply of $y$ thus obtained contains a neighborhood of
   the edge $(c',c]$ in $T$, with the face $cd'd$ removed.

   To choose a possible $z$, we can similarly draw the interval from any $c_3\in (c',c]$ to $d'$ and take any point on this interval
   sufficiently close to $c_3$. The supply of $z$ thus obtained covers a neighborhood of the edge $(c'c]$ in the triangle
   $\triangle(d'c'c)$. 
   
   From this description it is clear that one can start from any admissible  $x\in [c',d_1]$ and take $y\in[c_2, d_2]$ sufficiently
   close to $c_2$, very near the face $\triangle(d'c'c)$. Then the interval $[x,y]$, continued after $y$, will hit $\triangle(d'c'c)$ in a point $z$ very close
   to the edge $[c',c]\ni c_2$, so such $z$ will be obtained by the above construction, i.e., will lie in $C\circ D'$. 
   This proves Lemma \ref{lem:circ-collinear} and Proposition \ref{prop:delta-comm-gen}.

 \begin{prop}\label{prop:mon+tran=comp}
    The morphisms of sheaves $\underline\delta_{CC'}$,
   $C'\leq C$, commute with each other, i.e., give rise to a contravariant representation of $(\Cc, \leq)$
   in $\Sh_{\CC^n}$.
   \end{prop}

 \noindent {\sl Proof:} 
      We have to prove the identity
  \be\label{eq:commut-delta-from-Q}
  \phi_{K_1' ,K''}\, \phi_{K,K_1'}\,\ =\,\,  \phi_{K'_2, K''}\, \phi_{K,K_2'}:
  E_{K}\lra E_{K''}
  \ee
  for any four faces $K, K'_1, K'_2, K''$ with the following property:
  there exists a square of codimension 1 inclusions
  \[
  C >_1  C'_1, C'_2 >_1  C''
  \]  
  and a face $D\geq C''$ such that
  \be\label{eq:4-compos}
  \begin{gathered}
  K=C\circ D, \quad K''= C''\circ D,\\
  K'_1= C'_1\circ D, \quad K'_2= C'_2\circ D.
  \end{gathered}
  \ee
  This would correspond to the commutativity of the stalks of the square
  \[
  \xymatrix{
  \Ec_C\ar[d]_{\underline\delta_{CC'_2}}
   \ar[r]^{\underline\delta_{CC'_1}}& \Ec_{C'_1}
   \ar[d]^{\underline\delta_{C'_1 C''}}
  \\
 \Ec_{C'_2} \ar[r]_{\underline\delta_{C'_2 C''}}& \Ec_{C''} 
}  
  \]
 over the face $[C'', D]$. Note that the $[C'', D]$ for $D\geq C''$, form a cell decomposition
 of the tube cell $\RR^n+i C''$, and the sheaves in question are direct image extensions
 of some sheaves from tube cells to their closures, so checking the commutativity of the above square over
 \[
 \RR^n + iC'' \,\, \subset \,\,\on{supp} (\Ec_{C''}) \,\,=\,\, \RR^n + i\overline C''
 \]
 is enough. 

To prove \eqref{eq:commut-delta-from-Q}, we  first remark that for any faces $A,B$ the triple $A, A\circ B, B$ is
collinear, since $A\circ B$ is defined in terms of points on the interval  $[a,b]$ for $a\in A, b\in B$. 
We also note that in \eqref{eq:4-compos} we have $K''=D$.  Consider the diagram of
inclusions
\[
\xymatrix{
K''=D & K'_1 & K
\\
C''\ar[u]\ar[r]& C'_1\ar[u] \ar[r]& C.\ar[u]
} 
\]
In this diagram, $C'_1, K'_1, K''$ form a collinear triple, since $K'_1=C'_1\circ D = C'_1\circ K''$. Therefore by (Tran),
\[
\phi_{C'_1, K''} \,\,=\,\, \phi_{K'_1, K''} \, \phi_{C'_1, K'_1} \,\,=\,\,  \phi_{K'_1, K''} \, \gamma_{C'_1, K'_1}. 
\]
Now, the LHS of the putatve equality \eqref{eq:commut-delta-from-Q} is transformed as follows:
\[
\begin{gathered}
\phi_{K'_1, K''} \, \phi_{K, K'_1} \,\,=\,\,\phi_{K'_1, K''} \, \gamma_{C'_1, K'_1} \, \delta_{K, C'_1}\,\,=
\\
=\,\, \phi_{C'_1, K''} \, \delta_{K, C'_1} \,\,=\,\, \gamma_{C'', K''}\, \delta_{C'_1, C''} \, \delta_{K, C'_1} \,\,=
\\
= \,\, \gamma_{C'', K''} \, \delta_{K, C''} \,\, = \,\,  \phi_{K, K''}. 
\end{gathered}
\]
Considering a similar diagram but with $C'_2, K'_2$ instead of $C'_1, K'_1$, we find that the RHS of 
 \eqref{eq:commut-delta-from-Q} is also equal to $\phi_{K, K''}$. \qed

\vskip .3cm

\noindent {\bf D. The complex $\Ec^\bullet(Q)$ is $\Sc^{(0)}$-smooth: inclusions of type (1).}
Proposition \ref{prop:mon+tran=comp} implies that 
we have a complex of sheaves $\Ec^\bullet = \Ec^\bullet(Q)$
with
\[
\Ec^p(Q) \,\,=\,\,\bigoplus_{\on{codim}(C)=p} \Ec_C\otimes\OR(C)
\]
and the differential $\widetilde{\underline \delta}$ given by the $\underline\delta_{CC'}$
 and satisfying  $\widetilde{\underline \delta}^2=0$. By construction, each summand of $\Ec^\bullet$
 and thus $\Ec^\bullet$ itself, is $\Sc^{(1)}$-smooth. 
 
 \begin{prop}\label{prop:E-C-smooth}
 The complex $\Ec^\bullet$ is $\Sc^{(0)}$-smooth. 
 \end{prop}
 
 \noindent {\sl Proof:} By Proposition \ref{prop:S1-to-S0}, it is enough to show that for any elementary inclusion
 $[C_1, D_1]\leq [C_2, D_2]$ of $\Sc^{(1)}$-cells, the generalization map 
 \[
 \gamma_{[C_1, D_1], [C_2, D_2]}: \Ec^\bullet|_{[C_1, D_1]} \lra \Ec^\bullet_{[C_2, D_2]}
 \]
 is a quasi-isomorphism of complexes of vector spaces. We first consider an inclusion of type (1)
 \[
 [C_1, D]\leq [C_2, D], \quad C_1 < C_2\leq D. 
 \]
 Note that every such inclusion is a composition of inclusions with $C_1 <_1 C_2$
 (if this condition does not hold, it suffices to choose a maximal chain of faces  $C_1 = C'_0 <_1 C'_1 <_1 \cdots <_1 C'_p =C_2$
 between $C_1$ and $C_2$ and consider the inclusions $[C'_i, D] < [C'_{i+1}, D]$). 
 So we assume $C_1 <_1 C_2$.

 \begin{lem}
 $\gamma_{[C_1, D], [C_2, D]}$ is a surjective morphism of complexes of vector spaces, with kernel 
  \[
  K^\bullet = \biggl\{ \bigoplus_{ \substack {\on{codim}(C)=0 \\ C\geq C_1, C\not\geq C_2}} E_{C\circ D}\otimes\OR(C)
  \buildrel \widetilde{\underline\delta} \over\lra 
  \bigoplus_{ \substack {\on{codim}(C)=1 \\ C\geq C_1, C\not\geq C_2}} E_{C\circ D}\otimes\OR(C)
  \buildrel \widetilde{\underline\delta} \over\lra \cdots\biggr\}, 
  \]
  where $\widetilde{\underline\delta}$ has, as matrix elements, the $\delta_{CC'}$ for relevant $C, C'$
  tensored with the coorientations of the adjacent faces. 
 
 \end{lem}
 
 \noindent {\sl Proof:} Recall that we defined the sheaf $\Ec_C$ by formulas 
\eqref{eq:stalks-Q}. Therefore  the stalk $\Ec_C|_{[C_1, D]}$ is equal to $E_{C\circ D}$,
 if $C\geq C_1$ and to $0$, if $C\not\geq C_1$, and similarly for $\Ec_C|_{[C_2, D]}$.
 The generalization maps were also defined by formulas of
  Corollary \ref{cor:stalks-E-C}(c). This implies that the matrix elements of 
  $\gamma_{[C_1, D], [C_2, D]}$ are either $\on{Id}$ or $0$, whence our statement. \qed
  
  \vskip .2cm
  
  We need therefore to prove that $K^\bullet$ is an exact complex. 
  
  \begin{lem}
  Consider the increasing filtration $F$  of $K^\bullet$ by graded subspaces defined by
  \[
  F_d K^p \,\,=\,\,\bigoplus_{\substack{
  \on{codim}(C)=p\\ C\geq C_1, C\not\geq C_2 \\
  \dim(C\circ D)\leq d
  }} E_{C\circ D}\otimes\OR(C). 
  \]
  This filtration is compatible with the differential $\widetilde\delta$. 
  \end{lem}

\noindent {\sl Proof:} The matrix element $\widetilde \delta_{CC'}$ is nonzero only if $C'\leq C$. In this case
$\LL(C'\circ D)\subset \LL(C\circ D)$ by Proposition \ref{prop:circ-monotone}(b), and so
 $\dim(C'\circ D)\leq\dim(C\circ D)$.
\qed

\vskip .2cm

We are therefore reduced to proving that each $\on{gr}^F_d K^\bullet$ is exact. 

\begin{lem}
Fix $d\geq 0$ and consider all $d$-dimensional flats $M\in\Lc$ containing $C_1$. Let
\[
(\on{gr}^F_d K)^p_M \,\,=\,\, 
\bigoplus_{\substack{
\on{codim}(C)=p\\
C\geq C_1, C\not\geq C_2\\
\LL(C\circ D)=M
}}
E_{C\circ D}\otimes\OR(C). 
\]
Then $(\on{gr}^F_d K)^\bullet_M$ is a subcomplex in $\on{gr}_d^FK^\bullet$, and we have a decomposition
into  a direct sum of complexes
\[
\on{gr}_d^FK^\bullet \,\,=\,\, \bigoplus_{M\supset C_1} (\on{gr}^F_d K)^\bullet_M. 
\]

\end{lem}

\noindent {\sl Proof:} Look at a summand $E_{C\circ D}\otimes\OR(C)$ with
\[
\on{codim}(C)=p, \,\, C\geq C_1, \,\, C\not\geq C_2,\,\, \dim(C\circ D)=d. 
\]
The differential in $\on{gr}^F_dK^\bullet$ can take this summand only to summands of the form
$E_{C'\circ D}$, where
\[
C'<_1 C, \,\,\on{codim}(C)=p+1, \,\, C'\geq C_1, \,\, C'\not\geq C_2,\,\, \dim(C'\circ D)=d. 
\]
Since $\LL(C'\circ D)\subset \LL(C\circ D)$ by Proposition  \ref{prop:circ-monotone}(b) and
$\dim \LL(C'\circ D) = \dim \LL(C\circ D)$, we conclude that $\LL(C'\circ D)=\LL(C\circ D)$.
So the complex $\on{gr}^F_DK^\bullet$ splits into a direct sum of
sub complexes corresponding to all possible  values of $M=\LL(C\circ D)$. \qed

\vskip .2cm

  Denote  $G^\bullet_M = (\on{gr}^F_d K)^\bullet_M$.
We are reduced therefore to proving that $G^\bullet_M$ is exact for each $M$. So we fix $M$ and note that by the
above
\[
G^p_M \,\,=\,\,\bigoplus_{\substack{
\on{codim}(C)=p \\
M\supset C\geq C_1, \,\, C\not\geq C_2\\
C\circ D \text{ open in } M
}}
E_{C\circ D}\otimes\OR(C). 
\]
Note that we can represent $G^\bullet_M$ as the quotient $G^\bullet_{1,M}/G^\bullet_{2,M}$, where
$G^p_{1,M}$ is the direct sum as above but without the restriction $C\not\geq C_2$, and $G^p_{2,M}$
is a similar direct sum but with the additional restriction $C\geq C_2$. So we are reduced to proving
that the embedding $G^\bullet_{2,M}\hookrightarrow G^\bullet_{1,M}$ is a quasi-isomorphism. 

In order that the summand corresponding to a face $C$ to be present in $G^p_{1,M}$, it is  necessary that
not only $C\subset M$, but also that $D\subset M$, since $M=\LL(C\circ D)$ is the minimal flat of
$\Hc$ containing $\LL(C)$ and $\LL(D)$. This means that $G^\bullet_{1,M}$ (as well as $G^\bullet_{2,M}$)
is entirely described in terms of faces of $\Hc$ contained in $M$, i.e., in  terms of the restricted
double quver $\Qc^{\leq M}$. So by passing to the restricted configuration $\Hc\cap M$, if necessary,
we can and will assume that $M=\RR^n$, and write $G^\bullet_\nu=G^\bullet_{\nu,\RR^n}$, $\nu=1,2$.
Further, since $C$ runs over faces containing $C_1$, the $G^\bullet_\nu$ are entirely described
in terms of the restricted double quiver $\Qc^{\geq C_1}$. So by passing to the quotient configuration
if necessary, we can and will assume that $C_1=\{0\}$, and therefore $\dim(C_2)=1$.

Under all these assumptions, let us give a geometric interpretation of the $G^\bullet_\nu$. Let $\Cc_D\subset
\Cc$ be the set of faces $C$ such that $C\circ D$ is open in $\RR^n$, and $|\Cc_D|$ be the union of cells
from $\Cc_D$.

\begin{lem}\label{lem:open-equiv}
Let $C,D\in\Cc$. Then the following are equivalent:
\begin{enumerate}
\item[(i)] $C\circ D$ is open in $\RR^n$.

\item[(ii)] $C$ does not lie in any hyperplane $H\in\Hc$ containing $D$. 
\end{enumerate}
\end{lem}

\begin{cor}\label{cor:CD-contract}
\[
|\Cc_D| \,\,\,=\,\,\, \RR^n \setminus\,\,\bigcup_{\substack{
H\in\Hc \\H\supset D  
}} \, H
\]
is the union of  convex, hence contractible components labelled by the (open) chambers of the quotient
arrangement $\Hc/\LL(D)$. \qed
\end{cor}

\noindent {\sl Proof of Lemma \ref{lem:open-equiv}:} We recall the interpretation of $C\circ D$ in terms of the projection
$\pi_C: \RR^n\to \RR^n/\LL(C)$ given in Proposition \ref{prop:circ-pi}. That is, $C\circ D$ is the unique face from $\Cc^{\geq C}$ 
which projects onto $\sigma(\pi_C(D))$, the cell of $\Hc/\LL(C)$ containing $\pi_C(D)$. So for $C\circ D$
to be open in $\RR^n$, it is necessary and sufficient that $\sigma(\pi_C(D))$ be open in
$\RR^n/\LL(C)$. This means that $\pi_C(D)$ does not lie in any hyperplane of
$\Hc/\LL(C)$, i.e., $D$ does not lie in any hyperplane of $\Hc$ containing $C$.
In other words, the sets of hyperplanes of $\Hc$ containing $C$ and $D$, must be disjoint. 
\qed

 \vskip .2cm

 Let us look at an arbitrary matrix element of  the differential in $G^\bullet_1$,
\[
 \widetilde {\ul \delta}_{CC'}: E_{C\circ D} \otimes\orr(C) \lra  E_{C'\circ D} \otimes\orr(C').
\]
It can be nonzero only if $C_1\leq C' <_1 C \subset M$, and both $C\circ D$ and $C'\circ D$ are open in $M$. If these conditions are 
satisfied, then \eqref{eq:deltaCCD-phi}
\[
 \widetilde {\ul \delta}_{CC'} \,\,=\,\,\phi_{C\circ D, C'\circ D}\otimes \epsilon_{CC'},
\]
where $\epsilon_{CC'}$ is the identification of orientation torsors induced by $C'<_1 C$. 
We note that in this case $\phi_{C\circ D, C'\circ D}$ is an isomorphism. Indeed,
since both $C\circ D$ and $C'\circ D$ are open in $M$, we can choose generic points
$x\in C\circ D$, $x'\in C'\circ D$ so that the interval $[x,x']\subset M$ intersects only flats of $\Hc\cap M$
which have codimension 1 in $M$.  Therefore we can write all the faces intersecting $[x,x']$, in the order
from $x$ to $x'$, as
\[
F_0=C\circ D >_1 F'_0 <_1 F_1 >_1 F'_1 <_1 \cdots >_1 F'_{p-1} <_1 F_p = C\circ D', 
\]
where the $F_j$ are open in $M$, while the $F'_j$ have codimension 1 in $M$. Therefore each
$\phi_{F_j, F_{j+1}}$ is an isomorphism by (Inv), while
\[
\phi_{C\circ D, C'\circ D} \,\,=\,\,\phi_{F_0, F_p} \,\,=\,\, \phi_{F_{p-1}, F_p} \, \cdots \, \phi_{F_0, F_1}
\]
by (Tran) (long form of the transitivity relations).

We now note that the isomorphisms $\phi_{C\circ D, C'\circ D}$
 define a cellular locally constant (and therefore, by Corollary \ref{cor:CD-contract},
constant)
sheaf $\Gc$ on $|\Cc_D|$, and the above argument shows that
\[
\Gc^\bullet_1 \,\,=\,\, C^{\on{cell}}_{n-\bullet}(|\Cc_D|, \Gc)
\]
is the standard cellular chain complex of $\Gc$ shifted to so as to start in degree 0
(and to have differential of degree $-1$). Similarly, let $\Cc_D^{\geq C_2}$ be the subposet in
$\Cc_D$ formed by $C$ satisfying $C\geq C_2$ and $|\Cc_D^{\geq C_2}|$ be the union of
faces from this subposet. Then
\[
\Gc^\bullet_2 \,\,=\,\,C^{\on{cell}}_{n-\bullet}(|\Cc_D^{\geq C_2}|, \Gc). 
\]
Therefore, the acyclicity of $G^\bullet_1/G^\bullet_2$ will follow if we prove that the
embedding of spaces $|\Cc_D^{\geq C_2}|\hookrightarrow |\Cc_D|$ is a homotopy equivalence.

  Recall that $C_2\leq D$, $\dim(C_2)=1$ is, by our assumption,
a half-line. Therefore for any connected component $U\subset |\Cc_D|$, the intersection
$U\cap |\Cc_D^{\geq C_2}|$ is a non-empty convex set so it is contractible as well. This proves
that $|\Cc^{\geq C_2}_D|\hookrightarrow |\Cc_D|$ is a homotopy equivalence and so
$\gamma_{[C_1, D], [C_2, D]}$ is a quasi-isomorphism. 

\vskip .3cm

\noindent {\bf E. Inclusions of type (2).} To finish the proof of Proposition \ref{prop:E-C-smooth}, we need
to consider elementary inclusions of type (2). We write such an inclusion as $[C_1, D_1]\leq[D_2, D_2]$
where $D_1, D_2$ are two faces of the same dimension $m$, lying in the same $m$-dimensional flat
$L=\LL(D_1)=\LL(D_2)$ on the opposite side of an $(m-1)$-dimensional cell $C_1\leq D_1, D_2$. 

The argument runs very similarly to the case of an inclusion of type (1) considered in \S D. We indicate the
changes, using the same notation as in \S D for the intermediate complexes and  and treating
the grading in these complexes as implicit.

\begin{lem}
The morphism of complexes $\gamma_{[C_1, D_1], [D_2, D_2]}$ is surjective with kernel
\[
K^\bullet \,\,=\,\,\bigoplus_{\substack{
C\geq C_1\\ C\not\geq D_2
}}
E_{C\circ D_1}\otimes\OR(C),
\]
graded by $\on{codim}(C)$. 
\end{lem}

\noindent {\sl Proof:} By Corollary \ref{cor:stalks-E-C}(c), 
\[
\Ec^\bullet_{[C_1, D_1]} \,\,= \bigoplus_{C\geq C_1} E_{C\circ D_1}\otimes\orr(C), \quad
\Ec^\bullet_{[D_2, D_2]} \,\,= \bigoplus_{C\geq D_2} E_{C\circ D_2}\otimes\orr(C). 
\]
The set of admissible $C$ in the sum for $\Ec^\bullet_{[D_2, D_2]}$ is clearly a subset in the set of admissible
 $C$ in the sum for $\Ec^\bullet_{[C_1, D_1]}$, because $C_1< D_2$. Further, for $C\geq D_2$
the corresponding summand in $\Ec^\bullet_{[C_1, D_1]}$ is equal to the corresponding summand in
$\Ec^\bullet_{[D_2, D_2]}$. Indeed, since $D_1$ and $D_2$ are adjacent (have the same dimension $m$ and lie
on the opposite sides of $C_1$ in an $m$-flat), we have $D_2\circ D_1 = D_2$. We also have $C\circ D_2=C$
since $C\geq D_2$, so by the associativity of $\circ$ we find
\[
C\circ D_1 \,\,=\,\, (C\circ D_2)\circ D_1 \,\, = \,\, C\circ(D_2\circ D_1) \,\,=\,\, C\circ D_2. 
\]
Further, for $C\geq D_2$,  
Corollary \ref{cor:stalks-E-C}(c)  shows that
the matrix element of $\gamma^{\Ec^\bullet}_{[C_1, D_1], [D_2, D_2]}$ on the summand corresponding to
$C$ is the identity. Recalling that the matrix elements of  $\gamma^{\Ec^\bullet}_{[C_1, D_1], [D_2, D_2]}$
can act only between summands labelled by the same $C$, we conclude that it
is the projection onto  the direct sum of summands labelled by $C\geq D_2$, as claimed. \qed

\vskip .2cm

Now, as before, we have the increasing filtration $F$ in $K^\bullet$ by $\dim(C\circ D_1)$
with quotients that split into direct sum over $M\in \Lc$ of complexes
\[
G^\bullet_M \,\,=\,\,\bigoplus_{\substack{
C\geq C_1, \, C\not\geq D_2\\
\LL(C\circ D_1)=M
}}
E_{C\circ D_1}\otimes\OR(C).
\]
So we need to prove that each $G^\bullet_M$ is exact. As in \S D, we have $G^\bullet_M=G^\bullet_1/G^\bullet_2$
where $G^\bullet_1$ is the direct sum over $C\geq C_1$ such that $L(C\circ D_1)=M$ (so it is exactly the
same complex as in \S D), and $G^\bullet_2$ is the subcomplex formed by $E_{C\circ D_1}\otimes\OR(C)$
for $C\geq D_2$. As before, we can and will assume that $C_1=0$ and $M=\RR^n$. Thus $D_1, D_2$
are  two opposite half-lines of the same line: $D_2=-D_1$. So $G^\bullet_1$ is the cellular chain complex of
a local system $\Gc$ on $|\Cc_{D_1}|$ and $G^\bullet_2$ is the chain complex of $\Gc$ on
$|\Cc_{D_1}^{\geq D_2}|$. Since $D_2=-D_1$, the subspace $|\Cc_{D_1}^{\geq D_2}|$ is again
a union of convex connected components, one inside each convex connected component of 
$|\Cc_{D_1}|$, so the embedding $G^\bullet_2\hookrightarrow G^\bullet_1$ is a quasi-isomorphism,
and Proposition \ref{prop:E-C-smooth} is proved. 

\vskip .3cm

\noindent {\bf F. End of the proof of Reformulation \ref{refo:mainref}.  }
Given $Q\in\Rep^{(2)}(\Cc)$ satisfying (Mon), (Tran) and (Inv), we have associated
to it an $\Sc^{(0)}$-constructible complex $\Ec^\bullet(Q)$ on $\CC^n$. Note that the
orthogonality relations of Lemma  \ref{lem:orthogonality} apply to the $\Ec_C(Q)$ and imply 
that the linear data of $\Ec^\bullet(Q)$ are given by $Q$.  Since $Q$ consists of
single vector spaces (not just complexes), Proposition
\ref {prop:perversity-char} implies that $\Ec^\bullet(Q)$ is a perverse sheaf,
with double quiver $Q$. This proves Reformulation \ref{refo:mainref} and thus
Theorem \ref{thm:main}.

\section{Examples and complements}\label{sec:excom}

\noindent {\bf A. The case of dimension 1.} Suppose $n=1$. The real vector space $\RR$ has the unique
hyperplane $\{0\}$. The arrangement consisting of this hyperplane will be denoted simply by 0.
It has three faces
$\RR_+ = \{ x>0\}$, $\RR_-=\{x<0\}$ and $\{0\}$. The category $\Perv(\CC,0)$ consists of perverse
sheaves on $\CC$ smooth with respect to the stratification consisting of $\{0\}$ and $\CC\setminus \{0\}$. 
The classical description \cite{galligo-GM} identifies $\Perv(\CC,0)$  with the category $\Pc$  of diagrams of
finite-dimensional vector spaces
\[
  \xymatrix{
 \Phi \ar@<.7ex>[r]^v& \Psi \ar@<.7ex>[l]^u
  }, \quad vu+\on{Id}_\Psi \text{invertible}.
\]
The spaces $\Psi$ and $\Phi$ associated to $\Fc\in \Perv(\CC,0)$ are canonically identified with
the spaces of nearby and vanishing cycles of $\Fc$ with respect to the standard coordinate function on $\CC$,
while $vu+\on{Id}_\Psi$ is the monodromy on the space of nearby cycles, 
see \cite{beil-gluing}. 

On the other hand, our description from Theorem \ref{thm:main} identifies $\Perv(\CC,0)$ with the
category $\Ac$ formed by double representations
\[
\begin{gathered}
\Qc = \bigl\{
\xymatrix{
E_-\ar@<-.7ex>[r]_{\delta_-}& E_0
\ar@<-.7ex>[l]_{\gamma_-}
\ar@<.7ex>[r]^{\gamma_+}& E_+
\ar@<.7ex>[l]^{\delta_+}
}
\bigr\} \quad \text{such that}\\
\gamma_-\delta_- = \on{Id}_{E_-}, \,\,\,\gamma_+
\delta_+ = \on{Id}_{E_+},\\
\gamma_-\delta_+: E_+\lra E_-, \,\,\, \gamma_+\delta_-: E_-\to E_+ \quad \text {are invertible}. 
\end{gathered}
\]
Let us construct an equivalence between $\Pc$ and $\Ac$ directly. In fact, it is convenient to reformulate the
definition of $\Ac$ slightly. Given $\Qc\in\Ac$, consider endomorphisms
\[
P_+ = \delta_+\gamma_+, \,\,\, P_- = \delta_-\gamma_- \,\,\,\in\,\,\, \on{End}(E_0). 
\]
These endomorphisms are idempotent:
\[
P_+^2 \,\,=\,\, \delta_+\gamma_+\delta_+\gamma_+ \,\,=\,\, \delta_+ \on{Id} \gamma_+ \,\,=\,\, \delta_+\gamma_+ \,\,=\,\, P_+, 
\]
and similarly for $P_-$. The spaces $E_\pm$ are identified with $\on{Im}(P_\pm)$ via $\delta_\pm$. The
conditions of invertibility of $\gamma_\pm\delta_\mp$ is expressed by:
\be\label{eq:P-isom}
\begin{cases} P_-: \on{Im}(P_+)\lra \on{Im}(P_-), \\ P_+: \on{Im}(P_-)\lra \on{Im}(P_+)
\end{cases}
\quad \text{are isomorphisms.}
\ee
This establishes the following.

\begin{lem}
The category $\Ac$ is equivalent to the category $\Bc$ formed by data $(E_0, P_+, P_-)$
consisting of a finite-dimensional $\k$-vector space $E_0$ and two idempotents $P_+, P_-: E_0\to E_0$
satisfying \eqref{eq:P-isom}. \qed
\end{lem}

So we will construct an equivalence $\Bc\simeq \Pc$. Given $(E_0, P_+, P_-)\in\Bc$, we put
\[
\xymatrix{
\Phi:= \on{Ker}(P_-) \ar@<.7ex>[rr]^{v=P_+} && \Psi :=\on{Im}(P_+)
\ar@<.7ex>[ll]^{u=P_--\on{Id}}.
}
\]
Then $vu=P_+(P_-  - \on{Id})$ is, as an endomorphism of $\on{Im}(P_+)$, equal to $P_+P_- - P_+$. Now, on
$\on{Im}(P_+)$ we have $P_+ = \on{Id}$. So $vu+\on{Id}= P_+ P_-$ as an endomorphism of $\on{Im}(P_+)$
and so it is the composition of two invertible maps
\[
\on{Im}(P_+) \buildrel P_- \over\lra \on{Im}(P_-) \buildrel P_+ \over\lra \on{Im}(P_+)
\]
and hence invertible. This defines a functor
$F: \Bc\to\Pc$. Note that the two maps above have the meaning of half-monodromies from the
upper to the lower half plane, so $vu+\on{Id}$ is the full monodromy. 

Let us also define the functor $G:\Pc\to\Bc$ as follows. Given an object
$\{\xymatrix{
 \Phi \ar@<.7ex>[r]^v& \Psi \ar@<.7ex>[l]^u
  }\}\in\Pc$, we put
  \be\label{eq:P-pm-defined}
  E_0=\Phi\oplus\Psi, \quad P_+ = \begin{pmatrix} 0&0\\ v&1
  \end{pmatrix}, \quad P_- = \begin{pmatrix} 0&u\\ 0&1
  \end{pmatrix}.
  \ee
  Then $P_\pm$ are idempotents. 
  
  \begin{prop}
  The functors $F$ and $G$ are quasi-inverse to each other. 
  \end{prop}
  
  \noindent {\sl Proof:} Let us find $FG$. For $P_\pm$ defined above, we have
  \[
  \begin{gathered}
  \on{Im}(P_+) = \Psi, \quad \on{Ker}(P_+) = \biggl\{
  \begin{pmatrix} \phi\\ \psi
  \end{pmatrix}\biggl| \, v\phi + \psi = 0\biggr\}\,\,=\,\, \text{Graph of } (-v): \Phi\to\Psi, \\
  \on{Im}(P_-) = \text {Graph of } u: \Psi\to \Psi, \quad \on{Ker}(P_-)=\Phi. 
  \end{gathered}
  \]
  The map $P_+: \on{Ker}(P_-)\to\on{Im}(P_+)$ coincides with $v$. Further,
  $P_--\on{Id}$ restricted to $\Psi\subset \Phi\oplus\Psi$ gives $u: \Psi\to\Phi$. Therefore
  $FG$ is isomorphic to the identity functor of $\Pc$.
  
  Conversely, suppose $(E_0, P_+, P_-)\in\Bc$. We then have two direct sum decompositions of $E_0$:
  \[
  E_0 \,\,\,=\,\,\, \Phi\oplus \on{Im}(P_-) \,\,=\,\, \on{Ker}(P_+)\oplus\Psi. 
  \]
The condition that $P_-: \on{Im}(P_+) \to\on{Im}(P_-)$
is an isomorphism, implies that $\on{Im}(P_+) \cap \on{Ker}(P_-)=0$ which, by the dimension count implies that
we have a direct sum decomposition
\[
E_0 \,\,=\,\, \on{Im}(P_+) \oplus\on{Ker}(P_-) \,\,=\,\, \Psi\oplus\Phi.
\]
With respect to this decomposition, we find that $P_\pm$ are given by the matrices
in \eqref{eq:P-pm-defined}. So $GF$ is isomorphic to the identity functor of $\Bc$. \qed

\vskip .2cm

\begin{rem}
Note that the  composite equivalence $\Ac\to\Bc\buildrel F\over\to\Pc$ is  compatible with  (and so
can be considered as induced by) the identifications of
$\Ac$ and $\Pc$ with $\Perv(\CC,0)$, constructed in \cite{galligo-GM} and  in this paper, respectively. 
Indeed, let $\Fc\in\Perv(\CC,0)$. Then, in the original construction of  \cite{galligo-GM}, the object of $\Ac$
corresponding to $\Fc$ is given by:
\be\label{eq:phi-and-psi}
\Phi(\Fc) = \HH^1_{\RR_{\geq 0}} (\CC,\Fc) = \Gamma(\RR_{\geq 0}, \ul\HH^1_{\RR_{\geq 0}}(\Fc)), \,\,\,  \Psi(\Fc) 
=\Gamma( \RR_{>0}, \ul \HH^1_{\RR_{\geq 0}}(\Fc)),
\ee
and the map  $v$ is the generalization map for the $\RR$-constructible sheaf $\ul\HH^1_{\RR_{\geq 0}}(\Fc)$ on $\RR_{\geq 0}$. 
On the other hand, the object $Q(\Fc)\in \Pc$ corresponding to $\Fc$  by \eqref {eq:functor-Q} has
\be\label{eq:our-identific}
E_0(\Fc) = \HH^1_\RR(\CC, \Fc) = \Gamma(\RR, \ul\HH^1_\RR(\Fc)), \quad E_\pm(\Fc) = \Gamma(\RR_{\gtrless 0},
 \ul \HH^1_{\RR}(\Fc)),
\ee
and $\gamma_\pm$ is the generalization map for the $\RR$-constructible sheaf $\ul\HH^1_\RR(\Fc)$ on $\RR$. 
The functor $F$ sends $Q=Q(\Fc)$ into an object of $\Ac$ with
\[
\Phi=\Phi_Q :=\on{Ker}(P_-)=\on{Ker}(\gamma_-), \quad \Psi=\Psi_Q:= \on{Im}(P_+) = \on{Im}(\delta_+).
\]
Now, the identification
  $\on{Ker}(\gamma_-) \to \Phi(\Fc)$ is obtained from the  long exact sequence relating (hyper)cohomology with supports in $\RR$,
  $\RR_{\geq 0}$ and $\RR_{<0}$. The identification $\on{Im}(\delta_+)\to\Psi(\Fc)$
  is obtained by comparing \eqref{eq:phi-and-psi} and \eqref{eq:our-identific} and noting that $\ul\HH^1_{\RR_{\geq 0}}(\Fc)$
  and $\ul\HH^1_\RR(\Fc)$ are identified on $\RR_{> 0}$. 
We leave to the reader the identification of the arrows between the vector spaces thus identified. 

\end{rem}

\begin{rem} In addition to the  above equivalence which we denote here  
$F_+: \Bc\overset\sim\to\Pc$ there exists another one 
$F_-: \Bc\overset\sim\to\Pc$ where for $x = (E_0, P_+, P_-) \in \Bc$
we set 
$$
F_-(x) = 
\bigl\{\xymatrix{
\Phi:= \on{Ker}(P_+) \ar@<.7ex>[rr]^{v=P_-} && \Psi :=\on{Im}(P_-)
\ar@<.7ex>[ll]^{u=P_+ -\on{Id}}
}\bigr\}
$$
Furthermore, there are two invertible natural transformations  
("half-monodromies") $t_\pm: F_\pm \overset\sim\to F_\mp$ given by
$$
t_+(x) = (P_+ - \on{Id}, P_-):\ F_+(x) \overset\sim\to F_-(x)
$$
and
$$
t_-(x) = (P_- - \on{Id}, P_+):\ F_-(x) \overset\sim\to F_+(x)
$$
These data define a $\CC^*$-local system of equivalences 
$\Bc \to \Pc$; this is a 
\newline particular case of the {\it microlocalization}, similar to \cite{fs}.   
 
\end{rem}

\vskip .2cm

Finally, let us describe the effect of the geometric Fourier transform $\Fen:\Perv(\CC,0)\to\Perv(\CC,0)$
in terms of our description. In the classical description, as well known,
\[
\Fen \{\xymatrix{
 \Phi \ar@<.7ex>[r]^v& \Psi \ar@<.7ex>[l]^u
  }\}\,\,\simeq \,\,\{\xymatrix{
 \Psi \ar@<.7ex>[r]^u& \Phi \ar@<.7ex>[l]^v
  }\}
\]
We leave to the reader the verification of the following.

\begin{prop} In terms of the identification $\Perv(\CC, 0)\simeq\Bc$,
\[
\Fen (E_0, P_1, P_2) \,\,\simeq \,\, (E_0, 1-P_1, 1-P_2), 
\]
and in terms of  the identification $\Perv(\CC,0)\simeq \Ac$,
\[ \Fen\bigl\{
\xymatrix{
E_-\ar@<-.7ex>[r]_{\delta_-}& E_0
\ar@<-.7ex>[l]_{\gamma_-}
\ar@<.7ex>[r]^{\gamma_+}& E_+
\ar@<.7ex>[l]^{\delta_+}
}
\bigr\} \,\,\,\simeq \,\, \,
\bigl\{
\xymatrix{
\on{Ker}(\gamma_-) \ar@<-.7ex>[rr]_{\text{embedding}}&& E_0
\ar@<-.7ex>[ll]_{\on{Id}-\delta_-\gamma_-}
\ar@<.7ex>[rr]^{\on{Id}-\delta_+\gamma_+}&& \on{Ker}(\gamma_+)
\ar@<.7ex>[ll]^{\text{embedding}}
}
\bigr\}.  
\]
\qed

\end{prop}

\noindent {\bf B. Real affine arrangements.}
 Let $\Hc$ be an arrangement of {\em affine} hyperplanes in $\RR^n$ and $\Hc_\CC$
 be the complexified arrangement of affine hyperplanes in $\CC^n$. We have then the
 category $\on{Perv}(\CC^n, \Hc)$ of $\Hc$-smooth perverse sheaves on $\CC^n$,
 similarly to the case of linear arrangements. 
 
 Such categories are important for the geometric description 
 of tensor structures on the categories of quantum group 
 representations, ~\cite{BFS}.

 Theorem \ref{thm:main} is extended to this case as follows. We have the poset of faces
 $(\Cc, \leq)$, defined similarly to the linear case. A triple of faces $(A,B,C)$ is called
 {\em collinear}, if:
 \begin{itemize}
 \item[(C1)] There exists a face $D\leq A,B,C$.
 \item[(C2)] There exists points $a\in A, b\in B, c\in C$ such that $b\in [a,c]$. 
 \end{itemize}
The condition (C1) holds automatically in the linear case (take $D =0$). 
As with any poset,  we have the category $\on{Rep}^{(2)}(\Cc)$ of double representations
$Q= (E_C, \gamma_{C'C}, \delta_{CC'})$ of $\Cc$. 

The condition of {\em monotonicity} on $Q$ is defined just as in the linear case:
$\gamma_{C'C}\delta_{CC'}=\on{Id}$ for any $C'\leq C$. This allows us to define the transition
maps $\phi_{AB}=\gamma_{CB}\delta_{AC}:E_A\to E_B$ for any two faces $A,B$ such that
there is a face $C\leq A,B$. 

The condition of {\em transitivity} is defined by requiring that $\phi_{AC}=\phi_{BC}\phi_{AB}$ for
any triple of cases $(A,B,C)$ collinear in the new sense above. 

Finally, the condition of {\em invertibility} of $Q$ is defined completely similarly to the linear case. 

We denote by $\Ac=\Ac_\Hc$ the full subcategory in $\on{Rep}^{(2)}(\Cc)$ formed by
$Q$ which are monotone, transitive and invertible.

\begin{thm}\label{thm:affine}
The category $\on{Perv}(\CC^n, \Hc_\CC)$ is equivalent to $\Ac_\Hc$. 
\end{thm}

\noindent {\sl Proof:} Locally (in a neighborhood of any point of $\RR^n$), an affine arrangement
looks like a linear one. Therefore, applying known properties of the linear case, we establish the
following statements:

\begin{itemize}
\item[(1)] For $\Fc\in\on{Perv}(\CC^n, \Hc)$ the complex $\Rc_\Fc = \underline{R\Gamma}_{\RR^n}(\Fc)[n]$
is reduced to one sheaf in degree 0, constructible with respect to the (cellular) stratification $\Cc$

\item[(2)] Denoting by $E_C(\Fc)$ the stalk of $\Rc_\Fc$ at $C\in\Cc$, we have a canonical identification
$\Ec_C(\Fc^\bigstar) = E_C(\Fc)^*$. 

\item[(3)] Denoting  $\gamma_{C'C}^\Fc: E_{C'}(\Fc)\to E_C(\Fc)$, $C'\leq C$, the generalization maps
for $\Rc_F$, and putting $\delta_{CC'}^\Fc = (\gamma_{C'C}^{\Fc^\bigstar})^*$, we get an object
\[
\Qc(\Fc) \,\,=\,\,\bigl( E_C(\Fc), \,\gamma_{C'C}^\Fc, \, \delta_{CC'}^\Fc\bigr) \,\,\in \,\,\Rep^{(2)}(\Cc).
\]
This object lies in $\Ac$. 
\end{itemize}

We get therefore a functor $\Qc: \on{Perv}(\CC^n, \Hc_\CC)\to \Ac$ taking $\Fc$ to $\Qc(\Fc)$.
 To prove that it is an equivalence,
we use the fact that perverse sheaves form a stack. That is, we upgrade $\Qc$ to a morphism of
stacks of abelian categories on $\RR^n$
\[
\underline{\Qc}: \underline{\on{Perv}}_{(\CC^n, \Hc)} \lra \underline\Ac
\]
defined as follows. For an open $U\subset \RR^n$, the category $\underline{\on{Perv}}_{(\CC^n, \Hc)}(U)$
consists of perverse sheaves on $U+i\RR^n$ smooth with respect to the stratification cut out on $U+i\RR^n$
by $\Hc$. We further denote by $\Cc_U\subset \Cc$ the subset of faces meeting $U$ and extend the
concept of collinearity to $\Cc_U$ be requiring in (C1) that $D\in \Cc_U$. Then we define
$\underline\Ac(U)\subset\Rep^{(2)}(\Cc_U)$ to be the full subcategory specified by the
conditions of monotonicity, transitivity and invertibility. The same arguments as before define a functor
$\underline{\Qc}(U): \underline{\on{Perv}}_{(\CC^n, \Hc)}(U)\to \underline\Ac(U)$, and these functors
unite into a morphism of stacks $\underline{\Qc}$. Note that $\Qc$ is obtained from
 $\underline{ \Qc}$ by
passing to the categories of global sections. So it is enough to show that $\underline{ \Qc}$ is an equivalence
of stacks, a statement that can be checked locally, at the level of stalks. But the stalk functor of 
$\underline{\Qc}$
over any point $x\in\RR^n$ is a similar functor $\Qc$ for the (essentially) linear configuration formed by
hyperplanes from $\Hc$ containing $x$. So it is an equivalence by Theorem \ref{thm:main}. \qed

\vskip .2cm

The method of using the stacky nature of perverse sheaves to obtain descriptions in new situations
was applied by Dupont to the case of smooth toric varieties \cite{dupont}. 

\vskip .3cm

\noindent {\bf C. The fundamental groupoid of the open stratum.}
Let us write $V=\RR^n$, so  $V_\CC=\CC^n$, and let $V_\CC^\circ\subset V_\CC$ be the open stratum (the complement of all the hyperplanes
$L_\CC$, for $L\in \Hc$). Let $\Cc_0$ be the set of  open faces (chambers) of $\Hc$, so each $A\in\Cc_0$ is 
a contractible set contained in $V_\CC^\circ$. Denote $\pi_1(V_\CC^\circ, \Cc_0)$ the fundamental groupoid of $V_\CC^\circ$
with respect to a set of base points consisting of one point in each $A\in\Cc_0$. The construction of this paper
leads to a new description of this groupoid. 

\begin{prop}\label{prop:GK-groupoid}
 $\pi_1(V_\CC^\circ, \Cc_0)$ is isomorphic to the groupoid $\Gen$ defined by generators and relations as follows:
 \begin{itemize}
 \item[(0)] Objects $x_A$, $A\in \Cc_0$. 
 
 \item[(1)] Generating morphisms $\varphi_{AB}: x_A \to x_B$ for each ordered pair $(A,B)$ of chambers.
 We assume that $\varphi_{AA}=\on{Id}_{x_A}$. 
 
 \item[(2)] Relations $\varphi_{AC}=\varphi_{BC} \varphi_{AB}$ for any collinear triple $(A,B,C)$ of chambers. 
 \end{itemize}
\end{prop}

A more familiar description of  $\pi_1(V_\CC^\circ, \Cc_0)$  is the one  from the work of Salvetti \cite{salvetti}
which we now recall.

\begin{prop}\label{prop:salvetti}
 $\pi_1(V_\CC^\circ, \Cc_0)$ is isomorphic to the groupoid $\Sen$ defined by generators and relations as follows:
 \begin{itemize}
 \item[(0)] Objects $x_A$, $A\in \Cc_0$. 
 
 \item[(1)] Generating morphisms $\psi_{AB}: x_A \to x_B$ for each ordered pair $(A,B)$ of chambers
 which are adjacent, i.e., lie on opposite sides of a codimension 1 face $C$.
 
 \item[(2)]  The  Zifferblatt relations
 \[
 \psi_{B_{m}, C}\, \psi_{B_{m-1}, B_{m}} \, \cdots \, \psi_{A,  B_2} \,\, = \,\,
 \psi_{B_{m+2}, C} \, \psi_{B_{m+3}, B_{m+2}} \, \cdots \, \psi_{B_{2m}, B_{2m-1}} \, \psi_{A, B_{2m}}
 \]
 for any codimension 2 face $F$, and any chamber $A>F$. Here we number all the chambers $>F$ around the circle
 as $A=B_1, B_2, \cdots, B_{m+1}=C, B_{m+1}, \cdots, B_{2m}$, as in Fig. \ref{fig:Zifferblatt}.
 \end{itemize}
\end{prop}

Proposition \ref{prop:salvetti} can be obtained by noticing that $\Sc^{(1)}$-cells $[C,A]$ for $A$ being a chamber, form
a cell decomposition of $V_\CC^\circ$.  
Among them the cells $[A,A]$ are precisely the open ones. Therefore one can form the dual CW-complex.
denote it $\Sc al$, 
homotopy equivalent to $V_\CC^\circ$,  in which   each cell $[A,A]$ will give rise to a vertex,  each cells $[C,A]$ with
$\on{codim}(C)=1$ will give rise to an edge, and each cell $[F,A]$ with $\on{codim}(C)=2$ will give rise to
a $2m$-gon, with $2m = \# \{ B\in \Cc_0 |\,  B>F\}$ and so on. The groupoid $\Sen$ is the groupoid whose presentation  is obtained
in the standard way, from
the 2-skeleton of $\Sc al$, see \cite{BZ} \cite{salvetti} for more details. 

\vskip .2cm

We now prove Proposition \ref{prop:GK-groupoid}. Define a functor $F: \Sen\to\Gen$ to be identity on the objects
and to send $\psi_{AB}$ (with $A,B$ adjacent) to $\varphi_{AB}$. Note that the relations of $\Sen$ are 
satisfied in $\Gen$, see Example \ref{ex:zifferblatt}, so $F$ is well defined. Let us prove that $F$ is an isomorphism
of groupoids.  

For this, we define a functor $G: \Gen\to \Sen$, also identical on objects,  as follows. Let $A,B$ be two chambers.
Choose any two generic points $a\in A, b\in B$ so that the interval $[a,b]\subset \RR^n$ does not intersect
any faces of codimension $\geq 2$. Denote the chambers intersecting $[a,b]$, if written in the direction from
$a$ to $b$, by $C_1 = A, C_2, \cdots, C_r=B$. Then each $(C_i, C_{i+1})$ form an adjacent pair, so
the generator $\psi_{C_i, C_{i+1}}$ of $\Sen$ is defined, and we put
\be\label{eq:g-AB}
G(\varphi_{AB}) \,\,=\,\, \psi_{C_{r-1}, B} \,\psi_{C_{r-2}, C_{r-1}} \, \cdots \, \psi_{A, C_2}. 
\ee

\begin{lem}
The RHS of  \eqref{eq:g-AB}, considered as an element of $\Hom_{\Sen}(x_A, x_B)$,  is independent
on the choice of generic points $a\in A, b\in B$.
\end{lem}
 
 \noindent {\sl Proof:} It is enough to prove that if we keep $a$ and replace $b$ by another generic point $b'\in B$,
 or if we keep $b$ and replace $a$ by another generic point $a'\in A$, then the RHS of \eqref{eq:g-AB}
 will give the same morphism. Let us consider the first situation, the second one is treated similarly.
 
 Given $a,b, b'$, consider the plane triangle $\Delta  = \on{Conv}\{a,b,b'\}$ 
 inside the affine 2-plane $P$ spanned by $a,b,b'$. The side $[b,b']$ of $\Delta$ lies inside $B$.  
 Note that $\Hc$ induces an affine arrangement of lines inside $P$, and by our assumption, flats of $\Hc$
 of codimension $\geq 3$ do not meet $P$.
 Each codimension 2 face  $L$ of $\Hc$ intersecting $\Delta$,
 does so at an interior point of  $\Delta$. Now, around each such  point $L\cap \Delta$ we have a 
 Zifferblatt situation. That is, we can deform the path $\xi_0 = [a,b']\cup [b',b]$ into $\xi_1 = [a,b]$ 
 in a family of paths $(\xi_t)_{t\in[0,1]}$,  keeping the endpoints
 $a,b$ fixed so that at every moment $t$ we cross at most one of the points $L\cap\Delta$. Associating to each intermediate
 path $\xi_t$ the product of generators $\psi$ similar to \eqref{eq:g-AB}, we see that after crossing each $L\cap\Delta$,
 the product remains unchanged in virtue of the relations of $\Sen$. \qed
 
 \vskip .2cm
 
 With the lemma establshed, we see that $G$ preserves the relations of $\Gen$ by its very definition:
 for collinear chambers $A,B,C$ with points $a,b,c$ such that $c\in [a,b]$, we use the intervals
 $[a,b], [b,c], [a, c]$ to define the values of $G$ on $\varphi_{AB}, \varphi_{BC}, \varphi_{AC}$.
 So $G$ is indeed a functor and we see that it is inverse to $F$ by looking at the action of $F$ and $G$ on the generators.
 Proposition \ref{prop:GK-groupoid} is proved.

 {\small

M.K.:  Kavli IPMU, University of Tokyo, 5-1-5 Kashiwanoha, Kashiwa, Chiba,
277-8583 Japan,  {\tt mikhail.kapranov@yale.edu }

V.S.: Institut de Math\'ematiques de Toulouse, Universit\'e Paul Sabatier, 118 route de Narbonne, 
31062 Toulouse, France, 
 {\tt schechtman@math.ups-tlse.fr }

}

\end{document}